\documentclass[a4paper,11pt]{article}%
\usepackage{amsmath,amssymb,amsthm,latexsym}
\usepackage{hyperref}
\usepackage{epsfig,subfigure}
\usepackage{fullpage}
\usepackage{amsmath}
\usepackage{amsfonts}
\usepackage{amssymb,bbold}
\usepackage{graphicx}%
\usepackage{tikz}
\usepackage{array}
\usepackage{afterpage}

\providecommand{\U}[1]{\protect\rule{.1in}{.1in}}
\theoremstyle{plain}
\newtheorem{thm}{Theorem}[section]

\newtheorem{remark}[thm]{Remark}

\newcommand{\nn}{\nonumber}
\newcommand{\R}{{\mathbb R}}

\newcommand{\N}{{\mathbb N}}

\newcommand{\dx}{ \, {\rm d} x}
\newcommand{\dy}{ \, {\rm d} y}

\newcommand{\Om}{\Omega}
\newcommand{\barrho}{\bar{\rho}}
\newcommand{\trho}{\tilde{\rho}}
\newcommand{\rhooneout}{\rho_1^{\text{out}}}
\newcommand{\rhotwoout}{\rho_2^{\text{out}}}

\newcommand{\as}{a_s}
\newcommand{\bs}{b_s}
\newcommand{\ac}{a_c}
\newcommand{\bc}{b_c}

\newcommand{\A}{A}
\newcommand{\B}{B}
\newcommand{\M}{M}
\newcommand{\D}{D}
\newcommand{\f}{f}

\newcommand{\Rtwo}{R_2}
\newcommand{\Rone}{R_1}
\newcommand{\Rzero}{R_0}
\newcommand{\Rthree}{R_0}

\newcommand{\lm}{\lambda_m}

\newcommand{\CalA}{{\mathcal{A}}}
\newcommand{\CalB}{{\mathcal{B}}}
\newcommand{\CalF}{{\mathcal{F}}}

\newcommand{\Ord}{\mathcal{O}}

\newcommand{\C}{\mathbb{C}}
\newcommand{\Z}{\mathbb{Z}}

\newcommand{\supp}{\operatorname*{supp}}
\newcommand{\sgn}{\operatorname*{sgn}}
\newcommand{\eps}{\varepsilon}
\newcommand{\leqs}{\leqslant}
\newcommand{\geqs}{\geqslant}

\begin{document}

\title{Equilibria for an aggregation model with two species}
\author{Joep H.M.~Evers \and Razvan C.~Fetecau \and Theodore Kolokolnikov}
\maketitle

\begin{center}
\textbf{Abstract}
\end{center}

\begin{quote}
We consider an aggregation model for two interacting species. The coupling between the species is via their velocities, that incorporate self- and cross-interactions. Our main interest is categorizing the possible steady states of the considered model. Notably, we identify their regions of existence and stability in the parameter space. For assessing the stability we use a combination of variational tools (based on the gradient flow formulation of the model and the associated energy), and linear stability analysis (perturbing the boundaries of the species' supports). We rely on numerical investigations for those steady states that are not analytically tractable. Finally we perform a two-scale expansion to characterize the steady state in the limit of asymptotically weak cross-interactions.
\end{quote}



\section{Introduction}
\label{sect:intro}

In this paper we consider a two-species aggregation model in the form of a system of partial differential equations:
\begin{subequations}
\label{eqn:model}
\begin{gather}
\frac{\partial \rho_{1}}{\partial t}+\nabla\cdot(\rho_1 v_1)=0, \quad v_1=-\nabla K_s\ast\rho_1 - \nabla K_c \ast \rho_2, \label{eqn:v1} \\
\frac{\partial \rho_{2}}{\partial t}+\nabla\cdot(\rho_2 v_2)=0, \quad v_2=-\nabla K_c\ast\rho_1 - \nabla K_s \ast \rho_2. \label{eqn:v2}
\end{gather}
\end{subequations}
Here, $\rho_1$ and $\rho_2$ are the densities of the two species, $K_s$ and $K_c$ are self- and cross-interaction potentials, and the asterisk $\ast$ denotes convolution. The self-interaction potential $K_s$ models inter-individual social interactions within the same species, while $K_c$ models interactions between individuals of different species. Typically, interaction potentials are assumed to be symmetric, and also, to model long-range attraction and short-range repulsion. Model \eqref{eqn:model} applies to arbitrary spatial dimension.

Our main motivation for studying this model is the self-organization occurring in aggregates of biological cells. Experimentally observed sorting of embryonic cells has been documented in various works \cite{ArmPaiShe2006, FotPflForSte1996}; we also refer to Figure 1 in \cite{Brodland2002} for a schematic overview of possible patterns of two species of certain embryonic cells. These patterns range from complete mixing at the cell level, to full separation of the two species. The major goal of the present paper is to investigate equilibrium solutions, along with their stability, using model \eqref{eqn:model} in two spatial dimensions. The nonlocality in our model resembles that in e.g.~\cite{ArmPaiShe2006, PaiBloSheGer2015}.

The one-species analogue of \eqref{eqn:model} is a mathematical model for collective behaviour that has seen a surge of attention in recent literature. A variety of issues have been studied for the one-species model, such as the  well-posedness of solutions \cite{BodnarVelasquez2, BertozziLaurentRosado,Figalli_etal2011,EvHiMu2016}, equilibria and long-time behaviour \cite{LeToBe2009,FeHuKo11,Brecht_etal2011,FeHu13},  blow-up (in finite or infinite time) by mass concentration~\cite{FeRa10,BertozziCarilloLaurent, HuBe2010},  existence and characterization of global minimizers for the associated interaction energy \cite{Balague_etalARMA,ChFeTo2015,CaCaPa2015,SiSlTo2015}, and passage from discrete to continuum by mean-field limits \cite{CarrilloChoiHauray2014}. A particularly appealing aspect of the model is that despite its simplicity, its solutions can exhibit complex behaviour and can capture a wide variety of ``swarm" behaviour. Provoking and motivational galleries of solutions that can be obtained with the one species model can be found for instance in \cite{KoSuUmBe2011,Brecht_etal2011}. They include aggregations on disks, annuli, rings, soccer balls, and many others.

Despite the intensive activity on the one-species model, its extensions to multiple species have remained largely unexplored. From an analysis viewpoint, the well-posedness of solutions to multi-species aggregation models of type \eqref{eqn:model} has been recently considered in various works \cite{DiFrancescoFagioli2013, CrippaLecureux2013}. The general setup is to investigate existence and uniqueness of weak measure solutions using tools from optimal mass transportation such as Wasserstein distance(s). In particular, the authors of  \cite{DiFrancescoFagioli2013} consider a more general model where the two species have distinct self-interaction potentials, and cross-interaction potentials that are a scalar multiple of each other. They show, under certain assumptions on the potentials, that the two-species aggregation model represents a gradient flow with respect to a modified Wasserstein distance of an interaction energy that comprises self- and cross-interaction terms. 

Two-species models similar to \eqref{eqn:model} have been studied recently in the context of predator-prey dynamics  \cite{ChKo2014, DiFrancescoFagioli2016}. To model such interactions one needs to consider cross-interaction potentials that have opposite signs, so that the predator is attracted to prey, while the prey is repelled by it. Both \cite{ChKo2014} and 
\cite{DiFrancescoFagioli2016} show very intricate patterns that form dynamically (and at equilibrium) with such predator-prey systems. Note that in model \eqref{eqn:model}, the cross-interaction potential $K_c$ is not restricted to be exclusively attractive or repulsive; the particular form considered in the sequel includes in fact both attractive and repulsive cross-interactions between the two species.
 
To investigate the equilibria and their stability for model \eqref{eqn:model}, we present in the current paper two approaches, each with their own merits and limitations. One approach is to consider the variational interpretation \cite{DiFrancescoFagioli2013}, and investigate equilibrium solutions to \eqref{eqn:model} as stationary points of the interaction energy. To establish whether such stationary points are energy minimizers we adopt and extend the framework developed in \cite{BeTo2011} for the one-species analogue of model \eqref{eqn:model}. In such setup, the problem reduces to investigating the first and second variations of the energy for various perturbations that may occur. 

The other approach considered in this work is quite different in spirit, and it consists in a linear stability analysis of the boundaries that enclose the two aggregations. Here, we take advantage of the choice of potentials considered in this paper, consisting of Newtonian repulsion and quadratic attraction, for which the equilibria are made of compactly supported aggregations of constant densities \cite{FeHuKo11}.  While this approach to stability is novel, it relates to previous works of one of the authors of the current paper. In \cite{ChKo2014}, a similar technique is used to perturb the boundaries of a compactly supported swarm density. In \cite{KolKevCar2014}, the authors perturb a finite number of discrete points on a circle. The algebraic description can however be generalized to a (continuous) circular curve, and the formula can moreover be rewritten exactly in the form we use (including the Fourier decomposition; cf.~\eqref{eqn:perturb}). We apply this procedure to a particular equilibrium that resembles the image of a target used for shooting or archery (see e.g.~Figure \ref{fig:target}); henceforth, we simply call such configuration a ``target". We find excellent agreement between the stability analysis and the numerical simulations. 

\afterpage{
\begin{figure}[t!]
\centering
\begin{tikzpicture}[>= latex]
\begin{scope}[xscale=2.5,yscale=1.6]
\pgfmathsetmacro\M {2};

\draw[domain=0:1,smooth,variable=\b,line width=2]  plot ({(1+\M*\b)/(\b+\M)},{\b});
\draw[domain=1:4,smooth,variable=\b,line width=1, dashed]  plot ({(1+\M*\b)/(\b+\M)},{\b});
\draw (0.42,1.5) node {\includegraphics[height=.13\textwidth]{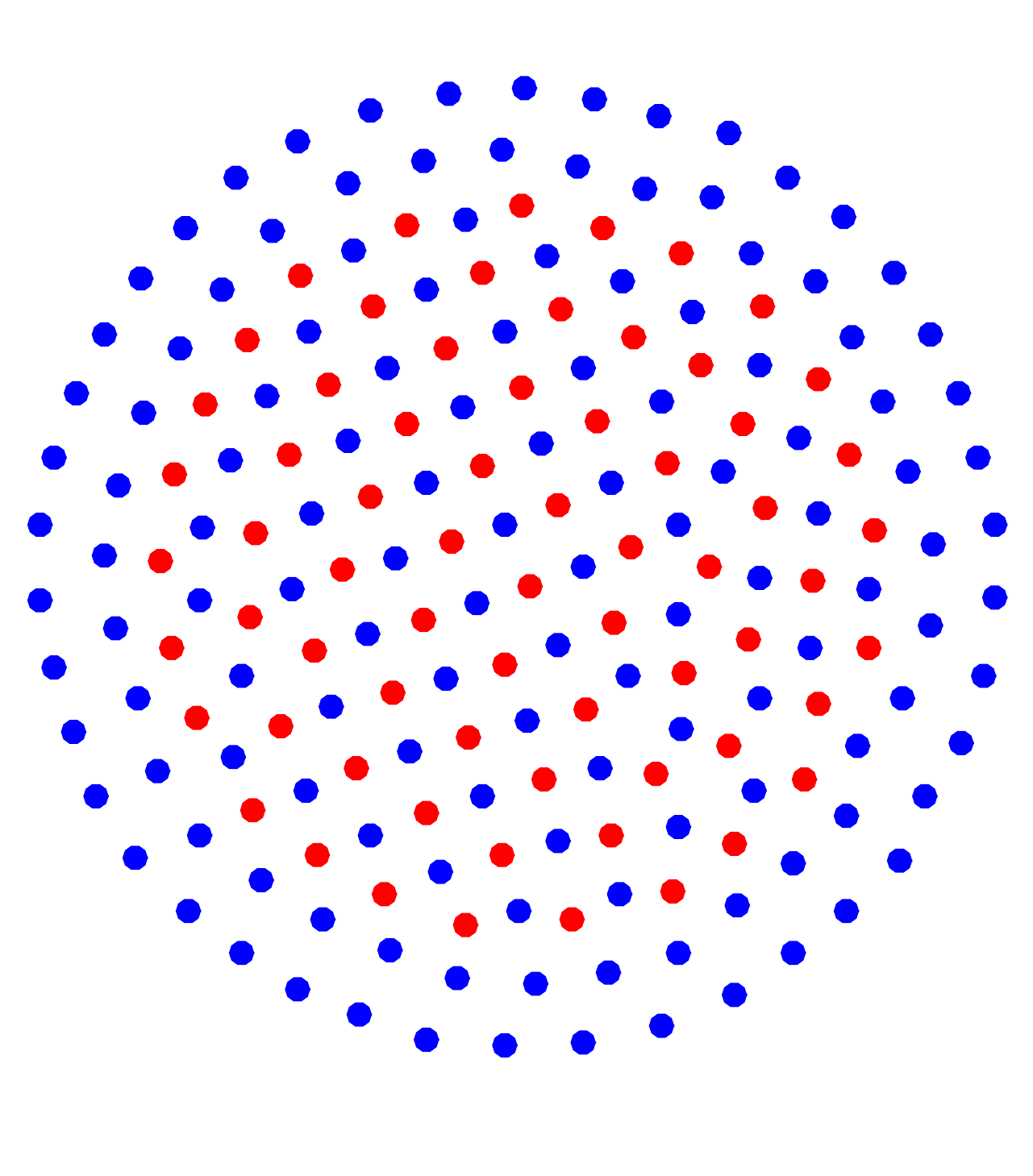}};
\draw (1.5,3.5) node {\includegraphics[height=.13\textwidth]{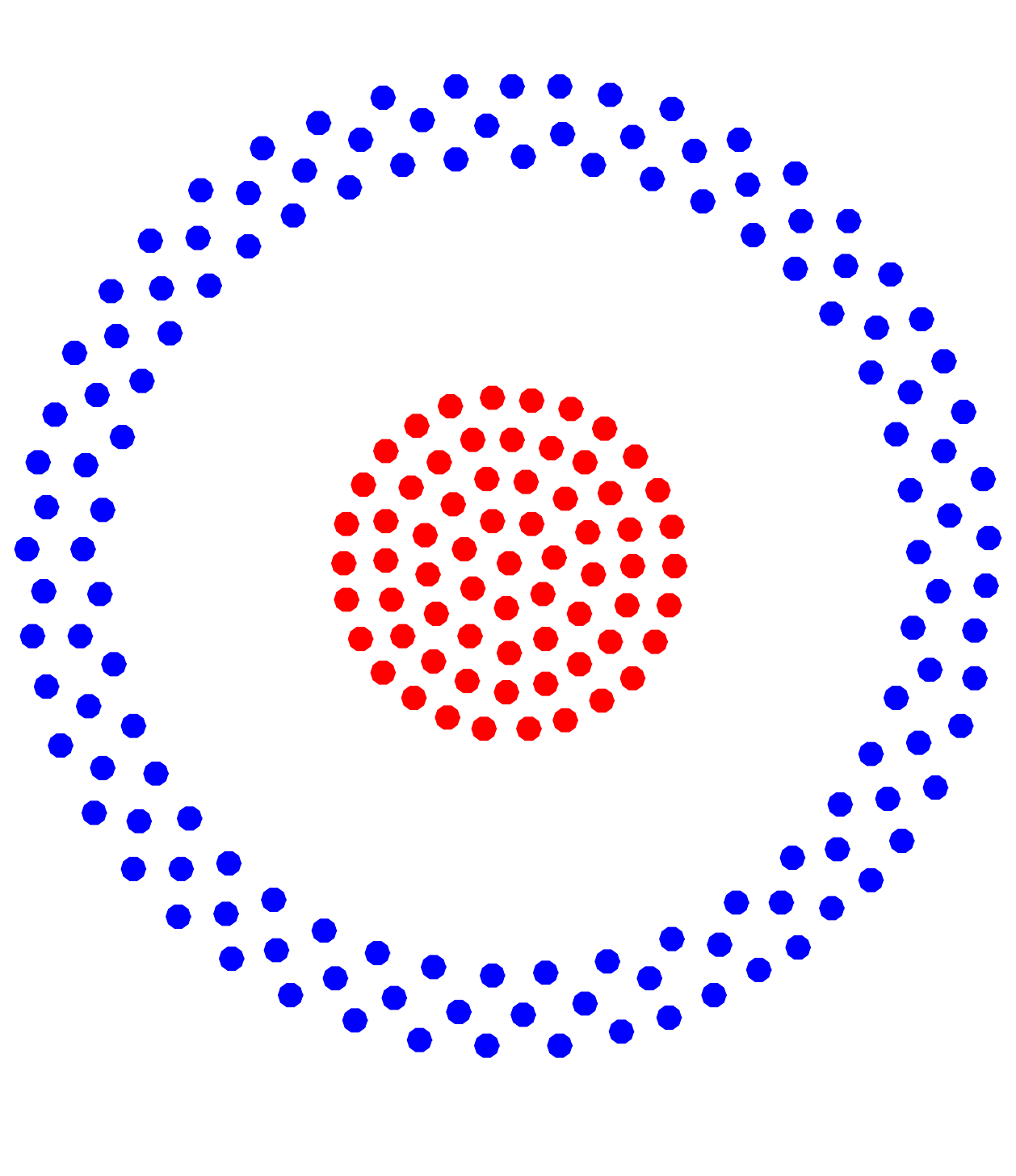}};
\draw (2.5,3.5) node {\includegraphics[height=.13\textwidth]{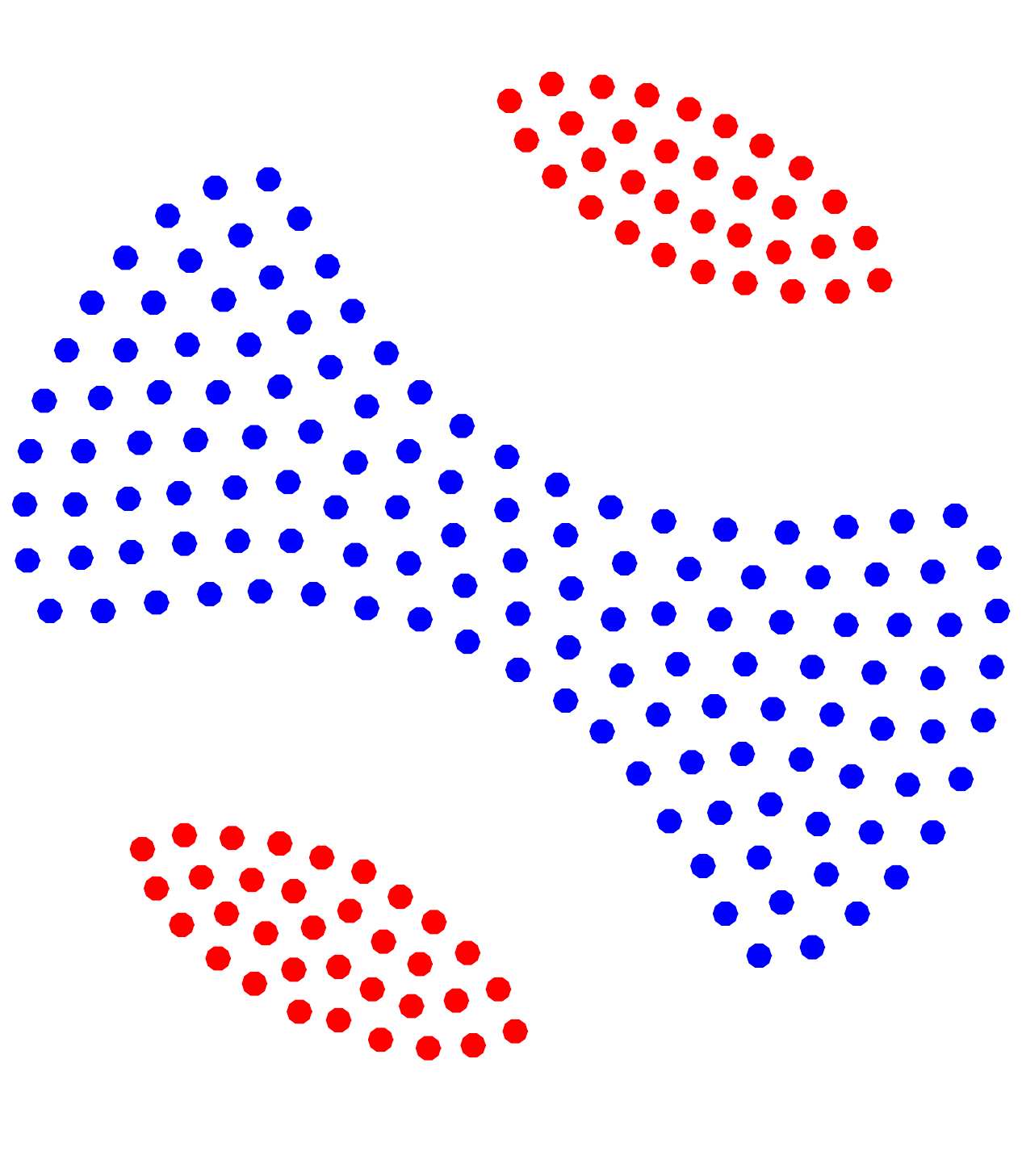}};
\draw (2.6,1.8) node {\includegraphics[height=.13\textwidth]{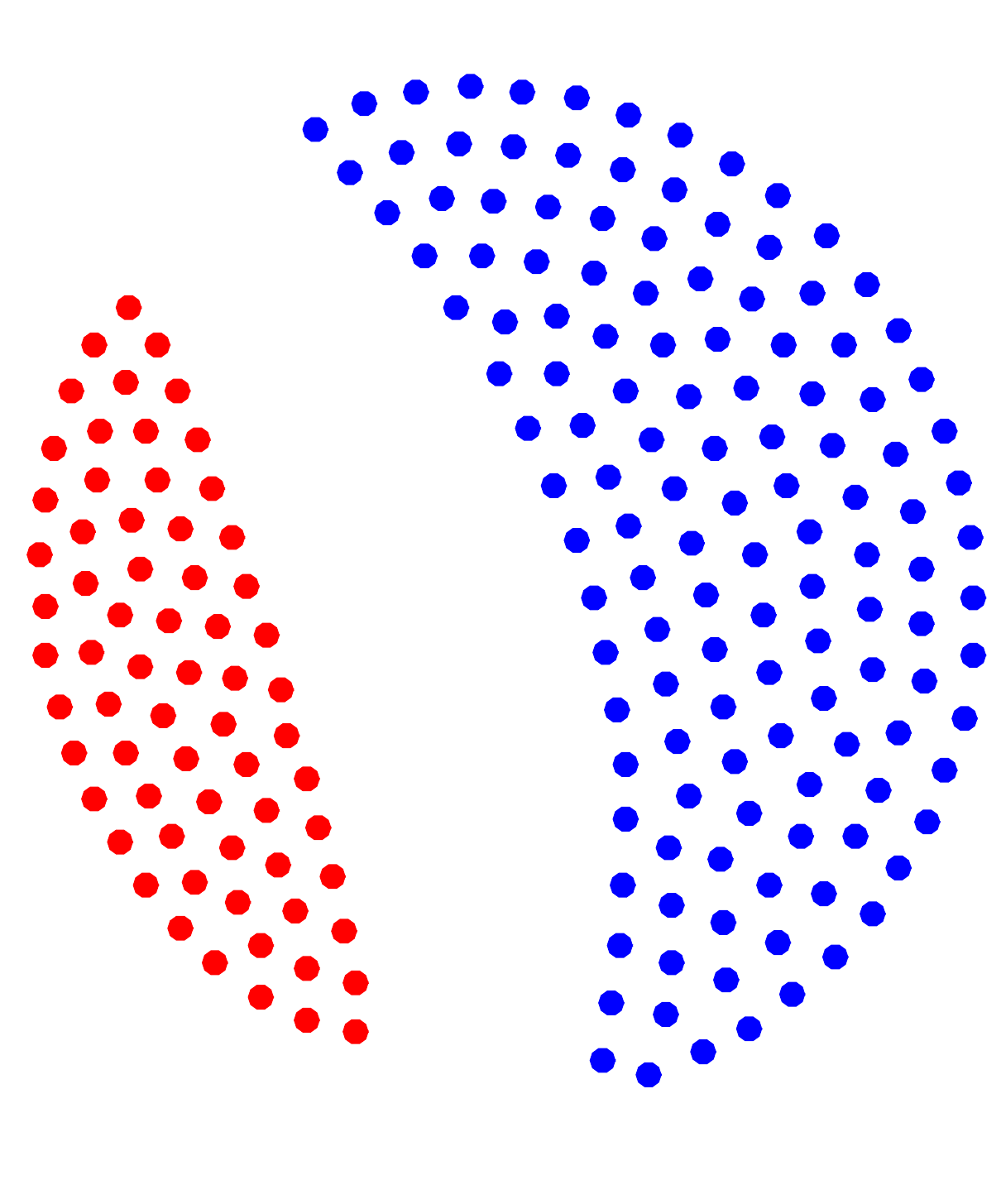}};
\draw (3.6,1.8) node {\includegraphics[height=.14\textwidth]{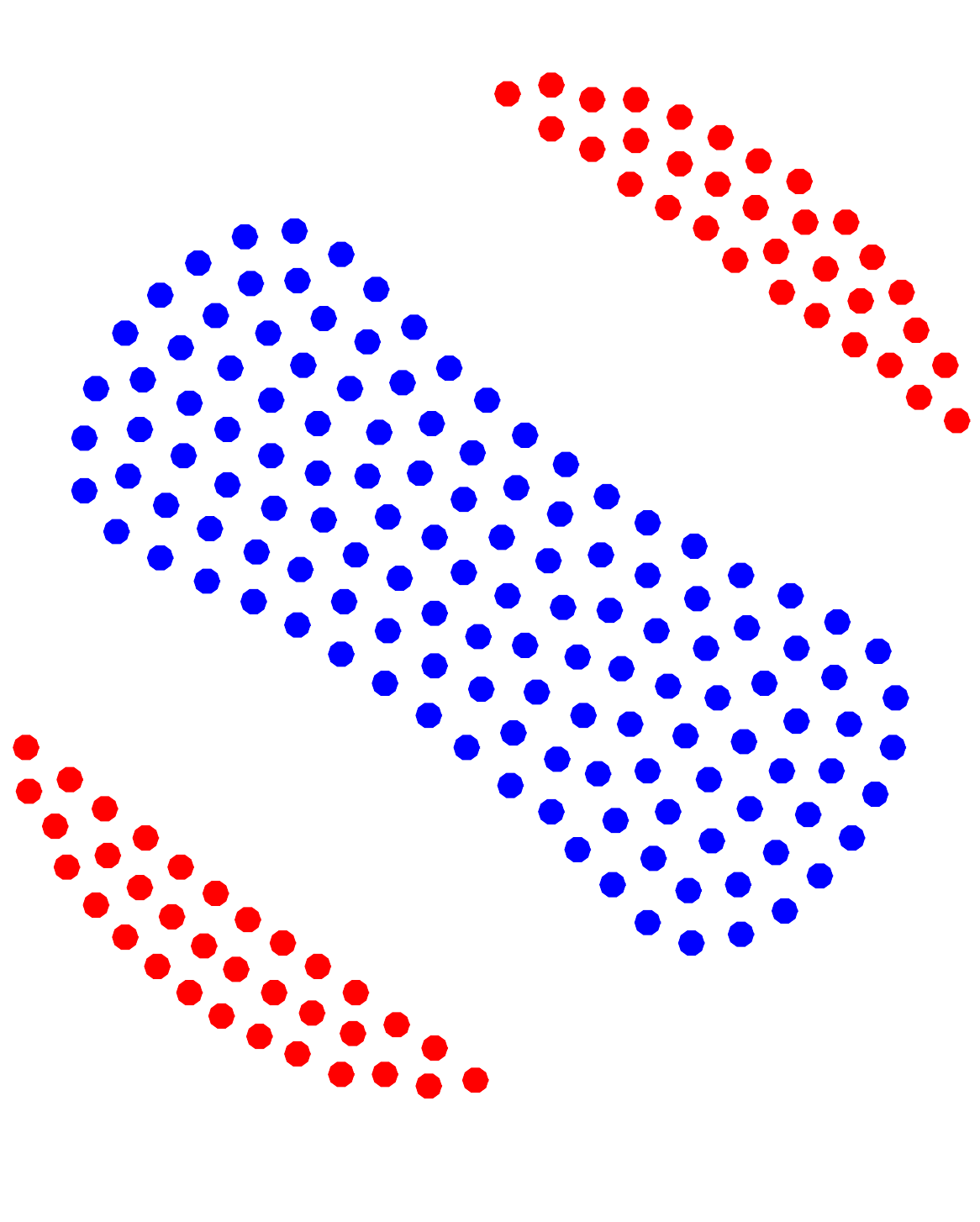}};
\draw (3.1,0.5) node {\includegraphics[height=.11\textwidth]{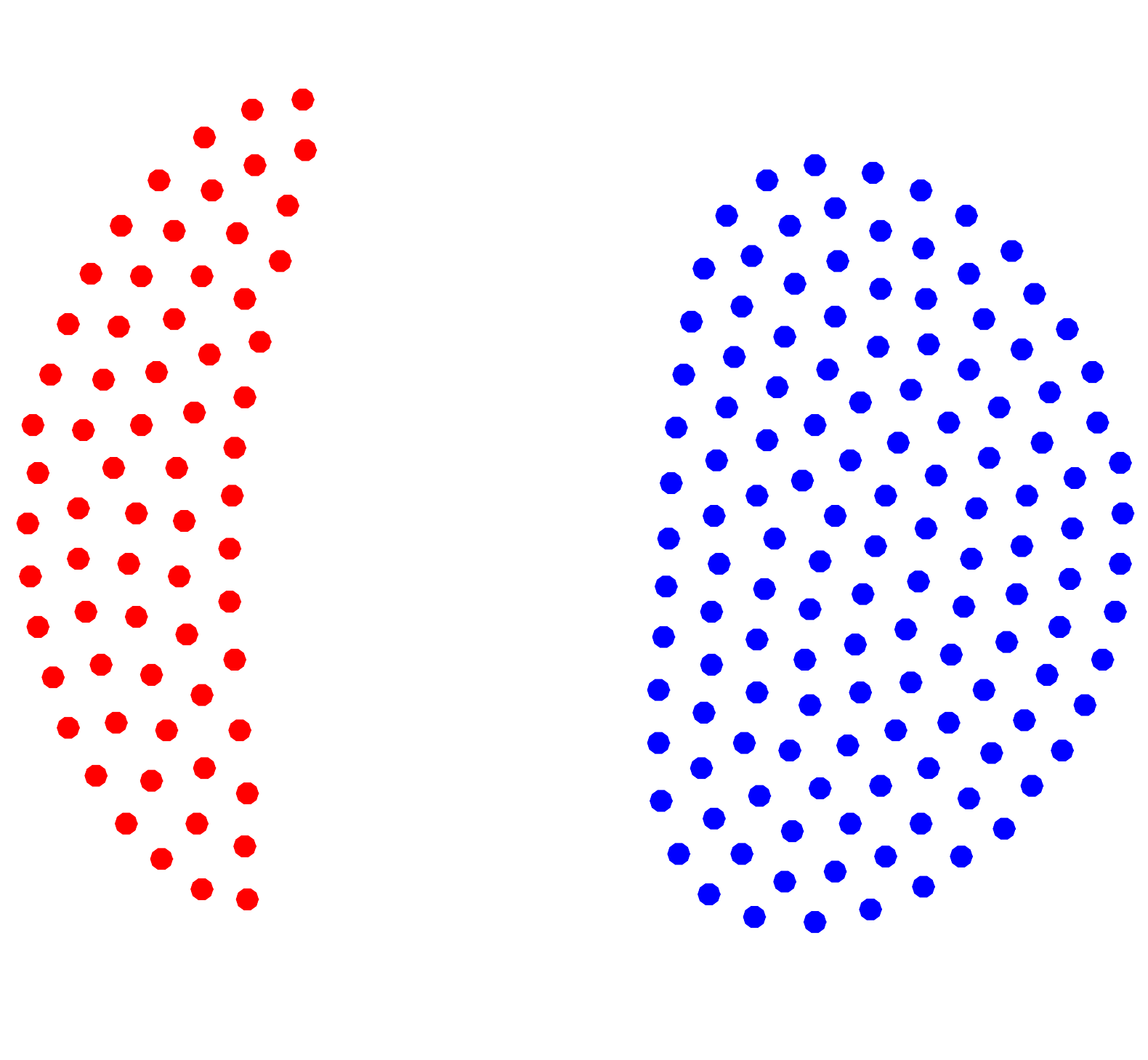}};
\draw (-0.8,0.75) node {\includegraphics[height=.14\textwidth]{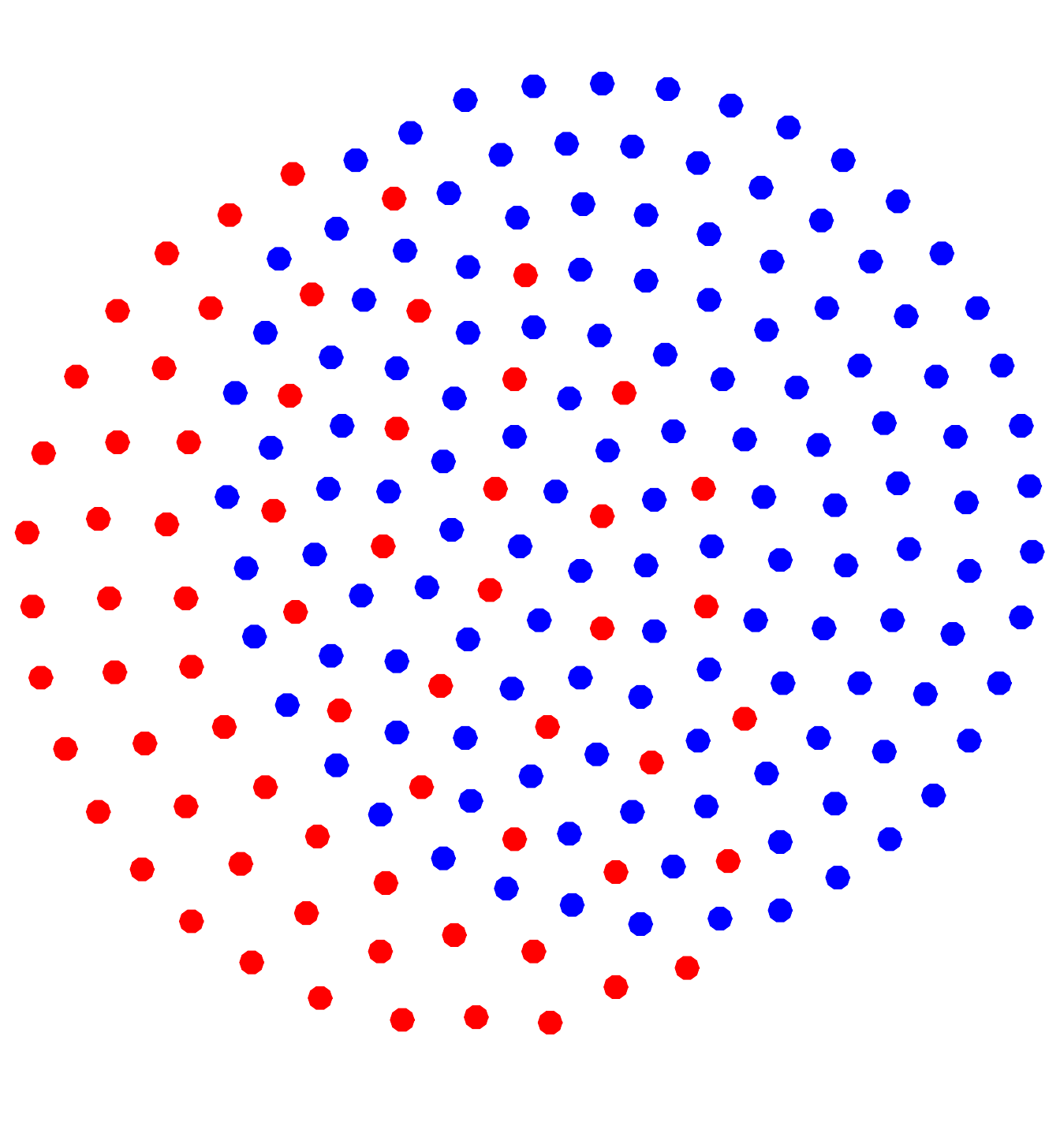}};
\draw [line width=0.75,->,color=gray!90] (-0.8,0.15)  to[out=-45,in=235] (0.4,0.2);

\draw[->,line width=2] (0,0) -- (0,4) node[left] {$B=\dfrac{b_c}{b_s}$};
\draw[->,line width=2] (0,0) -- (4,0) node[below] {$A=\dfrac{a_c}{a_s}$};

\draw[line width=0.25, dashed] (1,1) -- (4,1);
\draw[line width=1] (\M,0.04) -- (\M,-0.04) node[below] {$M$};
\draw[line width=1] (1/\M,0.04) -- (1/\M,-0.04) node[below] {$1/M$};
\draw[line width=1] (1,0.04) -- (1,-0.04) node[below] {$1$};
\draw[line width=1] (0.04,1) -- (-0.04,1) node[left] {$1$};

\draw[domain=0:1,smooth,variable=\b,line width=1,dashed]  plot ({(\b+\M)/(1+\M*\b)},{\b});
\draw[domain=1:4,smooth,variable=\b,line width=2]  plot ({(\b+\M)/(1+\M*\b)},{\b});
\draw[line width=2] (0,0) -- (4,4);
\begin{scope} 
\tikzstyle{every node} = [draw,rectangle, fill=gray!5, line width = 0.5,rounded corners=2pt]
\draw (-0.8,1.7) node {$\substack{A=0.5\\B=0.4}$} ;
\draw (0.42,2.4) node {$\substack{A=0.5\\B=1}$} ;
\draw (2,2.8) node {$\substack{A=3\\B=3.5}$};
\draw (3.1,2.5) node {$\substack{A=3\\B=2}$};
\draw (3.85,0.5) node {$\substack{A=3\\B=0.75}$};
\end{scope}
\end{scope}
\end{tikzpicture}
\caption{Overview of steady states of model \eqref{eqn:model} encountered in numerical investigations of the associated particle system -- cf.~Section \ref{sec:particle model and num}. The equilibria are placed in the $(A,B)$-plane according to their parameter values. 
For $(A,B)=(3,3.5)$ and for $(A,B)=(3,2)$ two distinct steady states are observed, depending on the initial data. Theoretical considerations using a variational approach and a linear analysis focus on the ``overlap solution" seen here for $(A,B)=(0.5,1)$ and the ``target equilibrium" illustrated in the left-hand plot for $(A,B)=(3,3.5)$.}
\label{fig:num st st}
\end{figure}
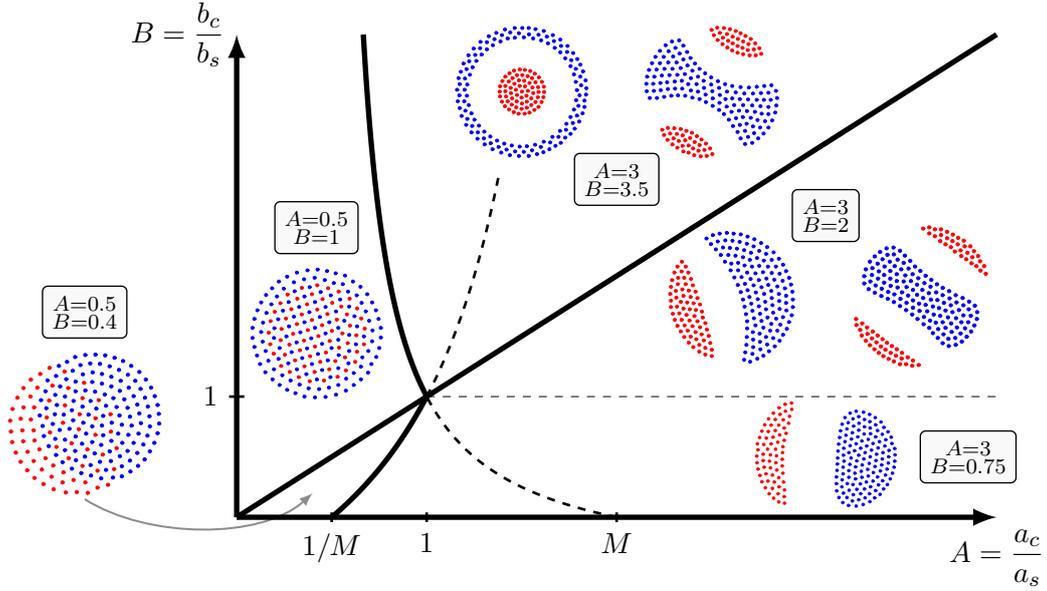


\begin{figure}[h!]
\centering
\begin{tikzpicture}[>= latex]
\begin{scope}[xscale=2.5,yscale=1.6]
\pgfmathsetmacro\M {2};

\draw[fill=gray!15,dashed,line width=1] (0,1) rectangle (1,4);

\draw[->,line width=2] (0,0) -- (0,3.5) node[left] {$B=\dfrac{b_c}{b_s}$};
\draw[->,line width=2] (0,0) -- (3.5,0) node[below] {$A=\dfrac{a_c}{a_s}$};

\draw[line width=0.25, dashed] (1/\M,4) -- (1/\M,-0.04);
\draw[line width=0.25, dashed] (\M,3.5) -- (\M,-0.04);
\draw[line width=1] (\M,0.04) -- (\M,-0.04) node[below] {$M$};
\draw[line width=1] (1/\M,0.04) -- (1/\M,-0.04) node[below] {$1/M$};
\draw[line width=1] (1,0.04) -- (1,-0.04) node[below] {$1$};
\draw[line width=1] (0.04,1) -- (-0.04,1) node[left] {$1$};

\draw[domain=0:3.5,smooth,variable=\b,line width=2]  plot ({(1+\M*\b)/(\b+\M)},{\b});
\draw[domain=0:3.5,smooth,variable=\b,line width=2]  plot ({(\b+\M)/(1+\M*\b)},{\b});

\begin{scope} 
\tikzstyle{every node} = [draw,rectangle, fill=gray!5, line width = 0.5]
\draw [line width=0.75,->] (1/\M,3.75) node[above] {$B=\dfrac{A-M}{1-M\,A}$} arc (180:220:0.6);
\draw [line width=0.75,->] (\M-0.3,3.75) node[above] {$B=\dfrac{1-M\,A}{A-M}$} arc (0:-43:0.6);
\draw[line width=2] (0,0) -- (3.5,3.5);
\draw [line width=0.75,->] (3.15,3.75) node[above] {$B=A$} arc (180:220:0.6);
\end{scope}
\begin{scope} 
\tikzstyle{every node} = [draw,circle, fill=gray!5, line width = 0.25]
\draw [line width=0.75,->] (0,-0.5) node {$D_1$} -- (0.35,0.15);
\draw (1,0.35) node {$D_2$};
\draw (2.1,1) node {$D_3$};
\draw (2.1,3) node {$D_4$};
\draw (1.05,2.6) node {$D_5$};
\draw (0.4,1.18) node {$D_6$};
\end{scope}

\end{scope}
\end{tikzpicture}
\caption{Regions of existence and stability of the equilibria. The grey shaded region is the one where $A<1$ and $B>1$; here we manage to identify \emph{global} minimizers.}
\label{fig:stability}
 \end{figure}
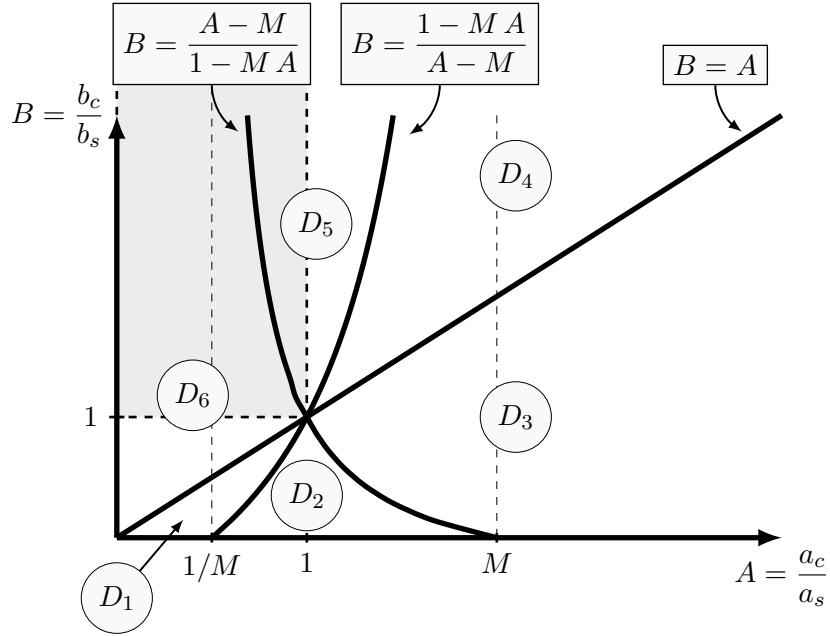
}

To further motivate our study, we introduce here Figure \ref{fig:num st st}, which shows a gallery of equilibria that can be obtained with model \eqref{eqn:model}. Numerically, we employ the particle system model associated to the PDE model \eqref{eqn:model}; see Section \ref{sec:particle model and num} for more details. The heavier species of mass $M_1$ is shown in blue and the lighter species of mass $M_2$ in red. The subindices $s$ and $c$ refer to self- and cross-interactions, respectively. The coefficients $\as,\ac$ and $\bs,\bc$ represent the strengths of the corresponding repulsive and attractive interactions respectively (for the specific form of the interaction kernels, see \eqref{eqn:intpot}). Note that in the figure we employ the notations:
\begin{equation}
\label{eqn:ABM}
\A = \frac{\ac}{\as}, \quad \B = \frac{\bc}{\bs}, \quad \M = \frac{M_1}{M_2},
\end{equation}
where by convention, $\M>1$.

In the figure, the obtained steady states appear in the $(A,B)$-phase plane according to the concerning parameter values. The solid and dashed lines and curves in the diagram are significant for the existence and stability of steady states. For their exact descriptions, see also Figure \ref{fig:stability}.

\paragraph{Overview of the results of the paper.}
Our theoretical approach focuses first on a state in which the two species segregate and we have a disk of one species surrounded concentrically by an annular region of the other species. As noted above, we call such states targets. One such equilibrium is shown for $(A,B)=(3,3.5)$, at the left-hand side, in Figure \ref{fig:num st st}. 

In Section \ref{sect:target} we show that a target with the \emph{lighter} species inside exists and it is stable in the parameter region $D_4\cup D_5$ (for notation, see Figure \ref{fig:stability}), which is exactly the region where we observe this target numerically as a steady state. However, a target with the \emph{heavier} species inside is, whenever it exists at all, not stable and hence it is not observed numerically. These theoretical conclusions follow from the linear stability analysis. In the variational approach, for most of the parameter space, we had to restrict ourselves to a subset of all possible perturbations (specifically, to perturbations which are not exclusively supported in the support of the equilibrium state;  see Section \ref{subsect:prelim-var} for details). General perturbations could only be dealt with for $A<1$ and $B>1$. Due to such limitations, for both types of targets (light and heavy inside, respectively), the variational approach in Section \ref{sect:target} leads to stability regions that are larger compared to what the linear stability analysis and the numerics predict. Nevertheless, an important result of the variational approach (see Remark \ref{rmk:target-gmin}), which was not possible to be derived by the other method,  is that the target equilibrium is a {\em global} minimizer for parameters $(A,B)$ in the subset of $D_5$ where $A<1$ (see shaded area in Figure \ref{fig:stability}).

Next, we investigate in Section \ref{sect:overlap} the state in which the supports of both species are two concentric disks of possibly different widths. Within the smallest of these disks, the two species coexist. We use the term ``overlap" equilibria to refer to such configurations; an example of such a state is shown in Figure \ref{fig:num st st} for $(A,B)=(0.5,1)$. In Section \ref{sect:overlap} only the variational approach is employed, as due to overlapping supports, perturbing boundaries to perform a linear stability, is very delicate. Similar to the target, while for parameters with $A<1$ and $B>1$ we were able to study general perturbations, for the rest of the parameter space only certain perturbations have been considered (details in Section \ref{subsect:prelim-var}). We show that the overlap state with the \emph{lighter} species inside (the one shown in Figure \ref{fig:num st st}), is a local minimizer of the energy for parameters $(A,B)$ in region $D_6$. This is in agreement with what we find numerically. A major result is that this overlap equilibrium  is a {\em global} minimizer in $D_6$ where $B>1$ (see shaded area in Figure \ref{fig:stability} and also, Remark \ref{rmk:overlap-gmin}). The overlap solution with the \emph{heavier} species inside is never observed numerically. Still it is a minimizer of the energy when $(A,B) \in \D_1$ for the specific class of perturbations considered. It is our conjecture in Section \ref{subsect:hinside} that this overlap state is \emph{not} in fact a minimizer for general perturbations.

In Section \ref{sec:num unstable modes}, we show by means of numerics how the non-radially symmetric states in the right-hand part of Figure \ref{fig:num st st}, relate to the target solution. We take parameters in certain regions where the target (either light, or heavy inside) is not a stable steady state. We initialize the system at the target (some small numerical deviation from the actual equilibrium being inevitably present), run the dynamics and wait until the small initial error triggers the instability. The whole process is shown in Figures \ref{fig:unstab modes targStab} and \ref{fig:unstab modes targUnstab}. It turns out that all non-symmetric states in Figure \ref{sec:num unstable modes} can be explained in this way, and that moreover we can identify distinct types (\emph{modes}) of instability.

Similarly, the non-symmetric state shown for $(A,B)=(0.5,0.4)$ arises by initializing the system in the overlap state. This is shown in Figure \ref{fig:unstab modes overlap} of Section \ref{sec:instab overlap D6 D1}. Apparently, the two species' centres of mass coincide in $D_6$, yet the symmetry is broken once we cross the boundary $B=A$ and take parameters in $D_1$.
The steady state for $(A,B)=(0.5,0.4)$ exhibits in the middle an 
area in which the two species coexist. In Section \ref{sec:weak cross} we quantify, as a function of $A$ and $B$, the size of this part of the support, and the shift of the centres of mass. We perform this analysis asymptotically in the limit of weak cross-interactions; that is, asymptotically close to $(A,B)=(0,0)$.

\section{Preliminaries}
\label{sect:prelim}

In this section we present the two approaches taken in this paper to study the stability of equilibria, i.e., variational and linear analysis, along with some general properties of model \eqref{eqn:model}. We refer specifically to the model in two dimensions, but the general properties and the variational formalism applies to any spatial dimension.

\subsection{General preliminaries}
It is common for aggregation models of type \eqref{eqn:model} to assume symmetry of the inter-individual interactions, or equivalently, the symmetry of the interaction potentials \cite{DiFrancescoFagioli2013}:
\begin{equation}
\label{eqn:Ksym}
K_s(x) = K_s(-x), \qquad K_c(x) = K_c(-x), \qquad \text{ for all } x.
\end{equation}
In the present work we only consider self- and cross-interaction potentials that satisfy \eqref{eqn:Ksym}.

\paragraph{Conservation of mass and centre of mass.} The dynamics of model \eqref{eqn:model} conserve the mass of each species:
 \begin{equation}
\label{eqn:massM}
\int \rho_i(x,t) \dx = M_i \qquad \text{ for all } t\geq 0, \quad i=1,2,
\end{equation}
as well as the total centre of mass:
\begin{equation}
\label{eqn:cmass}
\int x(\rho_1(x,t) + \rho_2(x,t)) \dx = \text{ const. } \qquad \text{ for all } t\geq 0.
\end{equation}
Indeed, conservation of mass follows immediately as the equations of motion \eqref{eqn:model} are in conservation law form. To get \eqref{eqn:cmass}, multiply each equation in \eqref{eqn:model} by $x_k$ ($x_k$ denotes the $k$-th spatial coordinate), and add the two equations. After integrating by parts and using vanishing at infinity boundary conditions, one finds
\begin{equation}
\label{eqn:dcmdt}
\frac{d}{dt} \int x_k (\rho_1(x,t)+\rho_2(x,t))\, dx = \int (\rho_1 v_1 + \rho_2 v_2) \cdot \mathbf{e}_k \, dx,
\end{equation}
where $\mathbf{e}_k$ denotes the unit vector in the direction of the $k$-th coordinate.  Furthermore, by the expressions of $v_1$ and $v_2$ from \eqref{eqn:v1} and \eqref{eqn:v2}, together with the symmetry of $K_s$, one gets
\begin{align*}
\int (\rho_1 v_1 + \rho_2 v_2) \, dx &= - \underbrace{\iint \nabla K_s(x-y) \rho_1(x) \rho_1(y) \, dx dy}_{=0 \text { by \eqref{eqn:Ksym}} } - \iint \nabla K_c(x-y) \rho_1(x) \rho_2(y) \, dx dy\\
& \quad - \underbrace{\iint \nabla K_s(x-y) \rho_2(x) \rho_2(y) \, dx dy}_{=0 \text { by \eqref{eqn:Ksym}} } -  \iint \nabla K_c(x-y) \rho_2(x) \rho_1(y) \, dx dy.
\end{align*}
Finally, the symmetry of $K_c$ yields the expression above zero, as
\[
\iint \nabla K_c(x-y) (\rho_1(x) \rho_2(y) + \rho_2(x) \rho_1(y)) \, dx dy = \iint (\underbrace{\nabla K_c(x-y) + \nabla K_c(y-x)}_{=0 \text{ by } \eqref{eqn:Ksym}}) \rho_1(x) \rho_2(y) \, dx dy.
\]
From \eqref{eqn:dcmdt} one can now derive the conservation of centre of mass \eqref{eqn:cmass}.


\paragraph{Interaction potentials, equilibria and phase plane.}
The present study focuses on the aggregation model \eqref{eqn:model} in two dimensions, and specific self- and cross-interaction potentials $K_s$ and $K_c$, given by 
\begin{subequations}
\label{eqn:intpot}
\begin{align}
\label{eqn:selfpot}
K_s(x) &= -\as \ln|x| + \frac{\bs}{2} |x|^{2}, \\
\label{eqn:crosspot}
K_c(x) &= -\ac \ln|x| + \frac{\bc}{2} |x|^{2}.
\end{align}
\end{subequations}
Here, the coefficients $\as, \ac, \bs$, and  $\bc$ set the magnitudes of the self- and cross-repulsion and -attraction, respectively.

Interaction potentials in the form \eqref{eqn:intpot}, consisting of Newtonian repulsion and quadratic attraction, have been studied in various works on the one-species model \cite{FeRa10, FeHuKo11, FeHu13, HuFe2013}.  It has been demonstrated in \cite{BertozziLaurentLeger, FeHuKo11} that solutions to the one-species aggregation equation, with an interaction potential of form  \eqref{eqn:intpot}, approach asymptotically a radially symmetric equilibrium that consists in a ball of constant density. 

For the two-species model studied in this paper, such potentials also have the remarkable property that they lead to equilibria of constant densities.  Indeed, from the equation for $v_1$ in \eqref{eqn:v1}, and \eqref{eqn:intpot}, we get
\begin{align}
\nabla\cdot v_1 &= - \Delta K_s \ast \rho_1 - \Delta K_c \ast \rho_2 \nonumber \\
&= 2 \pi \as \rho_1 + 2 \pi \ac \rho_2 - 2 \bs M_1 - 2 \bc M_2,
\label{eqn:divv1}
\end{align}
where for the second equality we also used $\Delta \left( \frac{1}{2 \pi} \ln |x|\right) = \delta_x$, and the mass constraint \eqref{eqn:massM}. Also, by a similar calculation,
\begin{equation}
\label{eqn:divv2}
\nabla\cdot v_2 = 2 \pi \as \rho_2 + 2 \pi \ac \rho_1 - 2 \bs M_2 - 2 \bc M_1.
\end{equation}

At equilibrium, the velocity $v_i$ (and consequently its divergence) must vanish on the support of the respective density $\rho_i$. By \eqref{eqn:divv1} and \eqref{eqn:divv2}, this leaves four possible combinations for the values $(\barrho_1,\barrho_2)$ which the densities of the two species can have at equilibrium:
\begin{equation}
\label{eqn:equilibria}
\begin{aligned}
(0,0), \quad & \left( 0, \frac{\bc M_1 + \bs M_2}{\pi \as} \right), \quad \left( \frac{\bs M_1 + \bc M_2}{\pi \as},0\right), \quad \text{ and } \\ 
& \qquad \left( \frac{(\as \bs - \ac \bc)M_1 + (\as \bc - \ac \bs)M_2}{\pi(\as^2 - \ac^2)}, \frac{(\as \bc - \ac \bs)M_1 + (\as \bs - \ac \bc)M_2}{\pi(\as^2 - \ac^2)} \right).
\end{aligned}
\end{equation}
The first pair corresponds to points outside the supports of both $\barrho_1$ and $\barrho_2$. The second and the third represent equilibrium densities at points that lie in the support of one density, but outside the support of the other. Finally, the fourth pair represents densities in regions where the two species overlap.

We use the parameters $A=a_s/a_c$ and $B=b_s/b_c$ introduced in \eqref{eqn:ABM}, to define the $(A,B)$ phase plane. This phase plane is given in Figure \ref{fig:stability} and subdivided in six regions named $D_1$ to $D_6$. The exact formulas for the boundaries of these regions are a result of our investigation of the existence and stability of steady states. More details will be given in Sections \ref{sect:target} and \ref{sect:overlap}.  Note that the same phase plane was already used in Figure \ref{fig:num st st} to present the gallery of numerical steady states.




\subsection{Variational approach}
\label{subsect:prelim-var}
\paragraph{Energy and gradient flow.} The interaction energy corresponding to model \eqref{eqn:energy} is given by
\begin{equation}
\label{eqn:energy}
E[\rho_1,\rho_2] = \frac{1}{2} \iint K_s(x-y)\left(  \rho_1(x) \rho_1(y) + \rho_2(x) \rho_2(y)\right)  \dx  \dy +   \iint K_c(x-y) \rho_1(x) \rho_2(y) \dx \dy,
\end{equation}
where the two terms represent the self-interaction and the cross-interaction energies, respectively. 
While there has been significant progress recently on the study of minimizers for the interaction energy of the one-species model \cite{Balague_etalARMA,BeTo2011,CaCaPa2015,ChFeTo2015,CaCaPa2015,SiSlTo2015}, there is only a handful of works on the two-species interaction energy. 

In \cite{DiFrancescoFagioli2013} the authors make precise the gradient flow structure (with respect to energy \eqref{eqn:energy}) of model \eqref{eqn:model} by generalizing the theory of gradient flows on probability spaces previously developed for the one-species model \cite{Figalli_etal2011}. The setup there is very general and allows for measure-valued solutions and mildly singular ($C^1$ except at origin) interaction potentials. Studying critical points and minimizers of the energy functional \eqref{eqn:energy} is central to our variational approach. We also point out that there are  related works on ground states for two-phase/two-species interaction energies such as \eqref{eqn:energy}, where interactions are assumed to be purely attractive, but with an additional requirement of boundedness being enforced \cite{Cicalese_etal2016,BerBurPie2016}. In \cite{BerBurPie2016} the setting of \cite{Cicalese_etal2016} is extended by including diffusion (entropy). Moreover, they study the connection between energy minimizers and the long-time dynamics of the gradient flow. 

\paragraph{Equilibria and energy minimizers.} The authors in \cite{BeTo2011} study the energy functional that corresponds to the one species model and find conditions for critical points to be energy minimizers. We adapt the setup from there to the two species model. 

Consider an equilibrium solution $(\barrho_1,\barrho_2)$ with masses $(M_1,M_2)$ and supports $(\Om_1,\Om_2)$, and take a small perturbation $\varepsilon (\trho_1,\trho_2)$:
\begin{equation}\label{eqn:pert dens var}
\rho_i(x) = \barrho_i(x) + \varepsilon \trho_i(x), \qquad i=1,2.
\end{equation}
Given the considerations above,  it is sufficient to consider perturbations that preserve the individual masses of the two species, as well as the total centre of mass. Hence, we have
\begin{equation}
\label{eqn:massc}
\int_{\Om_i} \barrho_i(x) \dx =M_i, \qquad 
\int_{\R^2} \trho_i(x) \dx =0, \qquad i=1,2,
\end{equation}
and
\begin{equation}
\label{eqn:cmassc}
\int_{\R^2} x (\trho_1(x) + \trho_2(x)) \dx =0.
\end{equation}

Since the energy functional is quadratic, one can write:
\begin{equation}
\label{eqn:Equad}
E[\rho_1,\rho_2] = E[\barrho_1,\barrho_2] + \varepsilon E_1[\barrho_1,\barrho_2,\trho_1,\trho_2] + \varepsilon^2 E_2[\trho_1,\trho_2],
\end{equation}
where $E_1$ denotes the first variation:
\begin{equation}
\label{eqn:1stvar}
\begin{aligned}
E_1[\barrho_1,\barrho_2,\trho_1,\trho_2] &=  \int \left[ \int K_s(x-y) \barrho_1(y) \dy + \int K_c(x-y) \barrho_2(y) \dy \right]  \trho_1(x) \dx \\
& \quad + \int \left[ \int K_s(x-y) \barrho_2(y) \dy + \int K_c(x-y) \barrho_1(y) \dy \right]  \trho_2(x) \dx,
\end{aligned}
\end{equation}
and $E_2$ the second variation, which in fact has the same expression as the energy itself:
\begin{equation}
\label{eqn:2ndvar}
E_2[\trho_1,\trho_2] =  E[\trho_1,\trho_2].
\end{equation}
Using the notation 
\begin{subequations}
\label{eqn:Lambda12}
\begin{align}
\Lambda_1(x) &= \int_{\Om_1} K_s(x-y) \barrho_1(y) \dy + \int_{\Om_2} K_c(x-y) \barrho_2(y) \dy, \label{eqn:Lambda1}\\
\Lambda_2(x) &= \int_{\Om_2} K_s(x-y) \barrho_2(y) \dy + \int_{\Om_1} K_c(x-y) \barrho_1(y) \dy, \label{eqn:Lambda2}
\end{align}
\end{subequations}
one can also write the first variation as 
\begin{equation}
\label{eqn:1stvarL}
E_1 [\barrho_1,\barrho_2,\trho_1,\trho_2] =  \int_{\R^2} \Lambda_1(x) \trho_1(x) \dx + \int_{\R^2} \Lambda_2(x) \trho_2(x) \dx.
\end{equation}

In the sequel, we will consider two classes of perturbations, denoted as $\CalA$ and $\CalB$. Class $\CalA$ consists of perturbations $(\trho_1,\trho_2)$ such that each $\trho_i$ is supported in $\Om_i$ (here $i=1,2$). Class $\CalB$ is made of perturbations $(\trho_1,\trho_2)$ such that at least one $\trho_i$ has a support with a non-empty intersection with the complement $\Om_i^c$ of $\Om_i$ ($i=1,2$). Hence, $\CalA$ and $\CalB$ are disjoint and cover all possible perturbations $(\trho_1,\trho_2)$. These choices are inspired by the setup in \cite{BeTo2011}.
\smallskip

Start by taking perturbations of class $\CalA$.  Since $\trho_i$ changes sign in $ \Om_i$ ($i=1,2$), for $(\barrho_1,\barrho_2)$ to be a critical point of the energy, the first variation must vanish. From \eqref{eqn:1stvarL}, given that perturbations $\trho_i$ are arbitrary and satisfy  \eqref{eqn:massc}, one finds that $E_1$ vanishes provided $\Lambda_i$ is constant in $\Om_i$, i.e., 
\begin{equation}
\label{eqn:equilsup}
\Lambda_1(x) = \lambda_1 \quad \text{ for } x \in \Om_1, \quad \text{ and } \qquad \Lambda_2(x) = \lambda_2 \quad \text{ for } x \in \Om_2.
\end{equation}

Equation \eqref{eqn:equilsup} represents a necessary condition for $(\barrho_1,\barrho_2)$ to be an equilibrium. For $(\barrho_1,\barrho_2)$ that satisfy \eqref{eqn:equilsup} to be a local minimizer with respect to class $\CalA$ perturbations, the second variation \eqref{eqn:2ndvar} must be non-negative. In general, the sign of $E_2$ cannot be assessed easily.

Consider now perturbations of class $\CalB$. Since perturbations $\trho_i$ must be non-negative in the complement $\Om_i^c$ of $\Om_i$, one can extend the argument in \cite{BeTo2011} and show that a necessary and sufficient condition for $E_1 \geq 0$ is 
\begin{equation}
\label{eqn:equilcomp}
\Lambda_1(x) \geq \lambda_1 \quad \text{ for } x \in \Om_1^c, \quad \text{ and } \quad \Lambda_2(x) \geq \lambda_2 \quad \text{ for } x \in \Om_2^c.
\end{equation}
Indeed, suppose an equilibrium $(\barrho_1,\barrho_2)$ satisfies \eqref{eqn:equilsup} and \eqref{eqn:equilcomp}. Then,
\begin{align*}
E_1 &=  \int_{\Om_1} \underbrace{\Lambda_1(x)}_{=\lambda_1} \trho_1(x) \dx +  \int_{\Om_1^c} \underbrace{\Lambda_1(x)}_{\geq \lambda_1} \trho_1(x) \dx+ \int_{\Om_2} \underbrace{\Lambda_2(x)}_{=\lambda_2} \trho_2(x) \dx + \int_{\Om_2^c} \underbrace{\Lambda_2(x)}_{\geq \lambda_2} \trho_2(x) \dx \\
& \geq \lambda_1 \int_{\R^2}  \trho_1(x) \dx + \lambda_2 \int_{\R^2} \trho_2(x) \dx,
\end{align*}
where we also used that $\trho_i \geq 0$ in $\Om_i^c$. By \eqref{eqn:massc} one concludes $E_1 \geq 0$. 

Conversely, suppose that \eqref{eqn:equilcomp} does not hold; assume for instance that $\Lambda_1(x) < \lambda_1$, for $x$ in a set $A \subset \Om_1^c$ of non-zero Lebesgue measure. Then, by taking $\trho_2=0$ and perturbations $\trho_1$ that are supported on $\Om_1$ and $A$, we have
\begin{equation*}
E_1 =  \int_{\Om_1} \underbrace{\Lambda_1(x)}_{=\lambda_1} \trho_1(x) \dx +  \int_{A} \underbrace{\Lambda_1(x)}_{< \lambda_1} \trho_1(x) \dx. \\
\end{equation*}
Again, by \eqref{eqn:massc}, one finds $E_1<0$, which completes the argument.

The interpretation of \eqref{eqn:equilcomp} is that transporting mass from $\Om_i$ into its complement $\Om_i^c$ increases the total energy \cite{BeTo2011}. In summary, a critical point 
$(\barrho_1,\barrho_2)$ for the energy satisfies the Fredholm integral equation \eqref{eqn:equilsup} on its support. The critical point is a local minimum with respect to perturbations of class $\CalA$ if the second variation is non-negative for such perturbations. Also, $(\barrho_1,\barrho_2)$ is a local minimizer with respect to perturbations of class $\CalB$ if it satisfies \eqref{eqn:equilcomp}. Note however that the word {\em local} in this context refers to the small {\em size} of the perturbations, as the perturbations themselves are in fact nonlocal in space. 

To establish whether an equilibrium is a global minimizer, one needs to investigate closely the second variation $E_2$ for general perturbations. From \eqref{eqn:Equad} we see that a sufficient condition for a local minimizer to be global minimizer is that $E_2 \geq 0$. Such condition is not necessary though, as \eqref{eqn:Equad} is exact, and for a global minimum one needs in fact $\varepsilon E_1 + \varepsilon^2 E_2 \geq 0$, for arbitrary $\varepsilon>0$.

The variational framework above holds true for general interaction potentials $K_s$ and $K_c$. Next, we will discuss it further, for the specific choice \eqref{eqn:intpot}.

\paragraph{Second variation of the energy. } We elaborate briefly on the second variation of the energy $E_2$  (see \eqref{eqn:Equad}) that corresponds to the interaction potentials \eqref{eqn:intpot}. This calculation is used to show that certain equilibria are global minimizers. By \eqref{eqn:2ndvar} and \eqref{eqn:energy}, we can write
\begin{equation}
\label{eqn:E2-decomp}
E_2(\trho_1,\trho_2) =  I + II,
\end{equation}
where
\begin{align*}
I := & -\frac{1}{2} \as \iint \ln |x-y| (\trho_1(x) \trho_1(y) + \trho_2(x) \trho_2(y)) \, dx dy- \ac \iint \ln |x-y| \trho_1(x) \trho_2(y) \, dx dy, \\
II :=& \frac{1}{4} \bs \iint |x-y|^2 (\trho_1(x) \trho_1(y) + \trho_2(x) \trho_2(y)) \, dx dy+ \frac{1}{2} \bc \iint |x-y|^2 \trho_1(x) \trho_2(y) \, dx dy.
\end{align*}

The expression for $II$ can be easily simplified by expanding $|x-y|^2 = |x|^2 - 2 x \cdot y + |y|^2$ and using the conservation properties \eqref{eqn:massc} and \eqref{eqn:cmassc} of the perturbations. One finds
\begin{align}
\label{eqn:est-II}
II &= -\frac{1}{2} \bs \left( \int x \trho_1(x) \, dx \right)^2  -\frac{1}{2} \bs \left( \int x \trho_2(x) \, dx \right)^2 - \bc \left(\int x \trho_1(x) \, dx \right) \left( \int x \trho_2(x) \, dx \right) \nonumber \\
&= (\bc - \bs)  \left( \int x \trho_1(x) \, dx \right)^2.
\end{align}

For $I$ we use the Plancherel's theorem:
\begin{align*}
I &= -\frac{1}{4 \pi} \as \int_{\R^2} \CalF \{\ln |x| \ast \trho_1\}(k) \, \overline{\CalF \{\trho_1\}} (k) \, dk  -\frac{1}{4 \pi} \as \int_{\R^2} \CalF \{\ln |x| \ast \trho_2\}(k) \, \overline{\CalF \{\trho_2\}} (k) \, dk \\
& \quad - \frac{1}{2 \pi}\ac  \int_{\R^2} \CalF \{\ln |x| \ast \trho_1\}(k) \, \overline{\CalF \{\trho_2\}} (k) \, dk,
\end{align*}
where $\CalF$ represents the Fourier transform
\[
\CalF \{ f \}(k) = \int_{\R^2} f(x) e^{-i k \cdot x} dx,
\]
and overbar denotes the complex conjugate.

Using $ \CalF \{\ln |x| \ast \trho_i\}(k) = \CalF \{ \ln |x| \}(k) \cdot \CalF \{\trho_i\}(k)$, and $\CalF \{\ln|x|\} (k) = - 2\pi /|k|^2$, we further arrive at
\begin{equation}\label{eqn:termI}
I=  \frac{1}{2} \as \int_{\R^2} \frac{1}{|k|^2} \left( \left | \CalF \{\trho_1\} (k) \right |^2  + \left | \CalF \{\trho_2\} (k)\right |^2 \right) dk + \ac \int_{\R^2} \frac{1}{|k|^2} \,  \CalF \{\trho_1\} (k) \cdot \overline{\CalF \{\trho_2\}} (k) dk
\end{equation}
By Cauchy-Schwarz,
\begin{align*}
\left | \int_{\R^2} \frac{1}{|k|^2} \,  \CalF \{\trho_1\} (k) \cdot \overline{\CalF \{\trho_2\}} (k) \, dk \right | &\leq  \left( \int_{\R^2} \frac{1}{|k|^2}  \left | \CalF \{\trho_1\} (k) \right |^2 dk \right)^{\frac{1}{2}} \cdot  \left( \int_{\R^2} \frac{1}{|k|^2}  \left | \CalF \{\trho_2\} (k) \right |^2 dk \right)^{\frac{1}{2}} \\
& \leq \frac{1}{2} \int_{\R^2} \frac{1}{|k|^2}  \left | \CalF \{\trho_1\} (k) \right |^2 dk + \frac{1}{2} \int_{\R^2} \frac{1}{|k|^2}  \left | \CalF \{\trho_2\} (k) \right |^2 dk,
\end{align*}
and hence, 
\begin{equation}
\label{eqn:est-I}
I \geq  \frac{1}{2} (\as - \ac) \int_{\R^2} \frac{1}{|k|^2} \left( \left | \CalF \{\trho_1\} (k) \right |^2  + \left | \CalF \{\trho_2\} (k)\right |^2 \right) dk.
\end{equation}

We also note here that an alternative way (to using Fourier transforms) for investigating term $I$ is to identify  $\iint \ln |x-y| \trho_i(x) \trho_i(y) dx dy$ with the $H^{-1}$ norm of $\trho_i$ and the cross-term as $H^{-1}$ inner product \cite{ChMaWi2011}. Then, bounding the inner product by Cauchy-Schwarz to reach \eqref{eqn:est-I} follows similarly.

\paragraph{Global, local minimizers and the sign of $E_2$.} Given that the expression \eqref{eqn:Equad} for energy is exact, the considerations above lead to some immediate conclusions. Consider an equilibrium $(\barrho_1,\barrho_2)$ and parameters such that $\as > \ac$ and $\bc >\bs$ ($A<1$ and $B>1$; the shaded region in Figure \ref{fig:stability}). Here it is assumed, of course, that the equilibrium $(\barrho_1,\barrho_2)$ exists for parameters in this regime. Assume that $(\barrho_1,\barrho_2)$ is a local minimizer with respect to perturbations of class $\CalB$ (i.e., it satisfies \eqref{eqn:equilsup}). Take now an arbitrary perturbation $(\trho_1,\trho_2)$; by the local minimizer condition, one has $E_1 \geq 0$. Furthermore, by \eqref{eqn:E2-decomp}, \eqref{eqn:est-II} and \eqref{eqn:est-I} above (which apply to {\em any} perturbation), $E_2 \geq 0$.  This implies that the equilibrium under consideration is in fact a {\em global} minimizer.

For values of the parameters outside $A<1$ and $B>1$, establishing the sign of the second variation is a challenging task. The reason lies in the very different expressions of the terms $I$ and $II$ that comprise $E_2$; the two terms are not immediately comparable and also, the logarithmic potential is not sign-definite. Balancing term $I$ and $II$ to yield a definite sign for $E_2$ seems difficult and we do not pursue this direction here. The best we can do in such regimes, using the variational method above, is to check for the sign of $E_1$ (condition \eqref{eqn:equilsup}) and restrict our conclusions to local minimizers with respect to perturbations of class $\CalB$. 

To conclude, unless $A<1$ and $B>1$, we only establish, by the variational approach, whether a specific equilibrium is a local minimizer with respect to perturbations of class $\CalB$. We do not make any conclusion concerning minimization with respect to class $\CalA$ perturbations. To compensate for this limitation, we develop and use an alternative approach for studying the stability of equilibria, based on linear analysis. We present now the main features of this alternative approach.

  
\subsection{Linear stability analysis}\label{sect:prelim lin}
\paragraph{Class of perturbations and difference with variational approach.}
In \eqref{eqn:equilibria} we identified the values that the density can attain in a steady state. In particular, any steady state consists of regions where (a) only one of the species has nonzero density, (b) the two species coexist, or (c) none of the species is present. Within each of these subdomains, the density of each species is a constant given by \eqref{eqn:equilibria}.

These properties of steady states are a specific consequence of the choice of kernels in \eqref{eqn:intpot}. The fact that the steady states are ``piecewise" constant, makes us consider here perturbations in which the densities remain constant, but the boundaries of the supports are deformed; cf.~\cite{ChKo2014}. We consider deformations such that the total area enclosed remains unchanged (up to higher-order contributions), due to the constraint of fixed total mass.

The perturbations considered here are completely different in spirit from the ones used in the variational approach. The difference lies in the meaning of the word `small', when we speak about `small perturbations'. Here, `small' means that only \emph{close} to the boundaries of the supports, the density may change. That is, exactly at those points that were outside the support of the unperturbed equilibrium, but now fall within the support of the perturbed state (upon alteration of the boundaries), or \textit{vice versa}. The change in density is $\Ord(1)$ at those points, according to the discrete set of values allowed by \eqref{eqn:equilibria}. In the variational approach, perturbations change the value of the density only \emph{slightly} (cf. multiplication by $\eps$ in \eqref{eqn:pert dens var}), but on the other hand, are allowed to be {\em nonlocal} in space (in particular, not necessarily close to the equilibrium's support). 

\paragraph{Mathematical description.}
In Figure \ref{fig:num st st} we identified several steady states numerically, among which is the one we called the `target' (see also equilibria shown schematically in Figures \ref{fig:target} and \ref{fig:target heavy inside}). We will apply the linear stability analysis to such equilibrium states. A short explanation of why we only consider these states follows at the end of this section. As we focus on target states, we only need to perturb \emph{circular} boundaries, which simplifies the exposition here.

For mathematical convenience, we identify the domain $\R^2$ with the complex plane.  Similarly to \cite{KolKevCar2014}, we consider the following perturbations (corresponding to Fourier mode $m\in\N^+$) of the circle with radius $R_j$:
\begin{equation}\label{eqn:perturb}
p_j(\theta):= R_j\,e^{i\,\theta}\left( 1 + \eps_{j,N}\,\cos(m\,\theta) + i\,\eps_{j,T}\,\sin(m\,\theta) \right),\qquad \text{ for } 0\leqs\theta< 2\pi,
\end{equation}
where $\eps_{j,N}$ and $\eps_{j,T}$ are assumed to be small parameters that control the normal (`N') and transversal (`T') deformation. The index $j$ takes values $0,1$ and $2$; the exact numbering is configuration-specific (see Figures \ref{fig:target} and \ref{fig:target heavy inside}). Let $\Omega_j^{\eps}$ denote the perturbed domain enclosed by $p_j(\theta)$; note that it depends on $m$. Roughly speaking, due to perturbations \eqref{eqn:perturb}, a number of ``oscillations" are superimposed on the unperturbed circular boundary. The number of oscillations is determined by the mode $m$. Here we have in mind the idea that any \textit{arbitrary} perturbation can be obtained by using its decomposition in Fourier modes \eqref{eqn:perturb}.

The area of $\Omega_j^{\eps}$ is $\pi R_j^2 + \frac12\pi R_j^2(\eps_{j,N}^2+\eps_{j,T}^2)$, independent of the mode $m$. For $m\geqs 2$, the perturbations in \eqref{eqn:perturb} preserve the centre of mass. For $m=1$, we have that
\[ \int_{\Omega_j^\eps}y\,dy=\pi R_j^3 \eps_{j,N} + \frac14 \pi R_j^3 (\eps_{j,N}^2-\eps_{j,T}^2)(\eps_{j,N}+\eps_{j,T}).
\]

Recall that in this approach we consider perturbations from equilibrium such that the densities remain constant (at the same values as at equilibrium) within the perturbed domains. To assess the stability, we investigate the dynamics of a generic point on the perturbed boundaries. We use the index $j$ here to denote that specific boundary, and the point we observe is $x=p_j(\theta_0)$ for some $0\leqs\theta_0<2\pi$. We need to calculate the velocity given in \eqref{eqn:model} at position $x=p_j(\theta_0)$, and thus we need to evaluate convolution integrals of $-\nabla K_s$ and $-\nabla K_c$ against the densities $\rho_1$ and $\rho_2$ over domains $\Omega_\ell^{\eps}$ ($\ell=0,1,2$). 
As densities have been assumed to remain constant, one can take the density values outside the integral, while the information about the steady state is accounted for by the specific integration domain. Consequently, the velocity at $x\in\R^2$ becomes a weighted sum over $\ell$ of integrals of the type
\begin{equation}\label{eqn:basic int lin}
\int_{\Omega_\ell^{\eps}} \nabla K(x-y)\,dy.
\end{equation}
Note for instance, that the integral over the annulus in the target state is obtained by subtracting integrals like these, for two different values of $\ell$. The potential $K$ in \eqref{eqn:basic int lin} either denotes $K_s$ or $K_c$. By some abuse of notation, we used $\Omega_\ell^{\eps}\subset\C$ for the integration domain in $\R^2$.


Due to the choice of potentials in \eqref{eqn:intpot}, the integral \eqref{eqn:basic int lin} can be written in terms of
\begin{equation}\label{eqn:pert int}
\int_{\Omega_\ell^\eps} \dfrac{x-y}{|x-y|^2}\,dy,
\qquad \text{and} \qquad
\int_{\Omega_\ell^\eps} (x-y)\,dy.
\end{equation}
The latter integral is relatively easy to evaluate exactly, using the aforementioned expressions for the area and centre of mass of $\Omega_\ell^\eps$. The outcome is given in \eqref{eqn:eval attr int} of Appendix \ref{app:pert integrals}, where higher-order terms in $\eps_{\ell,N}$ and $\eps_{\ell,T}$ are omitted. We emphasize that the specific choice $x=p_j(\theta_0)$ introduces a dependence on $\eps_{j,N}$ and $\eps_{j,T}$ as well.


By Gauss' theorem, the left-hand integral in \eqref{eqn:pert int} can be transformed into a contour integral over the boundary $\partial \Omega_\ell^\eps$. We have
\begin{align}\label{eqn:int ln Gauss}
\int_{\Omega_\ell^\eps} \dfrac{x-y}{|x-y|^2}\,dy = -\int_{\partial\Omega_\ell^\eps} \ln|x-y|\,\hat{n}\,dS,
\end{align}
where we used that $(x-y)/|x-y|^2=-\nabla_y \ln|x-y|$. The boundary of $\Omega_\ell^\eps$ is parameterized by $p_\ell(\theta)$, for $0\leqs\theta< 2\pi$. Consequently, $\hat{n}\,dS$ can be expressed in terms of $\theta$, $\eps_{\ell,N}$ and $\eps_{\ell,T}$; see \eqref{eqn:n dS} in Appendix \ref{app:pert integrals}. 

The resulting one-dimensional integral \eqref{eqn:ln int} depends in a nonlinear way on the small parameters $\eps_{j,N}$, $\eps_{j,T}$, $\eps_{\ell,N}$ and $\eps_{\ell,T}$. In our linear stability analysis, we expand the integrand in terms of the small parameters, we omit higher-order terms and end up with integrals that can be evaluated. See Appendix \ref{app:pert integrals} for more details. The first-order approximation of \eqref{eqn:int ln Gauss} that we find, depends on whether $x=p_j(\theta_0)$ is inside, on or outside $\partial \Omega^\eps_\ell$.

\paragraph{System of linearized ODE's.}
Assume that our point of interest $x=p_j(\theta_0)$ is on the boundary of the support of species $k$. For instance, for the target of Figure \ref{fig:target}, if $x$ is on the boundary of the inner disk, then $k=2$, and if $x$ is on one of the boundary of the annular region, then $k=1$. Any point $z$ in the support of species $k$ has velocity $v_k(z)$ and hence we know that
\begin{equation}\label{eqn:dpdt is v pert}
\frac{d}{dt}p_j(\theta_0)=v_k(p_j(\theta_0))
\end{equation}
holds. Note that both $p_j$ and $v_k$ depend on all $\eps$'s. Let $\underline\eps$ be the vector of all these $\eps$'s. In case of the target, $\underline\eps$ has six components. We want to investigate the stability of the state $\underline\eps=0$ (i.e.~the target equilibrium state in Figure \ref{fig:target}). We do so by considering $\underline{\eps}=\underline\eps(t)$ and linearizing \eqref{eqn:dpdt is v pert} around $\underline\eps=0$, for arbitrary $\theta_0\in[0,2\pi)$ and for all modes $m\in\N^+$. The same procedure is repeated for multiple indices $j$, taking into account all boundaries in the considered steady state. For a target, there are three such boundaries. 



On one hand, it follows from \eqref{eqn:perturb} that
\begin{equation}\label{eqn:dpdt}
\dfrac{d}{dt}p_j(\theta_0) = R_j\,e^{i\,\theta_0}\left(\eps'_{j,N}(t)\,\cos(m\,\theta_0) + i\,\eps'_{j,T}(t)\,\sin(m\theta_0) \right).
\end{equation}
On the other hand, $v_k(p_j(\theta_0))$ can be evaluated (up to higher-order terms in $\underline\eps$) using the results of Appendix \ref{app:pert integrals}. Combined, the $\Ord(1)$ terms vanish; this is a necessary condition for a steady state.

For each $j$, we combine \eqref{eqn:dpdt is v pert}, \eqref{eqn:dpdt} and the expression for $v_k(p_j(\theta_0))$ based on Appendix \ref{app:pert integrals}. 
Dividing by $R_j\exp(i\,\theta_0)$ and matching sine and cosine terms on both sides of the equation, we obtain a system of six ODE's of the form
\begin{equation}\label{eqn:ODE eps Q}
\dfrac{d}{dt}\underline\eps = Q_m\,\underline\eps,
\end{equation}
for $m\in\N^+$. 
In particular, this system turns out to be independent of the choice of $\theta_0$ (see Section \ref{sect:target}). The eigenvalues of the matrix $Q_m$ need to be found, and their sign yields the stability of the steady state. The inspection of the eigenvalues for arbitrary mode $m$ takes place in Sections \ref{sec:target1 lin} and \ref{subsect:lin t-hinside}.
\paragraph{Limitations of this approach.}
We restricted ourselves to circle-shaped domain boundaries primarily because we do not have an exact mathematical description for the other shapes appearing in Figure \ref{fig:num st st}. Such description is required to perform the analysis of this section, starting from a modified form of \eqref{eqn:perturb}.

Even if we did manage to find exact formulas for non-circular boundaries, it is not directly clear if we could obtain results analogous to those in Appendix \ref{app:pert integrals}. For circular boundaries, the final result in Appendix \ref{app:pert integrals} is based on integrals of the form \eqref{eqn:log integral} and \eqref{eqn:rational integral alpha not 1}, for which we have exact expressions. For non-circular boundaries a different parameterization in $\theta$ of \eqref{eqn:int ln Gauss} is needed. Linearization in $\eps$ of the right-hand side may again reduce the problem to the evaluation of certain basic integrals, but it is all but certain if these can be solved exactly.

We manage to analyze the target steady states by our linear perturbation method; cf.~Sections \ref{sec:target1 lin} and \ref{subsect:lin t-hinside}. These states have the advantage that there is an $\Ord(1)$ distance between the disk-shaped core and the annular region outside. Hence, for sufficiently small $\eps_{j,N}$'s and $\eps_{j,T}$'s, the three perturbed boundaries do not interfere. 

Now consider the `overlap' states (e.g.~the one shown for $A=0.5$ and $B=1$ in Figure \ref{fig:num st st}). The boundaries of the supports are circular, and hence \eqref{eqn:perturb} could in principle still be used. Examine however the boundary of the support of the lighter (red) species. For the heavier (blue) species, this \emph{same} circle is the separating curve between the region of coexistence with the red species (central circle), and the outer annulus in which only the blue species is present. This implies that a point $x=p_j(\theta_0)$ on this inner boundary needs to satisfy \emph{two} equations of the form \eqref{eqn:dpdt is v pert}: once for $k$ being the blue species, and once for $k$ being the red species. For the overlap state, we performed the corresponding calculations, but obtained two inconsistent equations. We did not manage to resolve this issue. It indicates however that interfering boundaries (even if they are circles) may lead to difficulties in our method.

An extra complication lies in the fact that for every boundary that appears in a certain steady state, system \eqref{eqn:ODE eps Q} contains two variables. Ultimately, (the signs of) the eigenvalues 
of the matrix $Q_m$ in \eqref{eqn:ODE eps Q} need to be determined. For $Q_m$ larger than $3\times3$, finding the eigenvalues analytically is in general not possible, and one needs to rely on other techniques. In this paper we draw conclusions for instance by inspecting the coefficients of the characteristic polynomial. An alternative is to use Gershgorin's theorem (which however does not give a conclusive answer about the sign in the cases treated in this paper) or eventually numerics.

\subsection{Discrete model and numerical investigation of equilibria}\label{sec:particle model and num}
There is a particle system formulation that follows immediately from \eqref{eqn:model} and this system of ODE's is the basis for the numerical investigations done in this paper. We consider a total number of $N$ particles, distributed over the two species, such that the two populations are $\{x^{(1)}_i\}_{i=1}^{N_1}\subset \R^2$ and $\{x^{(2)}_i\}_{i=1}^{N_2}\subset \R^2$, respectively, with $N_1+N_2=N$. 

The discrete analogue of model \eqref{eqn:model} is given by the following system of ODE's:
\begin{subequations}
\label{eqn:part syst}
\begin{gather}
\dfrac{dx^{(1)}_i}{dt} = -\dfrac{M_1}{N_1}\sum_{\substack{j=1\\j\neq i}}^{N_1}\nabla K_s\left(x^{(1)}_i-x^{(1)}_j\right) -\dfrac{M_2}{N_2}\sum_{j=1}^{N_2}\nabla K_c\left(x^{(1)}_i-x^{(2)}_j\right), \quad i=1,\ldots,N_1,\\
\dfrac{dx^{(2)}_i}{dt} = -\dfrac{M_1}{N_1}\sum_{j=1}^{N_1}\nabla K_c\left(x^{(2)}_i-x^{(1)}_j\right) -\dfrac{M_2}{N_2}\sum_{\substack{j=1\\j\neq i}}^{N_2}\nabla K_s\left(x^{(2)}_i-x^{(2)}_j\right), \quad i=1,\ldots,N_2.
\end{gather}
\end{subequations}
Note that the four summations in the right-hand sides correspond (except for the omission of the self-interaction term) to convolutions of the interaction kernels with the empirical measures $\mu_k:=1/N_k \sum_{i=1}^{N_k}\delta_{x^{(k)}_i}$.

We investigate the steady states of the particle system numerically, by performing long-time simulations of \eqref{eqn:part syst} starting from random initial data. The steady states of the particle system are expected to capture the steady states of the PDE model \eqref{eqn:model}. Particle simulations are in fact the main tool used to study numerically equilibria of the one species model \cite{Balague_etalARMA,KoSuUmBe2011,Brecht_etal2011,FeHuKo11}.






\section{Target equilibrium}
\label{sect:target}
In this section we investigate the radially symmetric state where the two species are supported concentrically on a disk and an annulus, respectively. See the left-hand picture for parameter values $(A,B)=(3,3.5)$ in Figure \ref{fig:num st st}. We also consider this state with the heavy and light species interchanged. Both of them we call ``targets".
\subsection{Lighter species inside}
\label{subsect:t-linside}
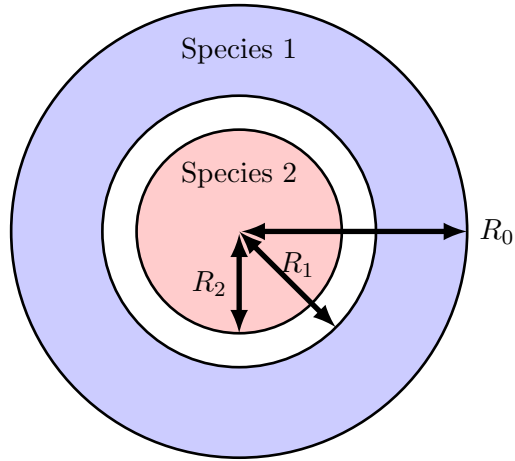
\begin{figure}[t]
\centering
\begin{tikzpicture}[>= latex]
\begin{scope}[scale=0.3] 

\draw[fill=blue!20
,even odd rule,line width=1]  (0,0) circle (10)
                                   (0,0) circle (6);
\draw[fill=red!20
,line width=1]  (0,0) circle (4.5);            
                
\draw[<->,line width=2] (0.1,0) -- (10,0) node[pos=1,right] {$R_0$};
\draw[<->,line width=2] (0,0) -- +(-45:6) node[pos=.6,above] {$R_1$};
\draw[<->,line width=2] (0,-0.1) -- (0,-4.5) node[pos=.5,left] {$R_2$};
\draw (0,2.5) node {Species 2};
\draw (0,8) node {Species 1};

\end{scope}
\end{tikzpicture}
\caption{Target equilibrium, with species $1$ supported on the annular region $\Rone<|x|<\Rzero$, and (the lighter) species $2$ supported on a disk of radius $\Rtwo$.}
\label{fig:target}
\end{figure}

Consider the target configuration sketched in Figure \ref{fig:target}, where the heavier species $1$ is supported on an annular region $\Rone<|x|<\Rzero$, and the lighter species $2$ is supported on a disk of radius $\Rtwo$, with $\Rtwo \leq \Rone \leq \Rzero$.  Following simple calculations, the radii are given by
\begin{equation}
\label{eqn:radii}
\Rtwo^2 = \frac{\as M_2}{\bc M_1 + \bs M_2}, \qquad \Rone^2 =  \frac{\ac M_2}{\bs M_1 + \bc M_2}, \qquad \Rzero^2 = \frac{\as M_1 + \ac M_2}{\bs M_1 + \bc M_2}.
\end{equation}
Within the respective (non-overlapping) supports $\Om_1 = \{ \Rone < |x| < \Rzero\}$ and $\Om_2 = \{ |x|<\Rtwo \}$, the equilibrium densities are  (cf. \eqref{eqn:equilibria}):
\begin{equation}
\label{eqn:rho12-target}
\barrho_1=\frac{\bs M_1 + \bc M_2}{\pi \as}, \qquad \barrho_2= \frac{\bc M_1 + \bs M_2}{\pi \as}.
\end{equation}

For consistency with the target solution ansatz, $(\A,\B)$ has to lie in $\D_3 \cup \D_4 \cup \D_5$.

\subsubsection{Variational approach}
\label{subsubsect:t-var-linside}
We now proceed to investigate whether the target equilibrium is a local minimizer with respect to class $\CalB$ perturbations, cf. \eqref{eqn:equilcomp}. The following elementary calculations will be needed in the sequel:
\begin{equation}
\label{eqn:intlog}
\int_{|y|<R} \ln |x-y| \dy = \begin{cases}
\frac{\pi}{2}|x|^2 + \pi R^2 \ln R - \frac{\pi}{2}R^2& \text{ if } |x|<R \\
\pi R^2 \ln |x| &  \text{ if } |x|>R,
\end{cases}
\end{equation}
and
\begin{equation}
\label{eqn:intderiv}
\int_{|y|<R} \frac{x-y}{|x-y|^2} \dy = \begin{cases}
\pi x & \text{ if } |x|<R \\ 
\pi R^2 \frac{x}{|x|^2} &  \text{ if } |x|>R.
\end{cases}
\end{equation}
Note in fact that \eqref{eqn:intderiv} can be derived from \eqref{eqn:intlog} by differentiation.

Calculate $\nabla \Lambda_1$ from \eqref{eqn:Lambda1}, with the potentials given by \eqref{eqn:intpot} and the equilibrium densities given by \eqref{eqn:rho12-target}. Using \eqref{eqn:intlog} and \eqref{eqn:intderiv} one can check indeed that $\nabla \Lambda_1(x)=0$ in $\Rone<|x|<\Rzero$, hence $\Lambda_1(x)$ is constant in the support $\Om_1$ (see \eqref{eqn:equilsup}). The calculation of $\nabla \Lambda_1(x) = \Lambda_1'(|x|)\, x/|x|$ (note the radial symmetry) outside the support $\Om_1$ yields the following:
\begin{equation}
\label{eqn:gradL1-target}
\frac{\Lambda_1'(|x|)}{|x|}= \begin{cases}
\bs M_1 + \bc M_2 - \frac{\ac}{\as} (\bc M_1 + \bs M_2)& \text{ if } |x|<\Rtwo \\ 
\bs M_1 + \bc M_2 - \ac M_2/|x|^2&  \text{ if } \Rtwo<|x|<\Rone \\
\bs M_1 + \bc M_2 - (\as M_1+ \ac M_2)/|x|^2  &  \text{ if } |x|>\Rzero.
\end{cases}
\end{equation}

Using the notations \eqref{eqn:ABM}, $\Lambda_1'$ in $|x|<\Rtwo$ can be written as:
\[
\Lambda_1'(|x|) =  \bs M_2 (M+B - A(BM+1)) |x|, \qquad \text{ in } |x|<\Rtwo.
\]
For all $(A,B)$ in $\D_3 \cup \D_4 \cup \D_5$ (see Figure \ref{fig:stability}), which is the entire parameter space where the assumed target equilibrium exists, the expression above is negative. Indeed, for such $(A,B)$, $B>(M-A)/(MA-1)$ (note that $MA-1>0$), and hence, $ABM - B - M +A >0$.

In $\Rtwo<|x|<\Rone$, $\Lambda_1(|x|)$ is also decreasing, as can be seen from the simple estimate below: 
\[
\bs M_1 + \bc M_2 - \frac{\ac M_2}{|x|^2} < \underbrace{\bs M_1 + \bc M_2 - \frac{\ac M_2}{{\Rone}^2}}_{=0 \text{ by } \eqref{eqn:radii}}.
\]
Finally, in $|x|>\Rzero$, $\Lambda_1(|x|)$ is increasing, as
\[
\bs M_1 + \bc M_2 - \frac{\as M_1+ \ac M_2}{|x|^2} > \underbrace{\bs M_1 + \bc M_2 - \frac{\as M_1+ \ac M_2}{\Rzero^2}}_{=0 \text{ by } \eqref{eqn:radii}}.
\]

In summary, for all $(A,B)$ in the relevant region $\D_3 \cup \D_4 \cup \D_5$, $\Lambda_1$ satisfies \eqref{eqn:equilcomp}; for an illustration see Figure \ref{fig:Lambdas}.

We now calculate $\nabla \Lambda_2$ from \eqref{eqn:Lambda2}. By \eqref{eqn:intpot} and \eqref{eqn:rho12-target}, also using \eqref{eqn:intlog} and \eqref{eqn:intderiv} one can check indeed that $\nabla \Lambda_2(x)=0$ in $|x|<\Rtwo$, hence $\Lambda_2(x)$ is constant in the support $\Om_2$ (see \eqref{eqn:equilsup}). The calculation of $\nabla \Lambda_2(x) = \Lambda_2'(|x|)\, x/|x|$ outside the support $\Om_2$ yields:
\begin{equation}
\label{eqn:gradL2-target}
\frac{\Lambda_2'(|x|)}{|x|}= \begin{cases}
\pi \as \barrho_2 (1-\Rtwo^2/|x|^2) & \text{ if } \Rtwo<|x|<\Rone \\ 
\pi \as \barrho_2 (1-\Rtwo^2/|x|^2) - \pi \ac \barrho_1 (1-\Rone^2/|x|^2) &  \text{ if } \Rone<|x|<\Rzero \\
\bc M_1 + \bs M_2 - (\ac M_1 + \as M_2)/|x|^2  &  \text{ if } |x|>\Rzero. 
\end{cases}
\end{equation}

From \eqref{eqn:gradL2-target} we infer that $\Lambda_2$ is increasing in the radial direction in the region $\Rtwo<|x|<\Rone$. On the other hand, $\Lambda_2'$ can become zero in $\Rone<|x|<\Rzero$, and consequently $\Lambda_2$ can decrease in this region. Let us investigate this scenario. The zero of $\Lambda_2'$ occurs at
\begin{equation}
\label{eqn:zero1}
|x|^2 = \frac{\pi \ac \barrho_1 \Rone^2 - \pi \as \barrho_2 \Rtwo^2}{\pi \ac \barrho_1-\pi \as \barrho_2}.
\end{equation}
For consistency, the expression above needs to be positive and also, it has to lie in the annular region (i.e., $\Rone^2<|x|^2<\Rzero^2$). By using \eqref{eqn:rho12-target} and notations \eqref{eqn:ABM}, the denominator in \eqref{eqn:zero1} reduces to:
\[
\pi \ac \barrho_1-\pi \as \barrho_2 = M_2 \bs (M(A-B)+ AB-1).
\]

Recall that the target solution only exists for parameters $(A,B)$ in the region $\D_3 \cup \D_4 \cup \D_5$. It is immediate to show that the expression above is negative in $\D_5$ and positive in $\D_3 \cup \D_4$. Consider first the case when it is negative, i.e., $(A,B) \in \D_5$. In this case, the zero of $\Lambda_2'$ from \eqref{eqn:zero1} does not lie in the relevant region $\Rone<|x|<\Rzero$. In other words, for $(A,B) \in \D_5$, $\Lambda_2'$ does not change sign in the annular region and remains positive throughout (one can check easily for instance that $\Lambda_2'(\Rone) >0$).

For $(A,B) \in \D_3 \cup \D_4$, where the denominator in \eqref{eqn:zero1} is positive, the location of the zero given by \eqref{eqn:zero1} is larger than $\Rone$ for all $(A,B)$. By requiring that it is also less than $\Rzero$, we arrive after some elementary calculations at the following condition:
\[
(M_1^2-M_2^2) (A-B) >0.
\]
Since $M_1>M_2$, this reduces simply to $A>B$. Hence, only for $(A,B) \in \D_3$, $\Lambda_2'$ can have a zero in the annular region between $\Rone$ and $\Rzero$.

To summarize the finding above, for $(A,B) \in \D_4 \cup \D_5$, $\Lambda_2'(|x|)$ does not change sign and remains positive in $\Rone<|x|<\Rzero$. For $(A,B) \in \D_3$, $\Lambda_2$ changes monotonicity in the annular region, and once it changes monotonicity, it stays decreasing through the rest of $\Rone<|x|<\Rzero$ -- see Figure \ref{fig:Lambdas} for an illustration. One can check in fact that indeed, at $|x|=\Rzero$,
\[
\Lambda_2'(\Rzero) = \frac{M_1^2-M_2^2}{M_1 \as + M_2 \ac} \as \bs (B-A) \Rzero <0 \qquad \text { for } (A,B) \in \D_3.
\]

In $|x|>\Rzero$, it can be shown easily that $\Lambda_2$ remains strictly increasing for $(A,B) \in \D_4 \cup \D_5$. Consequently, combining with the findings above, the condition for $\Lambda_2$ in \eqref{eqn:equilcomp} holds for all parameter values $(A,B) \in \D_4 \cup \D_5$ -- see Figure \ref{fig:Lambdas}(a). For $(A,B) \in \D_3$ however, it remains to be checked whether $\Lambda_2$ drops in $|x|>\Rzero$ below $\lambda_2$, the value it has on the support. From \eqref{eqn:gradL2-target}, one can infer immediately that $\Lambda_2$ changes monotonicity (again) in $|x|>\Rzero$, at 
\begin{equation}
\label{eqn:zero2}
|x|^2 = \frac{\ac M_1 + \as M_2}{\bc M_1 + \bs M_2} > \Rzero^2 \qquad \text { for } (A,B) \in \D_3.
\end{equation}
If $\Lambda_2$ evaluated at the minimum point above, drops below $\lambda_2$, then the condition for $\Lambda_2$ in \eqref{eqn:equilcomp} fails, and the target solution is not a minimizer.

Using \eqref{eqn:intlog} we calculate $\Lambda_2(x)$ in $|x|<\Rtwo$ (where it has constant value $\lambda_2$) and at the location \eqref{eqn:zero2}; call the latter $\lm$ (see Figures \ref{fig:Lambdas}(b) and (c)).  We find
\begin{equation}
\label{eqn:l2-target}
\begin{aligned}
\lambda_2 &= \frac{1}{2} (\ac M_1 +  \as M_2) - \barrho_1 \ac (\pi \Rzero^2 \log \Rzero - \pi \Rone^2 \log \Rone) - \barrho_2 \as \pi \Rtwo^2 \log \Rtwo \\
& \quad + \frac{1}{4} (\bc M_1 (\Rzero^2 + \Rone^2) + \bs M_2 \Rtwo^2),
\end{aligned}
\end{equation}
and
\begin{equation}
\label{eqn:lm-target}
\begin{aligned}
\lm &= \frac{1}{2} (\ac M_1 + \as M_2) - \frac{1}{2} (\ac M_1 + \as M_2) \log\left( \frac{\ac M_1 + \as M_2}{\bc M_1 + \bs M_2}\right) \\
& \quad + \frac{1}{4} (\bc M_1 (\Rzero^2 + \Rone^2) + \bs M_2 \Rtwo^2).
\end{aligned}
\end{equation}
Note that the first term in the right-hand-sides of \eqref{eqn:l2-target} and \eqref{eqn:lm-target}, as well as the expressions on the second lines of the respective right-hand-sides, are the same. Also, by adding and subtracting $\barrho_1 \ac \pi \Rone^2 \log \Rzero$ one can write 
\[
\barrho_1 \ac (\pi \Rzero^2 \log \Rzero - \pi \Rone^2 \log \Rone) +  \barrho_2 \as \pi \Rtwo^2 \log \Rtwo  = \ac M_1 \log \Rzero + \as M_2 \log \Rtwo + \barrho_1 \ac \pi \Rone^2 \log \left( \frac{\Rzero}{\Rone} \right).
\]

Use the above in the expression \eqref{eqn:l2-target} for $\lambda_2$. The target is not a minimizer if $\lm < \lambda_2$. By \eqref{eqn:l2-target} and \eqref{eqn:lm-target}, this occurs when
\[
 \frac{1}{2} (\ac M_1 + \as M_2) \log\left( \frac{\ac M_1 + \as M_2}{\bc M_1 + \bs M_2}\right) > \ac M_1 \log \Rzero + \as M_2 \log \Rtwo + \barrho_1 \ac \pi \Rone^2 \log \left( \frac{\Rzero}{\Rone} \right)
\]

By using \eqref{eqn:radii} and notations \eqref{eqn:ABM}, we can reduce the inequality above to
\[
AM \log \left( \frac{BM+1}{B+M}\right) - AM \log \left( \frac{AM+1}{A+M}\right) - \log(1+AM) < -A^2 \log \left( 1+ \frac{M}{A}\right),
\]
or, removing the $\log$,
\begin{equation}
\label{eqn:target-nonmin1}
\frac{BM+1}{B+M} < \frac{AM+1}{A+M} \cdot (1+AM)^{\frac{1}{AM}} \cdot \left( 1+ \frac{M}{A}\right)^{-\frac{A}{M}}.
\end{equation}
Denote the right-hand-side by $\f(A)$:
\[
\f(A) = \frac{AM+1}{A+M} \cdot (1+AM)^{\frac{1}{AM}} \cdot \left( 1+ \frac{M}{A}\right)^{-\frac{A}{M}}.
\]
By numerical inspection,  $\frac{1}{M}<\f(A)<M$, for $A>1$, and hence \eqref{eqn:target-nonmin1} can be written explicitly as 
\begin{equation}
\label{eqn:target-nonmin2}
B< \frac{M \f(A)-1}{M-\f(A)}.
\end{equation}

Figure \ref{fig:nonmin} shows (shaded area) the subset of $D_3$ where \eqref{eqn:target-nonmin2} holds; for parameters $(A,B)$ in this region, the target equilibrium is not a minimizer. For such $(A,B)$ a typical profile of $\Lambda_2$ is illustrated in Figure \ref{fig:Lambdas}(c). We also note that as $M$ increases to infinity, the subset seems to approach the entire domain $D_3$.

For $(A,B)$ in $D_3$ outside the shaded region, \eqref{eqn:target-nonmin2} is violated and hence, $\lambda_m>\lambda_2$ and \eqref{eqn:equilcomp} holds; a typical profile of $\Lambda_2$ in this case is shown in Figure \ref{fig:Lambdas}(b).

\afterpage{
\begin{figure}[h]
\includegraphics[width=0.3\textwidth]{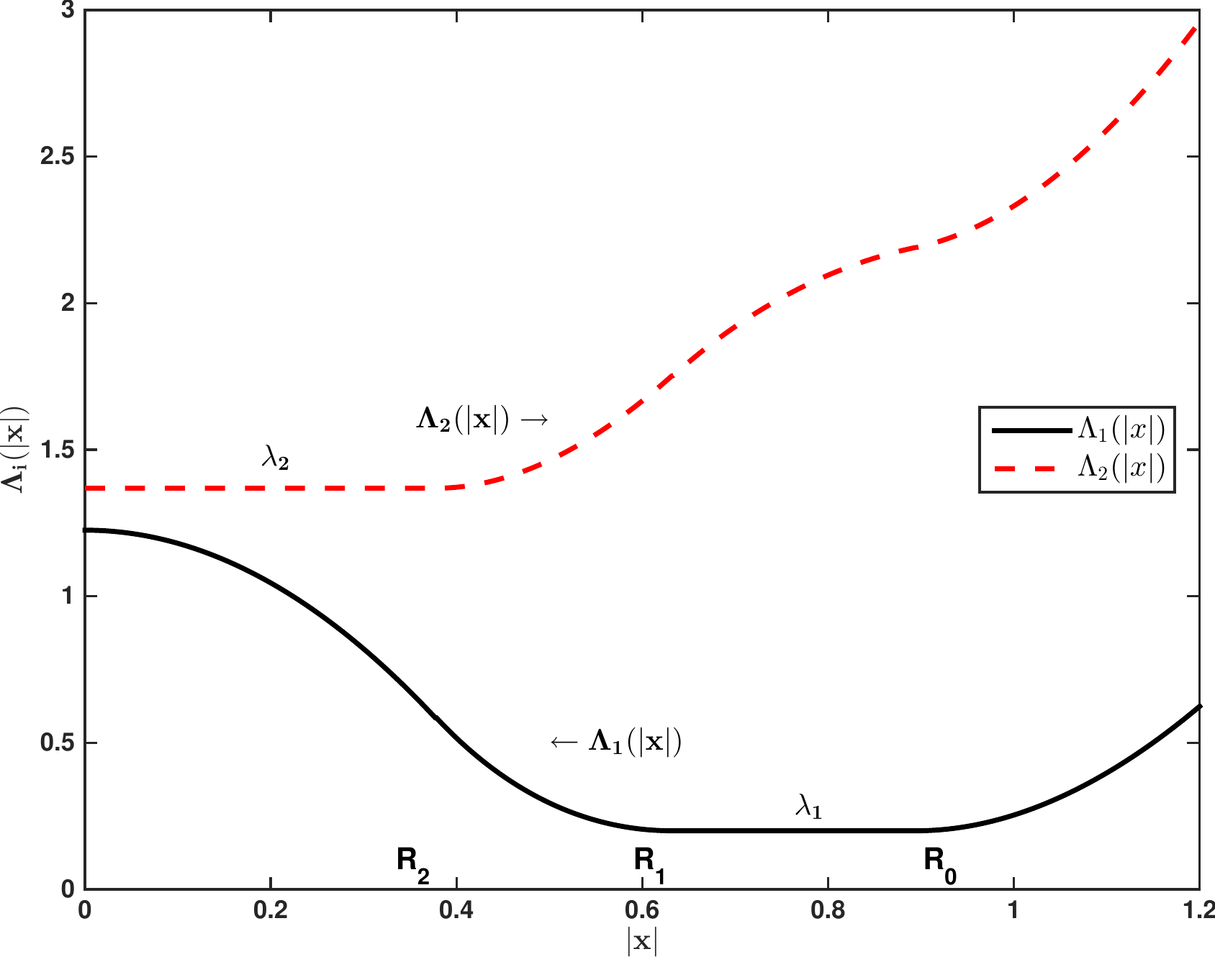} $~~~~$
\includegraphics[width=0.3\textwidth]{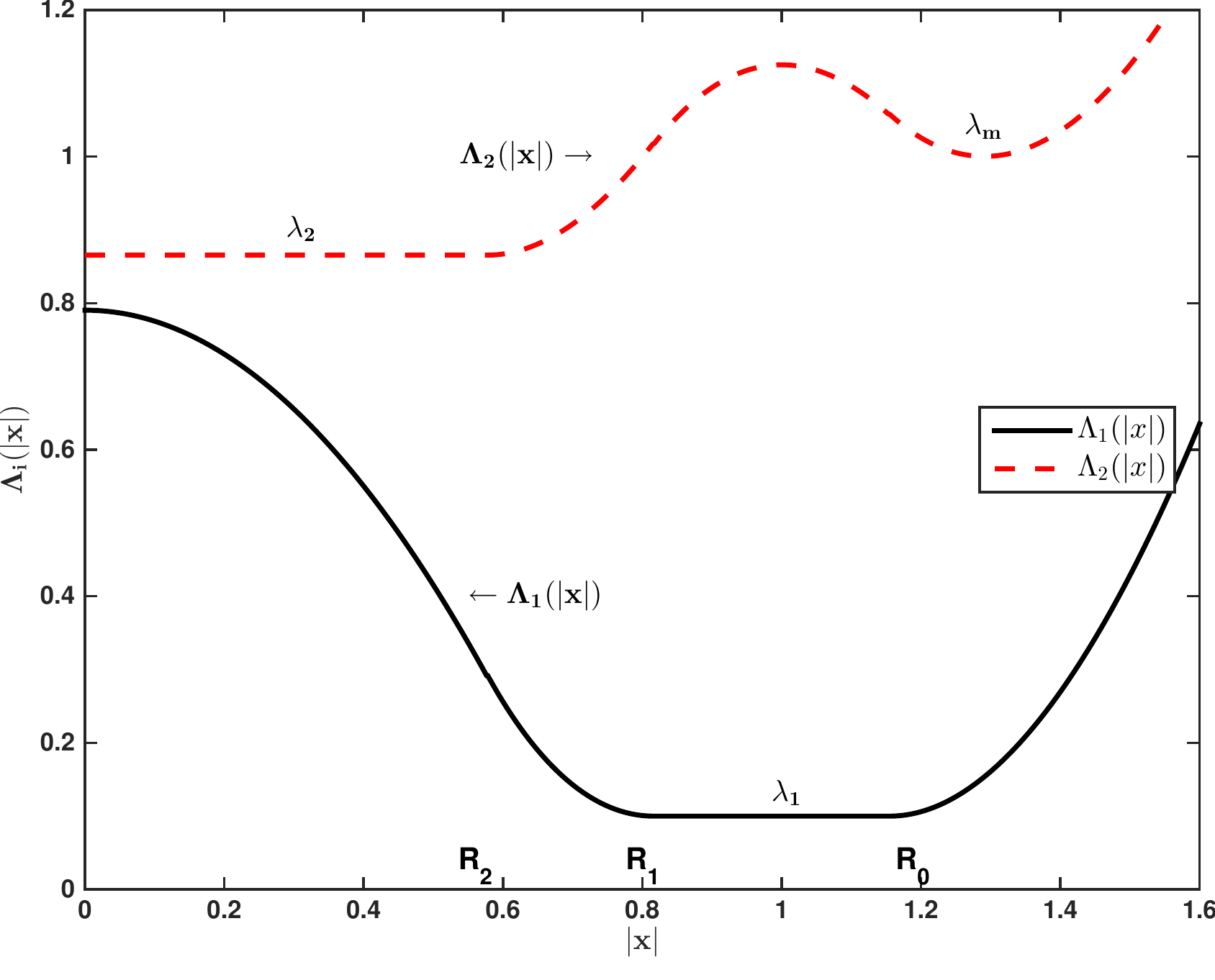} $~~~~$
\includegraphics[width=0.3\textwidth]{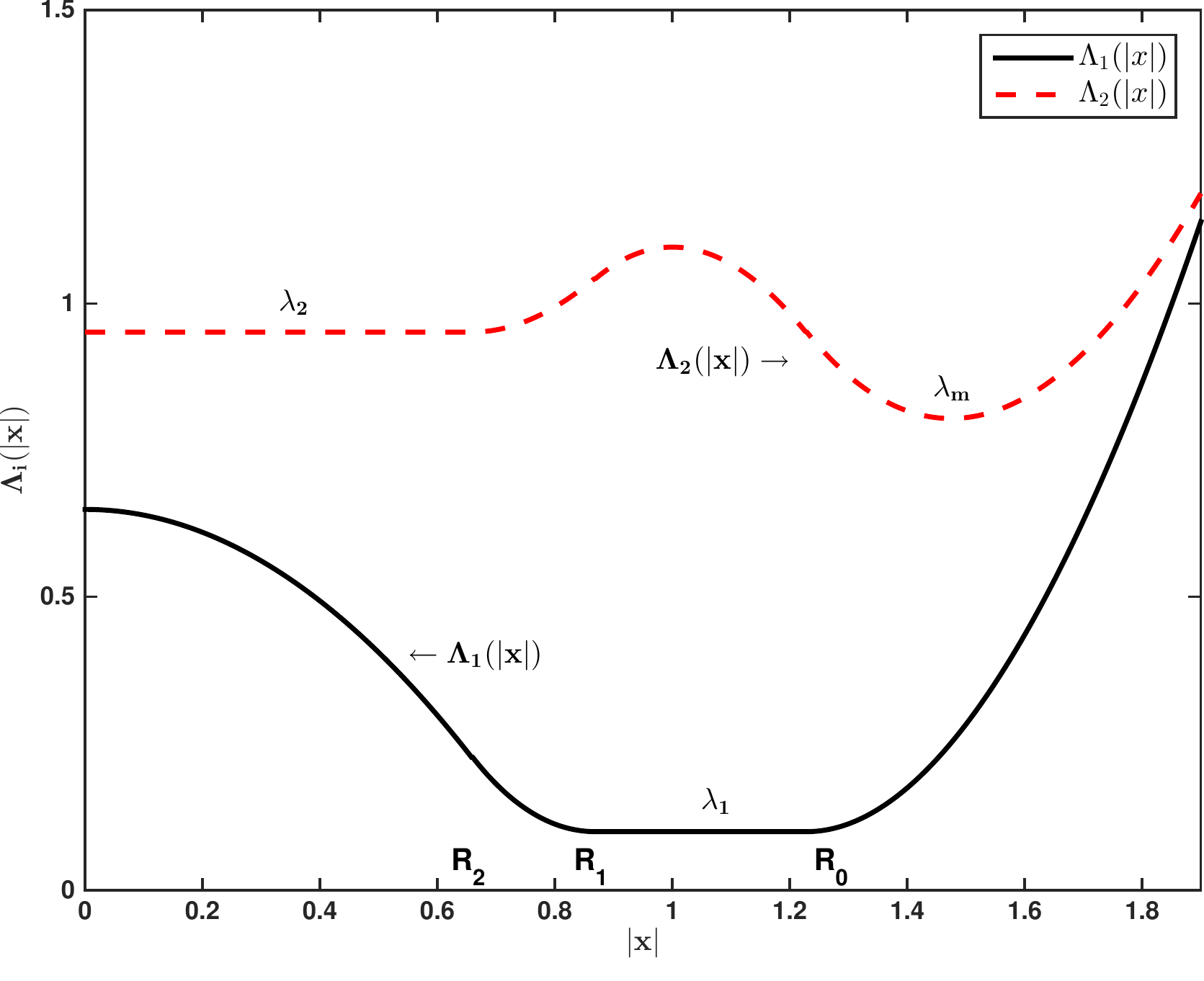} $~~~~$ 
\hspace{1.8cm}\textbf{(a)} \hspace{4.7cm} \textbf{(b)} \hspace{4.7cm} \textbf{(c)} 

\caption{Typical profiles of $\Lambda_1$ and $\Lambda_2$ corresponding to the target with lighter species inside. The profiles have been shifted vertically for a better visualization. (a) $(A,B)$ in $D_4 \cup D_5$, (b) $(A,B)$ in the subset of $D_3$ where \eqref{eqn:equilcomp} holds ($\lambda_m>\lambda_2$), (c) $(A,B)$ in $D_3$ where \eqref{eqn:equilcomp} fails ($\lambda_m<\lambda_2$) -- see also shaded areas in Figure \ref{fig:nonmin}. In (a) and (b) the target is a local minimizer with respect to perturbations of class $\CalB$. For (c) the equilibrium is not a minimizer, but it is a local minimizer with respect to perturbations that are also {\em local} in space -- see Remark \ref{rmk:target-localD3}.}
\label{fig:Lambdas}
\end{figure}

\begin{figure}[h]
\begin{center}
\includegraphics[width=0.45\textwidth]{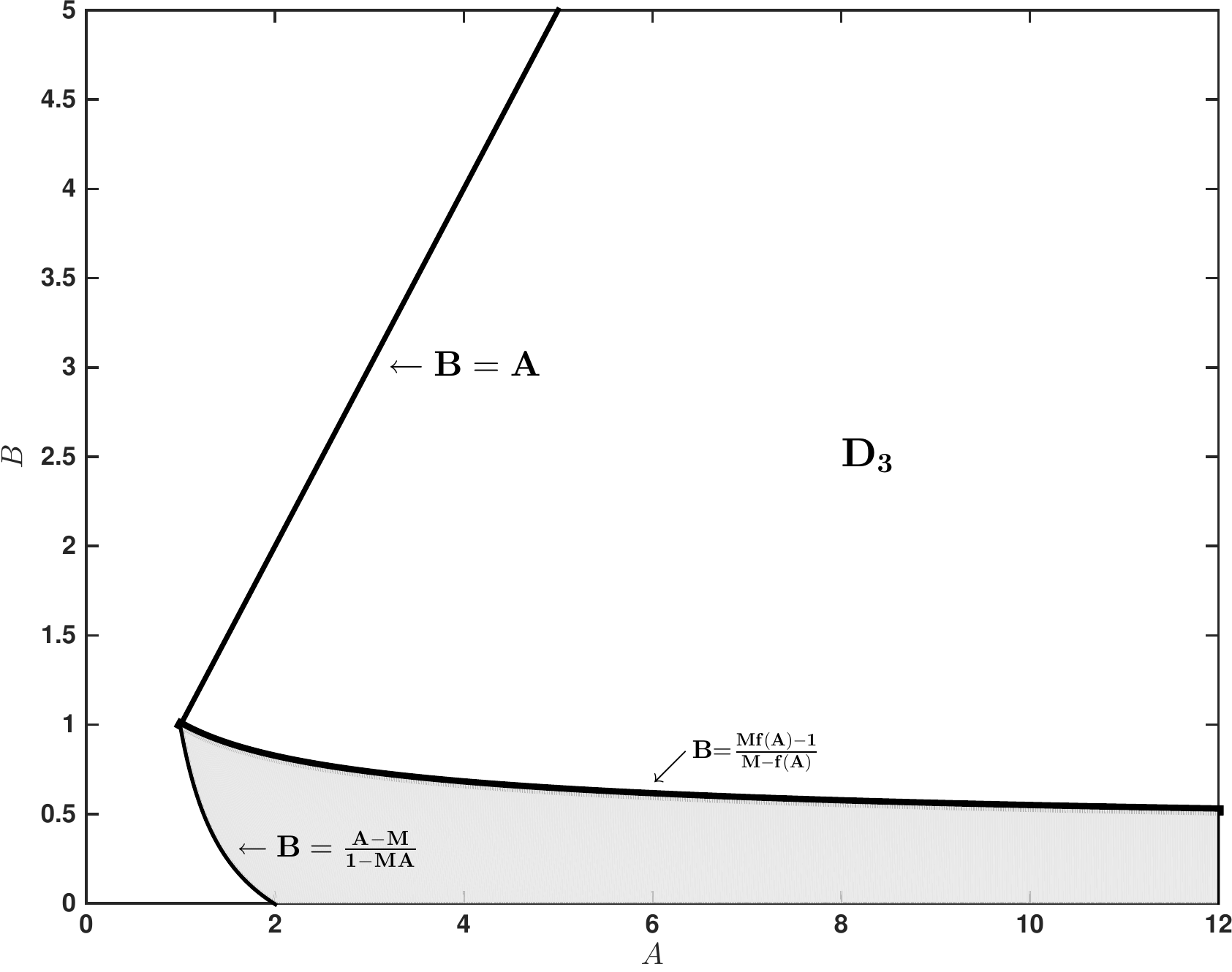} $~~~~$
\includegraphics[width=0.45\textwidth]{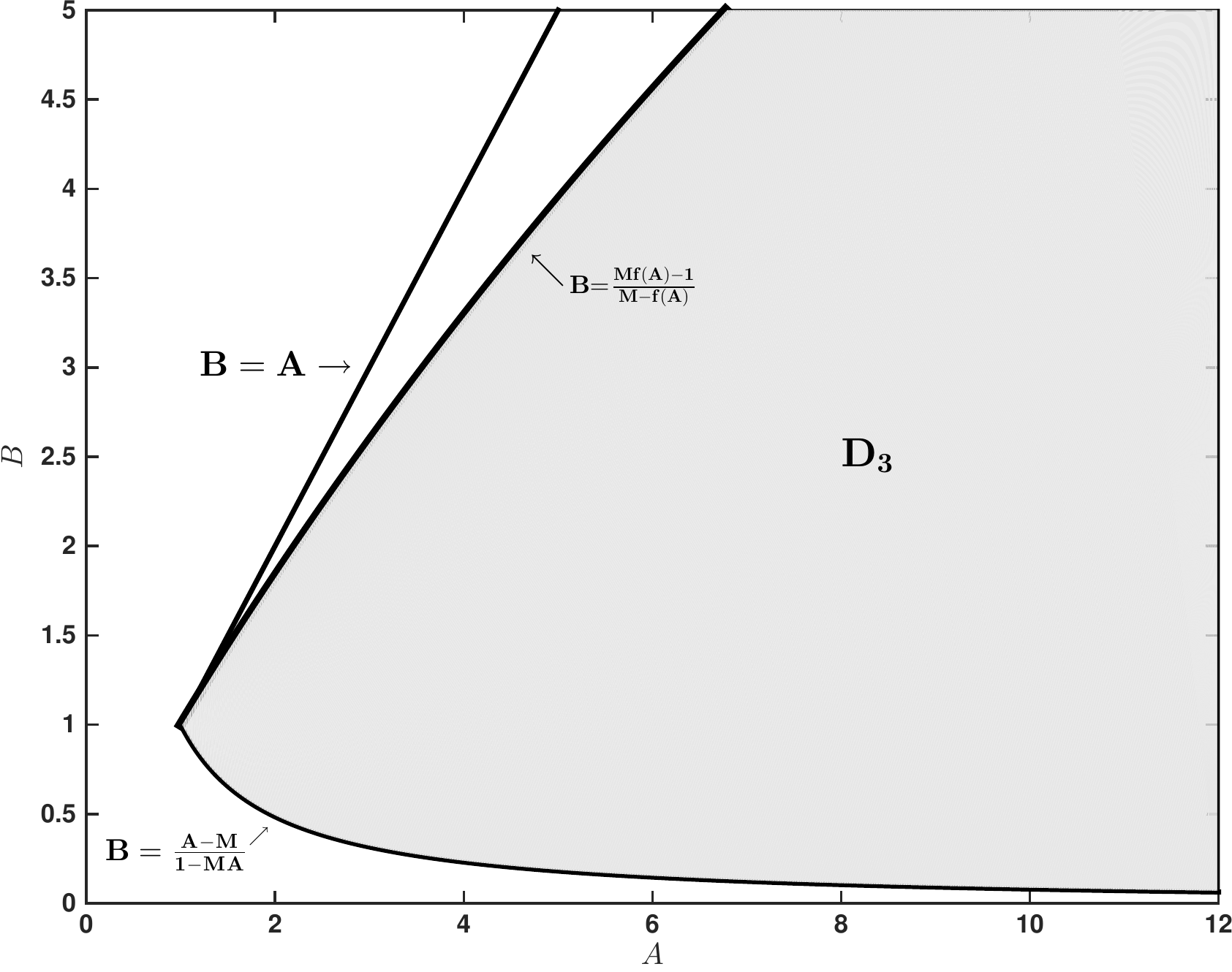} $~~~~$
\end{center}
\hspace{4cm}\textbf{(a)} \hspace{6.5cm} \textbf{(b)} 

\caption{Shaded areas represent the subsets of $D_3$ where \eqref{eqn:target-nonmin2} holds, for (a) $M=2$ and (b) $M=50$. For parameters in this region, the target equilibrium with the lighter species inside is not a local minimizer with respect to perturbations of class $\CalB$. Note that the regions where the target is not a minimizer tends to cover the entire $D_3$ as $M\to \infty$.}
\label{fig:nonmin}
\end{figure}
}

\begin{remark}
\label{rmk:target-localD3}
It has been noted in Section \ref{subsect:prelim-var} that the adjective ``local" refers to the size of the perturbations, as perturbations in class $\CalB$ are in fact nonlocal in space. While the target equilibrium is not a minimizer in this sense for all $(A,B)\in D_3$, it is however a minimizer with respect to perturbations that are also {\em local} in space, as \eqref{eqn:equilcomp} is indeed satisfied in a neighbourhood of the equilibrium -- see Figure \ref{fig:Lambdas}(c). The authors of \cite{BeTo2011} refer to such configurations as ``swarm minimizers".
\end{remark}

To conclude, we have derived that the target equilibrium in Figure \ref{fig:target} is a local minimizer with respect to perturbations of class $\CalB$  for all values of $(A,B)$ in $D_4$ and $D_5$. Moreover, it is also a local minimizer with respect to such perturbations in the subset of $D_3$ where \eqref{eqn:target-nonmin2} is false. We have not investigated perturbations of class $\CalA$ because of the difficulties pointed out in Section \ref{subsect:prelim-var}.

\begin{remark}[Global minimum]
\label{rmk:target-gmin}
By the considerations made in Section \ref{subsect:prelim-var}, the target equilibrium investigated here is a {\em global} minimizer for all $(A,B)$ in $D_5$ with $A<1$ (note that $B>1$ in $D_5$); i.e.~in the intersection of $D_5$ with the shaded region in Figure \ref{fig:stability}. This is a strong result which illustrates the relevance of this equilibrium. 
\end{remark}

Our numerical investigations of the particle system \eqref{eqn:part syst} suggest that the target equilibrium is \emph{only} stable in $D_4\cup D_5$. Some indication that this state is unstable in $D_3$ is given by the fact that the region \eqref{eqn:target-nonmin2} where the target is not a minimizer tends to cover the whole of $D_3$ as $M\to\infty$. Although we do not have the tools to show it, we conjecture that this state is a local minimizer with respect to perturbations of class $\CalA$ in $D_4\cup D_5$, but not in $D_3$.
It turns out that linear stability analysis supports this claim, as shown in the next section.


\subsubsection{Linear stability analysis}\label{sec:target1 lin}
We apply perturbations \eqref{eqn:perturb} to each of the three boundaries in Figure \ref{fig:target} and follow the lines of Section \ref{sect:prelim lin} to arrive at a linearized system \eqref{eqn:ODE eps Q}. To that aim, we have to find (linearized) expressions for the right-hand sides in
\begin{equation}\label{eqn:dpdt is v 3 times}
\frac{d}{dt}p_0(\theta_0)=v_1(p_0(\theta_0)), \quad \frac{d}{dt}p_1(\theta_0)=v_1(p_1(\theta_0)), \quad \text{and} \quad \frac{d}{dt}p_2(\theta_0)=v_2(p_2(\theta_0)).
\end{equation}
These right-hand sides involve integrals of the form \eqref{eqn:pert int} that are evaluated (up to higher-order terms) in Appendix \ref{app:pert integrals}. Suitably taking linear combinations of these integral evaluations, we obtain first-order approximations for the right-hand sides in \eqref{eqn:dpdt is v 3 times}. We use \eqref{eqn:dpdt}, divide by $R_j\exp(i\,\theta_0)$ and match sine and cosine terms on both sides of the equation, to obtain the system \eqref{eqn:ODE eps Q} for $m\in\N^+$. In particular, this system is independent of the choice of $\theta_0$. We use $\underline{\eps}:=(\eps_{0,N},\eps_{0,T},\eps_{1,N},\eps_{1,T},\eps_{2,N},\eps_{2,T})^T$ for the vector of all small perturbation parameters.

As anticipated in Section \ref{sect:prelim lin}, we indeed verified that the $\Ord(1)$ terms have zero contribution due to \eqref{eqn:radii} and \eqref{eqn:rho12-target}; indeed, $v_k(p_j(\theta_0))=0$ is a necessary condition for the target $\underline\eps=0$ to be a steady state.

The matrix $Q_m$ is given by
\begin{align}
Q_1:=\footnotesize
\left[\begin{array}{cccccc}
-a_s\,\pi\,\bar\rho_1+b_s\bar\rho_1\pi R_0^2 & 0 &  -a_s\bar\rho_1\pi \left(\dfrac{R_1}{R_0}\right)^3-b_s\bar\rho_1\pi \dfrac{R_1^3}{R_0} & 0 & \dfrac{M_2 a_c R_2}{R_0^3}+\dfrac{M_2 b_c R_2}{R_0} & 0\\
a_s\,\pi\,\bar\rho_1-b_s\,\bar\rho_1\,\pi\,R_0^2 &  0 & -a_s\,\pi\,\bar\rho_1\,\left(\dfrac{R_1}{R_0}\right)^3+b_s\,\bar\rho_1\,\pi\,\dfrac{R_1^3}{R_0} & 0 & \dfrac{M_2\,a_c\,R_2}{R_0^3}-\dfrac{M_2\,b_c\,R_2}{R_0} &  0 \\
-a_s\pi\bar\rho_1\dfrac{R_0}{R_1}+b_s\bar\rho_1\pi \dfrac{R_0^3}{R_1} & 0 & -a_s\,\pi\,\bar\rho_1-b_s \bar\rho_1\pi R_1^2 & 0 & \dfrac{M_2 a_c R_2}{R_1^3}+\dfrac{b_c M_2 R_2}{R_1} & 0\\
a_s\pi\bar\rho_1\dfrac{R_0}{R_1}-b_s\bar\rho_1\pi \dfrac{R_0^3}{R_1} & 0 & -a_s\pi\bar\rho_1+b_s\pi\bar\rho_1 R_1^2 & 0 & \dfrac{M_2 a_c R_2}{R_1^3} - \dfrac{M_2 b_c R_2}{R_1} & 0 \\
-a_c\pi\bar\rho_1\dfrac{R_0}{R_2}+b_c\pi\bar\rho_1\dfrac{R_0^3}{R_2} & 0 & a_c\pi\bar\rho_1\dfrac{R_1}{R_2}-b_c\bar\rho_1\pi \dfrac{R_1^3}{R_2} & 0 & -M_1\, b_c & 0 \\
a_c\pi\bar\rho_1\dfrac{R_0}{R_2}-b_c\bar\rho_1\pi\dfrac{R_0^3}{R_2} & 0 & -a_c\pi\bar\rho_1 \dfrac{R_1}{R_2}+b_c\bar\rho_1\pi\dfrac{R_1^3}{R_2} & 0 & M_1\, b_c & 0
\end{array}\right]
\label{eqn:Q1 target heavy outside}
\end{align}
for mode $m=1$, and for any other mode $m\geqs2$ by
\begin{align}
Q_m:=&\,   \footnotesize
\left[\begin{array}{cccccc}
-a_s\,\pi\,\bar\rho_1 & 0 &  -a_s\,\pi\,\bar\rho_1 \left(\dfrac{R_1}{R_0}\right)^{m+2} & 0 & a_c\,\pi\,\bar\rho_2\,\left(\dfrac{R_2}{R_0}\right)^{m+2} & 0\\
a_s\,\pi\,\bar\rho_1 &  0 & -a_s\,\pi\,\bar\rho_1\,\left(\dfrac{R_1}{R_0}\right)^{m+2} & 0 & a_c\,\pi\,\bar\rho_2\,\left(\dfrac{R_2}{R_0}\right)^{m+2} &  0 \\
-a_s\pi\bar\rho_1\left(\dfrac{R_1}{R_0}\right)^{m-2} & 0 & -a_s\,\pi\bar\rho_1 & 0 & a_c\,\pi\,\bar\rho_2\left(\dfrac{R_2}{R_1}\right)^{m+2}& 0\\
a_s\pi\bar\rho_1\left(\dfrac{R_1}{R_0}\right)^{m-2} & 0 & -a_s\pi\bar\rho_1& 0 & a_c\,\pi\,\bar\rho_2\left(\dfrac{R_2}{R_1}\right)^{m+2}  & 0 \\
-a_c\pi\bar\rho_1\left(\dfrac{R_2}{R_0}\right)^{m-2} & 0 & a_c\pi\bar\rho_1\left(\dfrac{R_2}{R_1}\right)^{m-2} & 0 & -a_s\pi\bar\rho_2 & 0 \\
a_c\pi\bar\rho_1\left(\dfrac{R_2}{R_0}\right)^{m-2} & 0 & -a_c\pi\bar\rho_1 \left(\dfrac{R_2}{R_1}\right)^{m-2} & 0 & a_s\pi\bar\rho_2 & 0
\end{array}\right].
\label{eqn:Qm target heavy outside}
\end{align}
The signs of the eigenvalues of $Q_m$ determine the stability of mode $m\geqs1$.\\
\\
To calculate $\det(Q_1-\lambda\,I)$ we first use a cofactor expansion with respect to the second, fourth and sixth row. We have that:
{\footnotesize
\begin{multline*}
 \det(Q_1-\lambda\,I)=\\ -\lambda^3 \det\left[\begin{array}{ccc}
b_sM_2(A-B) -\lambda & -b_sM_2(2M+A+B)\left(\frac{A}{M+A}\right)^{3/2}   &  b_sM_2(M(A+B)+2AB)\sqrt{\frac{C}{(M+A)^3}} \\
b_sM_2\sqrt{1+\frac{M}{A}}(A-B)  & -b_sM_2(M+B+A) -\lambda &  b_sM_2\sqrt{\frac{C}{A}}(M+2B)\\
b_sM_2\sqrt{\frac{M+A}{C}}M(B-A)  & b_sM_2MA\sqrt{\frac{A}{C}} &  -b_sM_2MB-\lambda  \\
\end{array}\right],
\end{multline*}}
where $C:=(M+B)/(1+MB)$. Subsequently, we compute and simplify the characteristic polynomial of the remaining $3\times3$ matrix and notice that its constant term is zero. We obtain:
\begin{align*}
\det(Q_1-\lambda\,I) =&\, \lambda^4 \left(\lambda^2   +b_sM_2(M+2B+MB)\,\lambda - (b_sM_2)^2M(M+1)(A-B)\dfrac{M+B}{M+A} \right).
\end{align*}
One can then show that the two nonzero eigenvalues are real and one of them is always negative. The other one is negative if and only if $B>A$. Recall that this target equilibrium only exists in $D_3\cup D_4\cup D_5$. Consequently, mode $m=1$ is unstable in the region $D_3$, and hence this target equilibrium itself is unstable in $D_3$.

\begin{figure}
\centering
\begin{tikzpicture}[>= latex]
\begin{scope}[xscale=2,yscale=2.8]
\pgfmathsetmacro\M {2};

\fill[fill=gray!5] (1,1) -- (4.25,1) -- (4.25,2.25) -- (2.25,2.25) -- cycle;
\fill[fill=gray!15] plot [domain=1:4.25,smooth,variable=\a]  ({\a},{(\M*\a^(2/3)-\a)/(\M*\a-\a^(2/3))}) -- (4.25,1) -- cycle;
\fill[fill=gray!30] plot [domain=1:4.25,smooth,variable=\a]  ({\a},{(\M*\a^(2/3)-\a)/(\M*\a-\a^(2/3))}) -- (4.25,0) -- plot [domain=4:1,smooth,variable=\a] ({\a},{(\M*\a^(2/4)-\a)/(\M*\a-\a^(2/4))}) -- cycle;
\fill[fill=gray!45] plot [domain=1:4,smooth,variable=\a] ({\a},{(\M*\a^(2/4)-\a)/(\M*\a-\a^(2/4))}) -- plot [domain=0:1,smooth,variable=\b] ({(\b+\M)/(1+\M*\b)},{\b}) -- cycle;

\draw[->,line width=2] (0,0) -- (0,2.25);
\draw[->,line width=2] (0,0) -- (4.5,0);

\draw[line width=1] (\M,0.04) -- (\M,-0.04) node[below] {$M$};
\draw[line width=1] (1/\M,0.04) -- (1/\M,-0.04) node[below] {$1/M$};
\draw[line width=1] (1,0.04) -- (1,-0.04) node[below] {$1$};
\draw[line width=1] (0.04,1) -- (-0.04,1) node[left] {$1$};

\draw[domain=0:2.25,smooth,variable=\b,line width=0.25,dashed]  plot ({(1+\M*\b)/(\b+\M)},{\b});
\draw[domain=0:2.25,smooth,variable=\b,line width=2]  plot ({(\b+\M)/(1+\M*\b)},{\b});

\begin{scope} 
\tikzstyle{every node} = [draw,rectangle, fill=gray!5, line width = 0.25]
\draw[line width=0.25, dashed, gray] (0,0) -- (1,1);
\draw[line width=2] (1,1) -- (2.25,2.25);
\draw [line width=0.75,->] (2.4,2) node[right] {$m=1$} arc (-90:-130:0.4);
\draw[line width=0.25, dashed, gray] (0,1) -- (1,1);
\draw[line width=2] (1,1) -- (4.25,1);
\draw [line width=0.75,->] (4.2,1.3) node[above] {$m=2$}  arc (0:-40:0.4);
\draw[domain=1:4.25,smooth,variable=\a,line width=2]  plot ({\a},{(\M*\a^(2/3)-\a)/(\M*\a-\a^(2/3))});
\draw [line width=0.75,->] (4.2,0.47) node {$m=3$} arc (0:-40:0.4);
\draw[domain=1:{\M^2},smooth,variable=\a,line width=2]  plot ({\a},{(\M*\a^(2/4)-\a)/(\M*\a-\a^(2/4))});
\draw [line width=0.75,->] (2.5,-0.2) node {$m=4$} arc (180:125:0.4);
\end{scope}
\draw[line width=1] (\M^2,0.04) -- (\M^2,-0.04) node[below] {$M^2$};
\end{scope}
\end{tikzpicture}
\caption{The boundaries of the regions $U_m$ for $m=1,2,3,4$. Note that $U_4\subset U_3\subset U_2 \subset U_1$.}%
\label{fig:stab regions target heavy outside}%
\end{figure}
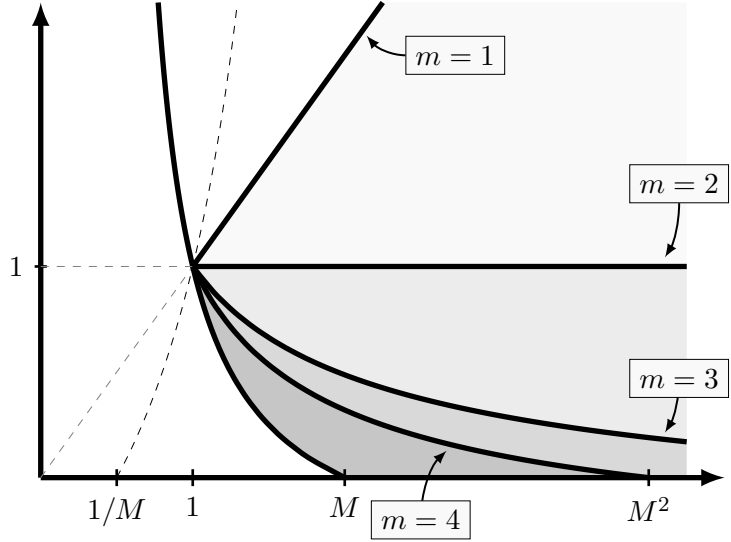

To further asses the stability in $D_4\cup D_5$, we are required to investigate the eigenvalues for all higher-order modes. For general $m\geqs2$ we have
\begin{align}\footnotesize
\nonumber \det(Q_m-\lambda\,I)=&\,
-\lambda^3\det\left[\begin{array}{ccc}
-a_s\,\pi\,\bar\rho_1-\lambda &  -a_s\,\pi\,\bar\rho_1 \left(\dfrac{R_1}{R_0}\right)^{m+2} & a_c\,\pi\,\bar\rho_2\,\left(\dfrac{R_2}{R_0}\right)^{m+2} \\
-a_s\pi\bar\rho_1\left(\dfrac{R_1}{R_0}\right)^{m-2}  & -a_s\,\pi\bar\rho_1-\lambda &  a_c\,\pi\,\bar\rho_2\left(\dfrac{R_2}{R_1}\right)^{m+2}\\
-a_c\pi\bar\rho_1\left(\dfrac{R_2}{R_0}\right)^{m-2}  & a_c\pi\bar\rho_1\left(\dfrac{R_2}{R_1}\right)^{m-2}  & -a_s\pi\bar\rho_2  -\lambda
\end{array}\right],
\end{align}
and the characteristic polynomial is
\begin{multline}\label{eqn:char pol m target heavy outside}
\lambda^3\Bigg[ \lambda^3 + \pi\,a_s\,\bar\rho_1\left(2+\dfrac1C\right)\,\lambda^2 +(\pi\,a_s\,\bar\rho_1)^2\left(\dfrac2C + \left(1-A\left(\dfrac CA\right)^{m-1}\right)\left(1-\left(\dfrac{A}{M+A}\right)^m\right)\right)\lambda\\ +(\pi\,a_s\,\bar\rho_1)^3\dfrac1C\left(1-A^2\left(\dfrac CA\right)^{m}\right)\left(1-\left(\dfrac{A}{M+A}\right)^m\right) \Bigg].
\end{multline}
One can show that in $D_3\cup D_4 \cup D_5$ it always holds that $1/M<C<A$. Let the phrase \textit{nontrivial eigenvalues} denote the roots of the cubic polynomial in square brackets in \eqref{eqn:char pol m target heavy outside}. The constant
\begin{equation}\label{eqn:const target}
-(\pi\,a_s\,\bar\rho_1)^3\dfrac1C\left(1-A^2\left(\dfrac CA\right)^{m}\right)\left(1-\left(\dfrac{A}{M+A}\right)^m\right)
\end{equation}
is the product of these three nontrivial eigenvalues. This constant is positive if and only if $C>A^{1-2/m}$. Hence, restricting ourselves to $ D_3\cup D_4\cup D_5$, we know that the constant is positive if and only if $(A,B)\in U_m$, with
\[ U_m := \left\{(A,B)\in D_3\cup D_4\cup D_5 \; : \; 1<A<M^{\frac{m}{m-2}} \text{ and } 0<B< \dfrac{MA^{\frac2m}-A}{MA-A^{\frac2m}} \right\}. \]
Here, it is understood that $M^{\frac{m}{m-2}}=\infty$ for $m=2$, and $U_2$ is $\{1<A<\infty, 0<B<1\}\cap (D_3\cup D_4\cup D_5)$.\\
\\
We repeat that the constant \eqref{eqn:const target} is the product of the three nontrivial eigenvalues, and it is positive in $U_m$. Hence, at least one of the eigenvalues must be (real and) positive.\footnote{This argument holds when all three eigenvalues are real, and also when the other two eigenvalues are a pair of complex conjugates.} Thus, we know that mode $m\geqs2$ is unstable in $U_m$.\\
\\
In view of our stability analysis for mode $m=1$, for consistency we define $U_1:=D_3$. If $(A,B)\in U_1$, then $B<A$, while in $U_2$ we have that $B<1$. Since furthermore $A>1$ in $U_2$, it holds that $U_2\subset U_1$. Next, let $n>m\geqs2$ and let $(A,B)\in U_n$. Then
\begin{align*}
A <&\, M^{\frac{n}{n-2}} < M^{\frac{m}{m-2}}, \text{ and }\\
B<&\, \dfrac{MA^{\frac2n}-A}{MA-A^{\frac2n}} < \dfrac{MA^{\frac2m}-A}{MA-A^{\frac2m}},
\end{align*}
so $(A,B)\in U_m$, and thus $U_n\subset U_m$. See Figure \ref{fig:stab regions target heavy outside} for an indication of the boundaries of the regions $U_m$ and the way in which $U_{m+1}$ is contained in $U_m$.\\

Now investigate the stability of mode $m\geqs2$ outside $U_m$. Write $\mu=\lambda/(\pi\,a_s\,\bar\rho_1)$ and for $m\geqslant2$ note that the nontrivial roots of \eqref{eqn:char pol m target heavy outside} correspond to the roots of the polynomial
\begin{multline}\label{eqn:char pol mu}
P(\mu):=\mu^3 + \left(2+\dfrac1C\right)\,\mu^2 +\left(\dfrac2C + \left(1-A\left(\dfrac CA\right)^{m-1}\right)\left(1-\left(\dfrac{A}{M+A}\right)^m\right)\right)\,\mu\\ +\dfrac1C\left(1-A^2\left(\dfrac CA\right)^{m}\right)\left(1-\left(\dfrac{A}{M+A}\right)^m\right).
\end{multline}
We observe that $\lim_{\mu\to-\infty}P(\mu)=-\infty$, and that
\begin{equation*}
P(-1)=\left(1-\dfrac{1}{C}\right)\left( \dfrac{A}{M+A} \right)^m, \quad \text{and} \quad P\left(-\dfrac{1}{C}  \right) =  (1-C)\left( \dfrac{C}{A} \right)^{m-2}\left(1-\left(\dfrac{A}{M+A}\right)^m\right).
\end{equation*}
Furthermore $P(0)>0$ holds; this follows from our previous observation that the expression in \eqref{eqn:const target} is negative in $S_m:=(D_3\cup D_4\cup D_5)\setminus U_m$. For $m\geqslant 2$ consider the cases:
\begin{description}
\item[$C<1$:] then $-\infty<-1/C<-1<0$ and furthermore $P(-1/C)>0$, and $P(-1)<0$, and $P(0)>0$. Since also $\lim_{\mu\to-\infty}P(\mu)=-\infty$, the intermediate value theorem implies that the polynomial $P$ has three real roots, and they are all negative.
\item[$C>1$:] then $-\infty<-1<-1/C<0$ and furthermore $P(-1)>0$, and $P(-1/C)<0$, and $P(0)>0$. It also holds that $\lim_{\mu\to-\infty}P(\mu)=-\infty$, and thus $P$ has three real, negative roots. Note that $C>1$ can not occur in $S_m$ for $m=2$.
\item[$C=1$:] this case is only relevant for $m\geqslant3$, since $C=1$ is the upper boundary of $S_2$. Note that $C=1$ is equivalent to $B=1$. For any $m\geqslant 3$, if $B=1$ in $S_m$, then $A>1$.\\
If $C=1$ then $P(-1)=0$. After isolating the factor $(\mu+1)$ in $P(\mu)$, the other two roots of $P$ can be found explicitly. They are real and both negative, which can be shown easily using that $m\geqslant3$ and $A>1$. 
\end{description}

Hence, all (nontrivial) eigenvalues of $Q_m$ are real and negative in $S_m$ for any $m\geqs2$. Thus mode $m\geqs2$ is stable in $S_m$ and moreover we have that $S_m\subset S_n$ if $n>m\geqs1$. We showed before that mode $m=1$ is stable if and only if $(A,B)\in D_4\cup D_5 =: S_1$. Hence we have showed the stability of the target (lighter species inside) for all modes $m\geqs1$ provided that $(A,B)\in S_1=D_4\cup D_5$.


\subsection{Heavier species inside}
\label{subsect:t-hinside}

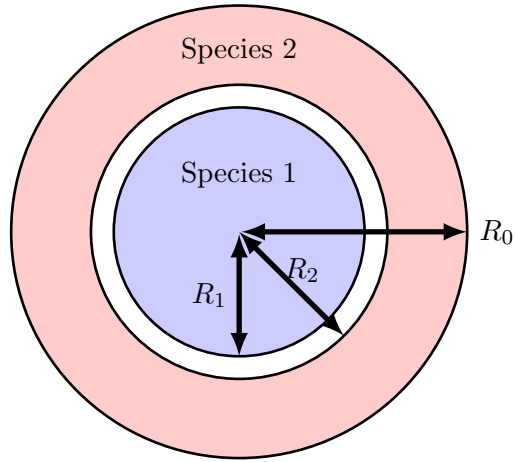
\begin{figure}[h]
\centering
\begin{tikzpicture}[>= latex]
\begin{scope}[scale=0.3] 

\draw[fill=red!20
,even odd rule,line width=1]  (0,0) circle (10)
                                   (0,0) circle (6.5);
\draw[fill=blue!20
,line width=1]  (0,0) circle (5.5);            
                
\draw[<->,line width=2] (0.1,0) -- (10,0) node[pos=1,right] {$R_0$};
\draw[<->,line width=2] (0,0) -- +(-45:6.5) node[pos=.6,above] {$R_2$};
\draw[<->,line width=2] (0,-0.1) -- (0,-5.5) node[pos=.5,left] {$R_1$};
\draw (0,2.5) node {Species 1};
\draw (0,8) node {Species 2};

\end{scope}
\end{tikzpicture}
\caption{Target equilibrium, with species $1$ supported on a disk of radius $\Rone$, and (the lighter) species $2$ supported on the annular region $\Rtwo<|x|<\Rzero$.}
\label{fig:target heavy inside}
\end{figure}

In this equilibrium state the (heavier) species $1$ is supported in the disk $|x|<\Rone$, and the (lighter) species $2$ is supported on the annulus  $\Rtwo<|x|<\Rthree$, with $\Rone \leq \Rtwo \leq \Rzero$ -- see Figure \ref{fig:target heavy inside}.  

Calculations lead to
\begin{equation}
\label{eqn:radii2}
\Rone^2 = \frac{\as M_1}{\bs M_1 + \bc M_2}, \qquad \Rtwo^2 =  \frac{\ac M_1}{\bc M_1 + \bs M_2}, \qquad \Rthree^2 = \frac{\ac M_1 + \as M_2}{\bc M_1 + \bs M_2}.
\end{equation}
The equilibrium densities are given by \eqref{eqn:rho12-target} (cf. \eqref{eqn:equilibria}). For consistency with the solution ansatz, $(\A,\B)$ has to lie in $\D_2\cup \D_3 \cup \D_4$.

\subsubsection{Variational approach}
\label{subsubsect:t-var-hinside}
The calculations of $\Lambda_1'$ and $\Lambda_2'$ mirror the ones in Section \ref{subsubsect:t-var-linside}. The results are:
\begin{equation}
\label{eqn:gradL1-target2}
\frac{\Lambda_1'(|x|)}{|x|}= \begin{cases}
\pi \as \barrho_1 (1-\Rone^2/|x|^2) & \text{ if } \Rone<|x|<\Rtwo \\ 
\pi \as \barrho_1 (1-\Rone^2/|x|^2) - \pi \ac \barrho_2 (1-\Rtwo^2/|x|^2) &  \text{ if } \Rtwo<|x|<\Rthree \\
\bs M_1 + \bc M_2 - (\as M_1 + \ac M_2)/|x|^2  &  \text{ if } |x|>\Rthree,
\end{cases}
\end{equation}
and
\begin{equation}
\label{eqn:gradL2-target2}
\frac{\Lambda_2'(|x|)}{|x|}= \begin{cases}
\bc M_1 + \bs M_2 - \frac{\ac}{\as} (\bs M_1 + \bc M_2)& \text{ if } |x|<\Rone \\ 
\bc M_1 + \bs M_2 - \ac M_1/|x|^2&  \text{ if } \Rone<|x|<\Rtwo \\
\bc M_1 + \bs M_2 - (\ac M_1+ \as M_2)/|x|^2  &  \text{ if } |x|>\Rthree.
\end{cases}
\end{equation}
Also note that $\Lambda_1'$ and $\Lambda_2'$ vanish on the respective supports of $\barrho_1$ and $\barrho_2$: $\{|x|<\Rone\}$ and $\{\Rtwo<|x|<\Rthree\}$.

It is immediate to check that $\Lambda_2$ (corresponding to the lighter species) satisfies \eqref{eqn:equilcomp}. Indeed,  $\Lambda_2'$ in $|x|<\Rone$ can be written as:
\[
\Lambda_2'(|x|) =  \bs M_2 (B(M-A)+ 1- AM) |x|, \qquad \text{ in } |x|<\Rone.
\]
For all $(A,B)$ in $\D_2 \cup \D_3 \cup \D_4$ (see Figure \ref{fig:stability}), the expression above is negative. 

In $\Rone<|x|<\Rtwo$, $\Lambda_2(|x|)$ is also decreasing, as can be seen from the simple estimate below: 
\[
\bc M_1 + \bs M_2 - \frac{\ac M_1}{|x|^2} < \underbrace{\bc M_1 + \bs M_2 - \frac{\ac M_1}{{\Rtwo}^2}}_{=0 \text{ by } \eqref{eqn:radii2}}.
\]
Finally, in $|x|>\Rthree$, $\Lambda_2(|x|)$ is increasing, as
\[
\bc M_1 + \bs M_2 - \frac{\ac M_1+ \as M_2}{|x|^2} > \underbrace{\bc M_1 + \bs M_2 - \frac{\ac M_1+ \as M_2}{\Rthree^2}}_{=0 \text{ by } \eqref{eqn:radii2}}.
\]
Hence, for all $(A,B)$ in the relevant region $\D_2 \cup \D_3 \cup \D_4$, $\Lambda_2$ satisfies \eqref{eqn:equilcomp}.

As in Section \ref{subsubsect:t-var-linside}, the calculations for the heavier species are slightly more involved and do not always lead to the minimization condition  \eqref{eqn:equilcomp}. From \eqref{eqn:gradL1-target2}, one concludes easily that $\Lambda_1'(|x|)>0$ in $\Rone<|x|<\Rtwo$. However, $\Lambda_1'$ can become zero in $\Rtwo<|x|<\Rthree$, and hence $\Lambda_1$ can decrease in this region. The zero of $\Lambda_1'$ occurs at
\begin{equation}
\label{eqn:zero12}
|x|^2 = \frac{\pi \as \barrho_1 \Rone^2 - \pi \ac \barrho_2 \Rtwo^2}{\pi \as \barrho_1-\pi \ac \barrho_2},
\end{equation}
where for consistency, the expression in \eqref{eqn:zero12} needs to be positive and also, it has to lie in the annular region $\{ \Rtwo^2<|x|^2<\Rthree^2\}$. The calculations of the numerator and denominator in \eqref{eqn:zero12} lead to:
\[
\pi \as \barrho_1 \Rone^2 - \pi \ac \barrho_2 \Rtwo^2 = M_1 \as (1-A^2),
\]
and
\[
\pi \as \barrho_1-\pi \ac \barrho_2 = M_2 \bs (M+B - A(BM+1)),
\]
respectively. 

By the above, it is immediate to show that the denominator \eqref{eqn:zero12} is positive in $\D_2$ and negative in $\D_3 \cup \D_4$. Consider first the case when it is positive, i.e., $(A,B) \in \D_2$. Since the denominator is positive, the numerator has to be positive as well (i.e., $A<1$). It is then a simple exercise to show that the zero of $\Lambda_1'$ from \eqref{eqn:zero12} does not lie in the relevant region $\Rtwo<|x|<\Rthree$. Hence, for $(A,B) \in \D_2$, $\Lambda_1'$ does not change sign in the annular region and remains positive throughout (as a check we found indeed that $\Lambda_1'(\Rthree) >0$).

For $(A,B) \in \D_3 \cup \D_4$, where the denominator in \eqref{eqn:zero12} is negative, the numerator is also negative (as $A>1$ there). Also, the location of the zero given by \eqref{eqn:zero12} is larger than $\Rtwo$ for all $(A,B) \in \D_3 \cup \D_4$. We require that it is also less than $\Rthree$ and we arrive after some elementary calculations to the following condition:
\[
(M_1^2-M_2^2) (A-B) <0.
\]
Since $M_1>M_2$, to have a zero of $\Lambda_1'$ in the annular region one needs $A<B$, i.e.,  $(A,B) \in \D_4$. Otherwise, for $(A,B) \in \D_3$, $\Lambda_1'$ does not change sign and remains positive in the annular region.

Finally, in $|x|>\Rthree$, it can be shown that $\Lambda_1$ remains strictly increasing for $(A,B) \in \D_2 \cup \D_3$, where $A>B$. Combined with the findings above (including the calculations for $\Lambda_2$), we infer that \eqref{eqn:equilcomp} holds for all parameter values $(A,B) \in \D_2 \cup \D_3$. 

On the other hand, for $(A,B) \in \D_4$, $\Lambda_1$ changes monotonicity in $|x|>\Rthree$, at 
\begin{equation}
|x|^2 = \frac{\as M_1 + \ac M_2}{\bs M_1 + \bc M_2} > \Rthree^2 \qquad \text { for } (A,B) \in \D_4.
\end{equation}

If $\Lambda_1$ evaluated at the minimum point above drops below $\lambda_1$, then the condition for $\Lambda_1$ in \eqref{eqn:equilcomp} fails, and this target solution is not a minimizer. We do not present these calculations, we only list the end result. Following calculations similar to the other target equilibrium (see  \eqref{eqn:target-nonmin1}), we find that the target with the heavier species inside is \emph{not} a local minimizer provided
\begin{equation}
\label{eqn:itarget-nonmin}
\frac{B+M}{BM+1} < \frac{A+M}{AM+1} \cdot \left(1+\frac{A}{M}\right)^{\frac{M}{A}} \cdot \left( 1+ \frac{1}{AM}\right)^{-AM}.
\end{equation}
The inequality \eqref{eqn:itarget-nonmin}, which can be rearranged to be explicit in $B$, describes a subset of $D_4$ which grows with the mass ratio $M$. For the rest of $(A,B) \in D_4$, $\lambda_m>\lambda_1$ and hence \eqref{eqn:equilcomp} holds.

In summary, the target with the heavier species inside is a minimizer with respect to class $\CalB$ perturbations for all  $(A,B) \in \D_2 \cup \D_3$ and for certain $(A,B) \in \D_4$ for which \eqref{eqn:itarget-nonmin} is violated. 

Our numerical investigations of the particle system \eqref{eqn:part syst} do \emph{not} agree with the conclusion of the variational approach: we never observe the target of Figure \ref{fig:target heavy inside} as a numerical steady state, and hence we conjecture that it is never stable. In particular, we conjecture that this state is unstable with respect to class $\CalA$ perturbations in $D_2$, $D_3$ and (the indicated part of) $D_4$, although we do not have the means to prove this.
But, as we show in the next section, the linear stability analysis is in agreement with this claim.

\subsubsection{Linear stability analysis}\label{subsect:lin t-hinside}
For this steady state, we can obtain the corresponding matrices $Q_m$ (for $m\geqs1$) directly from \eqref{eqn:Q1 target heavy outside} and \eqref{eqn:Qm target heavy outside} by writing them fully in terms of $a_s$, $a_c$, $b_s$, $b_c$, $M_1$ and $M_2$ and subsequently interchanging $M_1$ and $M_2$. These matrices  correspond to the system of ODE's \eqref{eqn:ODE eps Q} where in fact $\underline{\eps}:=(\eps_{0,N},\eps_{0,T},\eps_{2,N},\eps_{2,T},\eps_{1,N},\eps_{1,T})^T$ is reordered. Alternatively, we could have derived these matrices starting from the building blocks in Appendix \ref{app:pert integrals}, analogously to what we did in Section \ref{sec:target1 lin}.
After some further steps, we obtain for mode $m=1$:
{\footnotesize
\begin{multline*}
 \det(Q_1-\lambda\,I)=\\ -\lambda^3 \det\left[\begin{array}{ccc}
b_sM_2M(A-B) -\lambda & -b_sM_2(2+M(A+B))\left(\frac{AM}{G}\right)^{3/2}   & 
b_sM_2(M(A+B)+2M^2AB)\sqrt{\frac{M}{C\,G^3}} \\
b_sM_2\sqrt{\frac{G}{AM}}M(A-B)  & -b_sM_2(1+M(A+B)) -\lambda & b_sM_2\sqrt{\frac{1}{AC}}(1+2MB)\\ 
b_sM_2\sqrt{\frac{CG}{M}}(B-A)  & b_sM_2A\sqrt{AC} &  -b_sM_2B-\lambda  
\end{array}\right],
\end{multline*}
}
with $C:=(M+B)/(1+MB)$ and $G:=1+MA$. The characteristic polynomial is
\begin{align*}
\det(Q_1-\lambda\,I) =&\, \lambda^4 \left(\lambda^2  +b_sM_2(1+B+2MB)\,\lambda - (b_sM_2)^2(M+1)(A-B)\dfrac{1+MB}{1+MA} \right),
\end{align*}
and the nonzero roots can be calculated explicitly. These two eigenvalues are real. One of them is always negative and the other one is negative if and only if $B>A$; ti.e.~for $(A,B)\in D_4$. Hence, this state is unstable in $D_2\cup D_3$ where one eigenvalue is positive.

For mode $m=2$ we have:
\begin{align}\footnotesize
\nonumber \det(Q_2-\lambda\,I)=&\,
-\lambda^3\det\left[\begin{array}{ccc}
-a_s\,\pi\,\bar\rho_2-\lambda &  -a_s\,\pi\,\bar\rho_2 \left(\dfrac{R_2}{R_0}\right)^{4} & a_c\,\pi\,\bar\rho_1\,\left(\dfrac{R_1}{R_0}\right)^{4} \\
-a_s\pi\bar\rho_2 & -a_s\,\pi\bar\rho_2-\lambda &  a_c\,\pi\,\bar\rho_1\left(\dfrac{R_1}{R_2}\right)^{4}\\
-a_c\pi\bar\rho_2  & a_c\pi\bar\rho_2  & -a_s\pi\bar\rho_1  -\lambda
\end{array}\right],
\end{align}
and the characteristic polynomial is
\begin{multline}
\lambda^3\Bigg( \lambda^3 + \pi\,a_s\,\bar\rho_2\left(2+C\right)\,\lambda^2 +(\pi\,a_s\,\bar\rho_2)^2\left(2C + \left(1-\dfrac1C\right)\left(1-\left(\dfrac{AM}{1+MA}\right)^2\right)\right)\,\lambda\\ -(\pi\,a_s\,\bar\rho_2)^3\dfrac1C\left(1-C^2\right)\left(1-\left(\dfrac{AM}{1+MA}\right)^2\right) \Bigg).
\end{multline}
Compared to \eqref{eqn:char pol m target heavy outside}, there is in particular a change in sign for the constant term inside the brackets. Hence, since the constant
\[ (\pi\,a_s\,\bar\rho_2)^3\dfrac1C\left(1-C^2\right)\left(1-\left(\dfrac{AM}{1+MA}\right)^2\right)  \]
equals the product of the three nontrivial eigenvalues, we observe that there is at least one positive eigenvalue if this constant is positive.\footnote{This argument holds when all eigenvalues are real and also when two eigenvalues are complex conjugates. Here, it is thus not even necessary to verify whether the roots of the characteristic polynomial are real.} This is the case if $C<1$, or equivalently if $B>1$. Consequently, mode $m=2$ is unstable if $B>1$.\\
\\
We previously concluded that mode $m=1$ is unstable in $D_2\cup D_3$. It follows that modes 1 and 2 are never stable simultaneously, and hence this target steady state must be unstable for any choice of $(A,B)\in D_2\cup D_3 \cup D_4$. 


\subsection{Numerical illustration of the unstable modes}\label{sec:num unstable modes}

\subsubsection{Target: lighter species inside}
In Section \ref{subsect:t-linside} we derived that the target state (with the lighter species inside) is stable in $D_4\cup D_5$. Mode 1 is unstable in region $D_3$, while mode 2 is unstable for $B<1$. The instability regions for the higher-order modes are such that for mode $m$ this region is a subset of the instability region for mode $m-1$. See Section \ref{sec:target1 lin} for the full details.

Here, we further illustrate the instability. First we run a particle system of 200 particles with $M=2$ and $(A,B)=(3,3.5)$. The system approaches the target steady state shown in Figure \ref{fig:num st st} at the top. Next we choose two pairs of parameter values such that $1<B<A$, and $B<1$, respectively. Specifically, we take $(A,B)=(3,2)$ and $(A,B)=(3,0.75)$. For the latter parameter pair, mode 2 is unstable, but mode 3 and higher are still stable.

\begin{figure}[h]%
\begin{tabular}{ m{.1\textwidth} m{.9\textwidth} }
\begin{minipage}{.1\textwidth} %
$A=3$\\ 
$B=2$
\end{minipage}
& 
\includegraphics[width=0.9\textwidth]{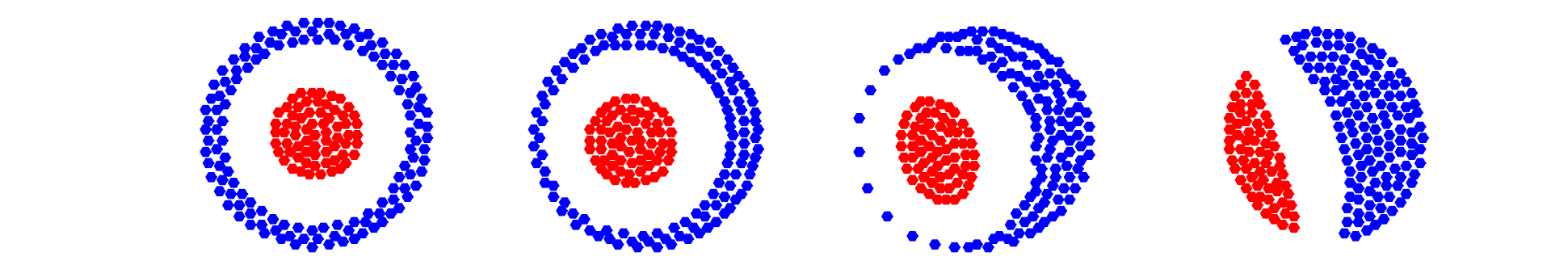}\\
&  \hspace{0.15\textwidth} $t=0$ \hspace{0.1\textwidth} $t=16.5$ \hspace{0.1\textwidth} $t=19$ \hspace{0.11\textwidth} $t=50$
\end{tabular}\\ 
\begin{tabular}{ m{.1\textwidth} m{.9\textwidth} }
\begin{minipage}{.1\textwidth} %
$A=3$\\ 
$B=0.75$
\end{minipage}
&  
\includegraphics[width=0.9\textwidth]{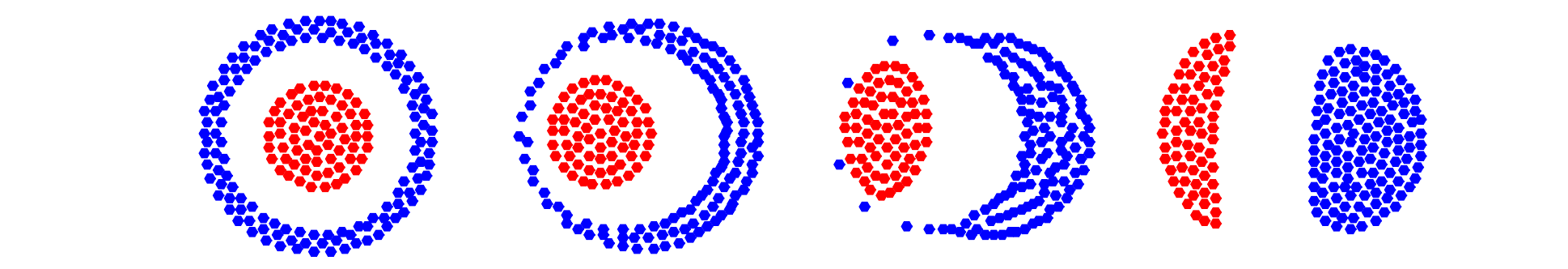}\\
& \hspace{0.15\textwidth} $t=0$ \hspace{0.1\textwidth} $t=6.6$ \hspace{0.11\textwidth} $t=7.4$ \hspace{0.1\textwidth} $t=50$
\end{tabular}
\caption{Steady states that arise, starting from the stable target configuration (lighter species inside). Apparently, for $B<A$, the mode 1 instability is dominant in the transition into a non-radially symmetric steady state. The is even the case if mode 2 is also unstable (see the plots for $B=0.75$).}%
\label{fig:unstab modes targStab}%
\end{figure}

We start from the target particle configuration that follows from the numerics for $(3,3.5)$. To obtain the correct target ansatz for our new choice of parameters $(A,B)$, we subsequently rescale this configuration using \eqref{eqn:radii}. We then perform a numerical run of \eqref{eqn:part syst}, both for $(A,B)=(3,2)$ and for $(A,B)=(3,0.75)$. Some snapshots are shown in Figure \ref{fig:unstab modes targStab}. As expected, for $(A,B)=(3,2)$ a mode 1 instability occurs (this is the only unstable mode), shifting the red core outside. For $(A,B)=(3,0.75)$, when modes 1 and 2 are both unstable, we again observe a mode 1 instability. Apparently, mode 1 dominates mode 2 here (larger eigenvalue in the linearized system of Section \ref{sec:target1 lin}).

\subsubsection{Target: heavy species inside}
In Section \ref{subsect:t-hinside} we focussed on the target state when the heavy species is inside the light species. We showed that is unstable. Hence we do not observe it numerically. In particular, mode 1 is unstable in region $D_2\cup D_3$; that is, for $B<A$. We also showed that mode 2 is unstable if $B>1$. See Section \ref{subsect:lin t-hinside}.

These different modes of instability will be illustrated here. We again design a particle system of 200 particles with $M=2$ and now we pick three pairs of parameter values such that first $B>A$, next $1<B<A$, and finally $B<1$. We construct a target configuration with the heavy species inside and with radii according to \eqref{eqn:radii2} to resemble the target ansatz. Since the target with the heavy species inside is unstable, it does not appear as a steady state in numerical simulations. Therefore, some manipulation is needed to obtain the target that we use as initial configurations. We omit further details.

\begin{figure}[h]%
\begin{tabular}{ m{.1\textwidth} m{.9\textwidth} }
\begin{minipage}{.1\textwidth} %
$A=3$\\ 
$B=3.5$
\end{minipage}
& 
\includegraphics[width=0.9\textwidth]{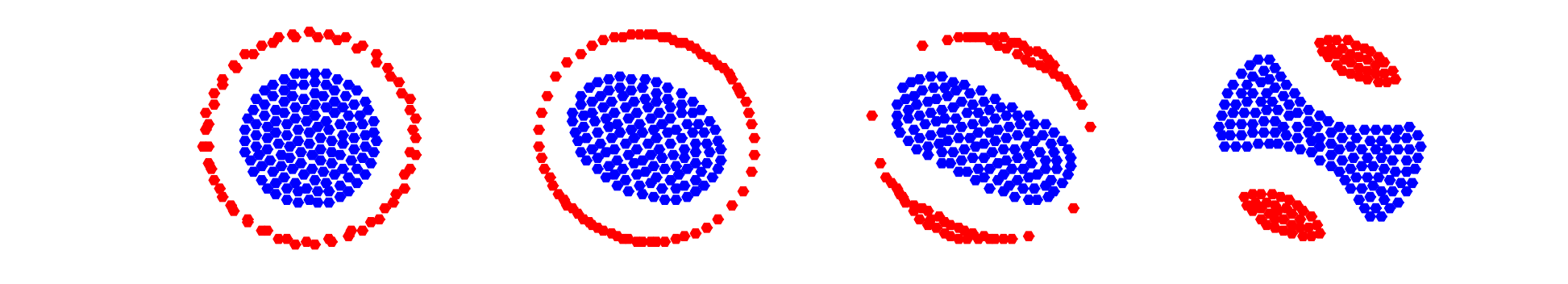}\\
& \hspace{0.15\textwidth} $t=0$ \hspace{0.1\textwidth} $t=6.5$ \hspace{0.11\textwidth} $t=7.5$ \hspace{0.11\textwidth} $t=50$
\end{tabular}\\ 
\begin{tabular}{ m{.1\textwidth} m{.9\textwidth} }
\begin{minipage}{.1\textwidth} %
$A=3$\\ 
$B=2$
\end{minipage}
& 
\includegraphics[width=0.9\textwidth]{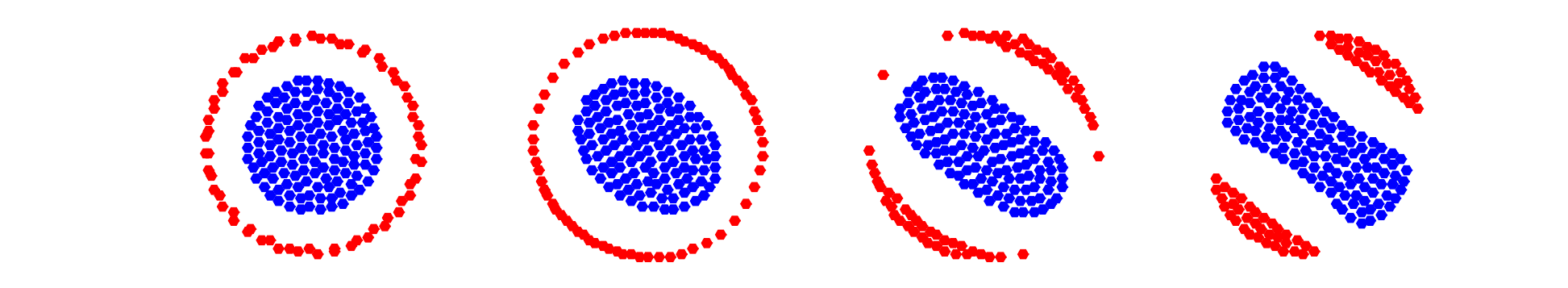}\\
&  \hspace{0.15\textwidth} $t=0$ \hspace{0.1\textwidth} $t=19.5$ \hspace{0.1\textwidth} $t=23$ \hspace{0.11\textwidth} $t=50$
\end{tabular}\\ 
\begin{tabular}{ m{.1\textwidth} m{.9\textwidth} }
\begin{minipage}{.1\textwidth} %
$A=3$\\ 
$B=0.75$
\end{minipage}
&  
\includegraphics[width=0.9\textwidth]{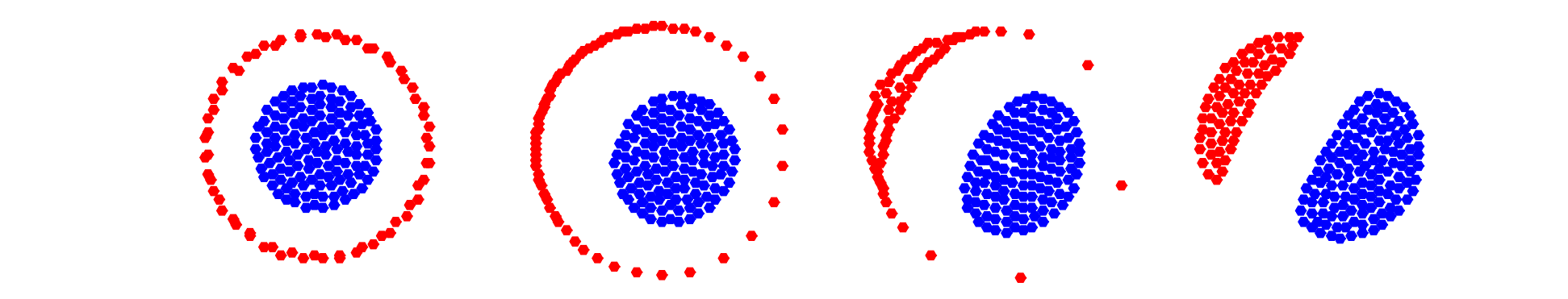}\\
& \hspace{0.15\textwidth} $t=0$ \hspace{0.1\textwidth} $t=16$ \hspace{0.11\textwidth} $t=18$ \hspace{0.11\textwidth} $t=50$
\end{tabular}
\caption{Steady states that arise, starting from the unstable target configuration (heavy inside). For $B>1$, mode 2 is either unstable while mode 1 is not (for $B>A$, see the plots for $B=3.5$), or it apparently dominates mode 1, that is also unstable (for $1<B<A$, see the plots for $B=2$). For $B<1$, mode 1 is stable, while mode 2 is stable (see the plots for $B=0.75$). Apparently the instability is driven by modes 1 or 2, and not by the higher-order modes.}%
\label{fig:unstab modes targUnstab}%
\end{figure}

In Figure \ref{fig:unstab modes targUnstab} we show snapshots of the time evolution for $(A,B)=(3,3.5)$, and $(A,B)=(3,2)$, and $(A,B)=(3,0.75)$. In each case we start from the (properly scaled) target with the heavy species inside. 

For $(A,B)=(3,3.5)$, mode 2 is unstable, while mode 1 is stable. In the top row of Figure \ref{fig:unstab modes targUnstab} we clearly see the mode 2 instability that elongates the blue core and triggers the system to evolve into a non-radially symmetric state. For $(A,B)=(3,2)$ both modes 1 and 2 are unstable. The middle row of Figure \ref{fig:unstab modes targUnstab} shows that apparently the mode 2 instability is dominant. The steady state that follows resembles the one in the top row. Mode 2 is stable for $(A,B)=(3,0.75)$, but mode 1 is not. In the bottom row of Figure \ref{fig:unstab modes targUnstab} the instability of mode 1 is visible as a translation of the blue core.

Note that the steady state on the bottom row of Figure \ref{fig:unstab modes targUnstab} is the same as the steady state at the bottom of Figure \ref{fig:unstab modes targStab}. However, the former arises due to a translation (mode 1 instability) of the blue core (species 1), while the latter arises due to a mode 1 instability of the red core (species 2).

The steady states obtained in Figures \ref{fig:unstab modes targStab} and \ref{fig:unstab modes targUnstab} were previously illustrated in Figure \ref{fig:num st st}.


\section{Overlap equilibrium}
\label{sect:overlap}
In this section we investigate the radially symmetric state where the two species are supported on concentric disks and thus there is a region in which the two species coexist (overlap) -- see the picture for parameter values $(A,B)=(0.5,1)$ in Figure \ref{fig:num st st}. We also consider two versions of this state: one where the coexistence region is surrounded by a ring of the heavier species (i.e. lighter species inside), and one with the heavy and light species interchanged (heavier species inside).

For overlap equilibria we were not able to develop a linear stability analysis as for the target solution (Sections \ref{sec:target1 lin} and \ref{subsect:lin t-hinside}). The difficulties are the ones pointed out in Section \ref{sect:prelim lin}: at the boundary at $R_2$ the velocities of both species 1 and species 2 need to be taken into consideration. That is, \eqref{eqn:dpdt is v pert} should be satisfied for $(j,k)=(2,1)$ and for $(j,k)=(2,2)$. These two equations simultaneously lead to inconsistencies in our approach. For this reason the considerations in this sections are limited to the variational approach.

\subsection{Lighter species inside}
\label{subsect:linside}
\begin{figure}[h]
\centering
\begin{tikzpicture}[>= latex]
\begin{scope}[scale=0.3] 

\draw[fill=blue!20
,line width=1]  (0,0) circle (10);
\draw[fill=purple!50
, line width=1]  (0,0) circle (5.5);            
                
\draw[<->,line width=2] (0.05,0) -- (10,0) node[pos=1,right] {$R_1$};
\draw[<->,line width=2] (0,-0.05) -- (0,-5.5) node[pos=.5,left] {$R_2$};
\draw (0,2) node {\begin{tabular}{c}Species 1 \& 2\end{tabular}}; 
\draw (0,7.7) node {Species 1 only};
\end{scope}
\end{tikzpicture}
\caption{Overlap equilibrium, with species $1$ and $2$ coexisting on a disk of radius $\Rtwo$, and (the heavier) species $1$ also being present in the annular region $\Rtwo<|x|<\Rone$.}
\label{fig:overlap light inside}
\end{figure}
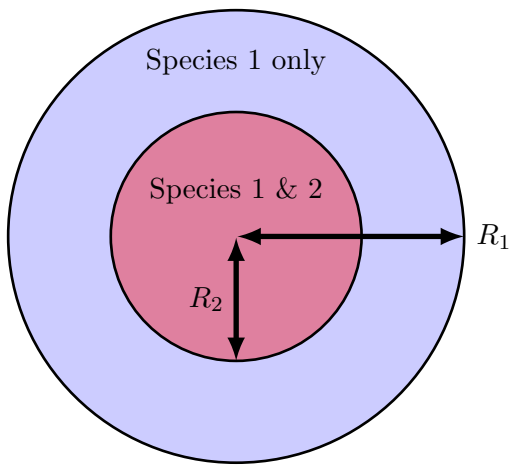

In this equilibrium state species $1$ and $2$ are supported in disks of radii $\Rone$ and $\Rtwo$, respectively, with $\Rtwo<\Rone$ -- see Figure \ref{fig:overlap light inside}. Within $|x|<\Rtwo$, where the two species coexist, the equilibrium densities are (see \eqref{eqn:equilibria}):
\begin{equation}
\label{eqn:equil-overlap}
\barrho_1 = \frac{(\as \bs - \ac \bc)M_1 + (\as \bc - \ac \bs)M_2}{\pi(\as^2 - \ac^2)}, \qquad \barrho_2 = \frac{(\as \bc - \ac \bs)M_1 + (\as \bs - \ac \bc)M_2}{\pi(\as^2 - \ac^2)}.
\end{equation}
In the annular region $\Rtwo<|x|<\Rone$, only species $1$ is present, with equilibrium density (also see \eqref{eqn:equilibria}):
\begin{equation}
\label{eqn:equil-out}
\rhooneout =  \frac{\bs M_1 + \bc M_2}{\pi \as}.
\end{equation}

By immediate calculations, the radii of the two disks are found to be
\begin{equation}
\label{eqn:overlap-radii}
\Rone^2 = \frac{\as M_1 + \ac M_2}{\bs M_1 + \bc M_2}, \qquad \Rtwo^2 = \frac{(\as^2-\ac^2)M_2}{(\as \bc - \ac \bs)M_1 + (\as \bs - \ac \bc) M_2}.
\end{equation}

Together with the consistency condition $\Rtwo<\Rone$, it can be shown that the overlap equilibrium above exists for $(A,B) \in \D_3 \cup \D_6$. Note that $A>B$ in $\D_3$, while in $\D_6$ one has 
$A<B$.

Calculate $\nabla \Lambda_1$ from \eqref{eqn:Lambda1}, with equilibrium densities given by \eqref{eqn:equil-overlap} and \eqref{eqn:equil-out}. Using \eqref{eqn:intderiv} one can check indeed that $\nabla \Lambda_1(x)=0$ in $|x|<\Rtwo$ and $\Rtwo<|x|<\Rone$, hence $\Lambda_1(x)$ is constant in the support of $\barrho_1$ (see \eqref{eqn:equilsup}). Outside the support, in $|x|>\Rone$, one finds

\begin{equation}
\label{eqn:gradL1-overlap}
\Lambda_1'(|x|) = \bigl( (\bs M_1 + \bc M_2)|x|^2 - (\as M_1 + \ac M_2) \bigr)/ |x|.
\end{equation}
From \eqref{eqn:gradL1-overlap} and the expression of $\Rone$ in \eqref{eqn:overlap-radii}, we conclude that $\Lambda_1$ is radially increasing in $|x|>\Rone$, and hence,  for all $(A,B)$ in the relevant region $\D_3 \cup \D_6$, $\Lambda_1$ satisfies \eqref{eqn:equilcomp}.

We now calculate $\nabla \Lambda_2$ from \eqref{eqn:Lambda2}. First, one can check that $\nabla \Lambda_2(x)=0$ in $|x|<\Rtwo$, hence $\Lambda_2(x)$ is constant in the support of $\barrho_2$. Then, the calculation of $\nabla \Lambda_2(x) = \Lambda_2'(|x|)\, x/|x|$ outside the support  yields:
\begin{equation}
\label{eqn:gradL2-overlap}
\frac{\Lambda_2'(|x|)}{|x|}= \begin{cases}
\bc M_1 + \bs M_2 - (\ac M_1 + \as M_2)/|x|^2 - \pi \ac \rhooneout (1-\Rone^2/|x|^2)  &  \text{ if } \Rtwo<|x|<\Rone \\
\bc M_ 1+ \bs M_2 - (\ac M_1 + \as M_2)/|x|^2& \text{ if } |x|>\Rone.
\end{cases}
\end{equation}

Consider case $(A,B) \in \D_6$ first. By \eqref{eqn:overlap-radii}, in $|x|>\Rone$ we have
\begin{align*}
\bc M_ 1+ \bs M_2 - (\ac M_1 + \as M_2)/|x|^2 &> \bc M_ 1+ \bs M_2 - (\ac M_1 + \as M_2)  \frac{\bs M_1 + \bc M_2}{\as M_1 + \ac M_2} \\
&= \bs M_2 \left( BM+1 - (AM+1) \frac{M+B}{M+A} \right) \\
& = \bs M_2 \frac{(B-A)(M^2-1)}{M+A}
\end{align*}
As $B>A$ for $(A,B) \in \D_6$ and $M>1$, the expression above is positive, and consequently, $\Lambda_2(x)$ is radially increasing in $|x|>\Rone$.

In $\Rtwo<|x| < \Rone$,  by \eqref{eqn:equil-out} and \eqref{eqn:ABM}, the expression on the right-hand-side of \eqref{eqn:gradL2-overlap} can be rewritten as:
\begin{equation}
\label{eqn:gradL2-mod}
\bc M_1 + \bs M_2 - (\bs M_1 + \bc M_2) A + \underbrace{\left( (\as M_1 + \ac M_2) A - \ac M_1 - \as M_2 \right)}_{=\as M_2(A^2-1) <0 \text{ in } \D_6}/|x|^2,
\end{equation}
and hence, also using \eqref{eqn:overlap-radii},
\begin{align*}
\frac{\Lambda_2'(|x|)}{|x|} &> \bc M_1 + \bs M_2 - (\bs M_1 + \bc M_2) A + \as M_2(A^2-1)/\Rtwo^2  \,=\, 0.
\end{align*}

We conclude that $\Lambda_2$ is radially increasing in $|x|>\Rtwo$. Since $\Lambda_1$ and $\Lambda_2$ satisfy \eqref{eqn:equilcomp}, the overlap solution is a local minimizer (with respect to perturbations of class $\CalB$) when $(A,B) \in \D_6$.

Next, consider the case $(A,B) \in \D_3$, and take the expression for the right-hand-side of \eqref{eqn:gradL2-overlap} that was derived in \eqref{eqn:gradL2-mod}:
\[
\bc M_1 + \bs M_2 - (\bs M_1 + \bc M_2) A + \underbrace{\left( (\as M_1 + \ac M_2) A - \ac M_1 - \as M_2 \right)}_{=\as M_2(A^2-1) >0 \text{ in } \D_3}/|x|^2.
\]
Then, in $\Rtwo<|x| < \Rone$,
\begin{align*}
\frac{\Lambda_2'(|x|)}{|x|} &< \bc M_1 + \bs M_2 - (\bs M_1 + \bc M_2) A + \as M_2(A^2-1)/\Rtwo^2 \,=\, 0.
\end{align*}

Also, $\Lambda_2'$ changes sign at 
\[
|x|^2 = \frac{\ac M_1 + \as M_2}{\bc M_1 + \bs M_2} > \Rone, 
\]
and stays positive beyond that. We conclude that the overlap solution is not a minimizer for $(A,B) \in \D_3$.

In summary, the overlap solution with the lighter species inside (Figure \ref{fig:overlap light inside}) is a local minimizer with respect to class $\CalB$ perturbations for $(A,B) \in D_6$, but not for $(A,B) \in D_3$.

\begin{remark}[Global minimum]
\label{rmk:overlap-gmin}
As for the target equilibrium (see Remark \ref{rmk:target-gmin}), we resort again to the calculations for the second variation of the energy in Section \ref{subsect:prelim-var} and conclude that the overlap equilibrium with the lighter species inside is a {\em global} minimizer for all $(A,B)$ in $D_6$ with $B>1$ ($A<1$ holds in all $D_6$); i.e.~in the intersection of $D_6$ with the shaded region in Figure \ref{fig:stability}.  Note that this restriction excludes only the bounded triangular region $0<A<B<1$, while $D_6$ is unbounded! Hence, by the variational approach we identified two global minimizers (overlap and target equilibria with lighter species inside) which exist in unbounded (and disjoint) subsets of the parameter space $(A,B)$.
\end{remark}


\subsection{Heavier species inside}
\label{subsect:hinside}
\begin{figure}[h]
\centering
\begin{tikzpicture}[>= latex]
\begin{scope}[scale=0.3] 

\draw[fill=red!20
,line width=1]  (0,0) circle (10);
\draw[fill=purple!50
, line width=1]  (0,0) circle (7);            

\draw[<->,line width=2] (0.05,0) -- (10,0) node[pos=1,right] {$R_2$};
\draw[<->,line width=2] (0,-0.05) -- (0,-7) node[pos=.5,left] {$R_1$};
\draw (0,3) node {\begin{tabular}{c}Species 1 \& 2\end{tabular}}; 
\draw (0,8) node {Species 2 only};
\end{scope}
\end{tikzpicture}
\caption{Overlap equilibrium, with species $1$ and $2$ coexisting on a disk of radius $\Rone$, and (the lighter) species $2$ also being present in the annular region $\Rone<|x|<\Rtwo$.}
\label{fig:overlap heavy inside}
\end{figure}
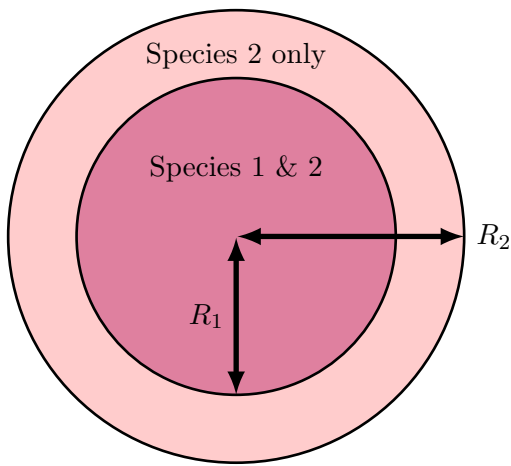

For this equilibrium species $1$ and $2$ are supported in disks of radii $\Rone$ and $\Rtwo$, respectively, with $\Rone<\Rtwo$ -- see Figure \ref{fig:overlap heavy inside}. The heavier species 1 is now {\em inside}. In $|x|<\Rone$, where the two species coexist, the equilibrium densities are also given by \eqref{eqn:equil-overlap} (cf. \eqref{eqn:equilibria}). In the annular region $\Rone<|x|<\Rtwo$, only species $2$ is present, with equilibrium density (also see \eqref{eqn:equilibria}):
\begin{equation}
\label{eqn:equil-out2}
\rhotwoout =  \frac{\bc M_1 + \bs M_2}{\pi \as}.
\end{equation}

The radii of the two disks are given by
\begin{equation}
\label{eqn:overlap-radii2}
\Rone^2 =\frac{(\as^2-\ac^2)M_1}{(\as \bs - \ac \bc)M_1 + (\as \bc - \ac \bs) M_2}, \qquad \Rtwo^2 = \frac{\ac M_1 + \as M_2}{\bc M_1 + \bs M_2}.
\end{equation}

Together with the consistency condition $\Rone<\Rtwo$, it can be shown that the overlap equilibrium above exists for $(A,B) \in \D_1 \cup \D_4$. Note that $A>B$ in $\D_1$, while in $\D_4$ one has 
$A<B$. 

Calculate $\nabla \Lambda_1$. We find $\nabla \Lambda_1(x)=0$ in $|x|<\Rone$, as required for equilibrium (see \eqref{eqn:equilsup}). Outside the support,
\begin{equation}
\label{eqn:gradL1-overlap2}
\frac{\Lambda_1'(|x|)}{|x|}= \begin{cases}
\bs M_1 + \bc M_2 - (\as M_1 + \ac M_2)/|x|^2 - \pi \ac \rhotwoout (1-\Rtwo^2/|x|^2)  &  \text{ if } \Rone<|x|<\Rtwo \\
\bs M_ 1+ \bc M_2 - (\as M_1 + \ac M_2)/|x|^2& \text{ if } |x|>\Rtwo.
\end{cases}
\end{equation}

Consider case $(A,B) \in \D_1$ first. By \eqref{eqn:overlap-radii2}, in $|x|>\Rtwo$ we have
\begin{align*}
\bs M_ 1+ \bc M_2 - (\as M_1 + \ac M_2)/|x|^2 &> \bs M_ 1+ \bc M_2 - (\as M_1 + \ac M_2)  \frac{\bc M_1 + \bs M_2}{\ac M_1 + \as M_2} \\
&= \bs M_2 \left( M+B - (M+A) \frac{BM+1}{AM+1} \right) \\
& = \bs M_2 \frac{(A-B)(M^2-1)}{AM+1}.
\end{align*}
As $A>B$ for $(A,B) \in \D_1$ and $M>1$, the expression above is positive, and consequently, $\Lambda_1(x)$ is radially increasing in $|x|>\Rtwo$.

In $\Rone<|x| < \Rtwo$,  by \eqref{eqn:equil-out2} and \eqref{eqn:ABM}, the expression on the right-hand-side of \eqref{eqn:gradL1-overlap2} can be rewritten as:
\begin{equation}
\label{eqn:gradL2-mod2}
\bs M_1 + \bc M_2 - (\bc M_1 + \bs M_2) A + (\underbrace{(\ac M_1 + \as M_2)A - \as M_1 - \ac M_2}_{\as M_2 (A^2-1)M < 0 \text{ in } \D_1 }) /|x|^2
\end{equation}
and hence, by \eqref{eqn:overlap-radii2},
\begin{align*}
\frac{\Lambda_1'(|x|)}{|x|} &> \bs M_1 + \bc M_2 - (\bc M_1 + \bs M_2) A  + \as M_2(A^2-1)M/\Rone^2 \, =\, 0.
\end{align*}
We conclude that $\Lambda_1$ is radially increasing in $|x|>\Rone$ (and satisfies \eqref{eqn:equilcomp}) for all $(A,B) \in \D_1$.

On the other hand, in $\D_4$ (where $A>1$), using \eqref{eqn:gradL2-mod2} and an argument analogous to the above, we find that $\Lambda_1'(|x|)<0$ and hence \eqref{eqn:equilcomp} is violated.

Finally, for $\Lambda_2$ we find $\nabla \Lambda_2(x)=0$ in $|x|<\Rtwo$ (as for equilibrium), and outside the support, in $|x| >\Rtwo$:
\begin{equation*}
\Lambda_2'(|x|) = \bigl( (\bc M_1 + \bs M_2)|x|^2 - (\ac M_1 + \as M_2) \bigr)/ |x|.
\end{equation*}
Consequently, by the expression of $\Rtwo$ in \eqref{eqn:overlap-radii2}, $\Lambda_2$ is radially increasing in $|x|>\Rtwo$.

In conclusion, the overlap solution with the heavier species inside is a minimizer with respect to perturbations of class $\CalB$ when $(A,B) \in \D_1$, but not for $(A,B) \in \D_4$. Based on particle simulations we conjecture however that this overlap solution is not a minimizer with respect to class $\CalA$ perturbations for any $(A,B) \in \D_1$; such an equilibrium has never been captured in simulations in fact. Unfortunately, unlike for the target solution, we do not have a linear stability analysis to support such a conjecture. 


\subsection{Numerical illustration: passing from $D_6$ to $D_1$}\label{sec:instab overlap D6 D1}
The overlap state with the lighter species inside (Section \ref{subsect:linside}) is observed numerically as a steady state in region $D_6$ and we indicated this accordingly in Figure \ref{fig:num st st}. Below we illustrate numerically how equilibria change when parameters cross the line $A=B$.

\begin{figure}[h]%
\begin{tabular}{ m{.1\textwidth} m{.9\textwidth} }
\begin{minipage}{.1\textwidth} %
$A=0.5$\\ 
$B=0.4$
\end{minipage}
& 
\includegraphics[width=0.9\textwidth]{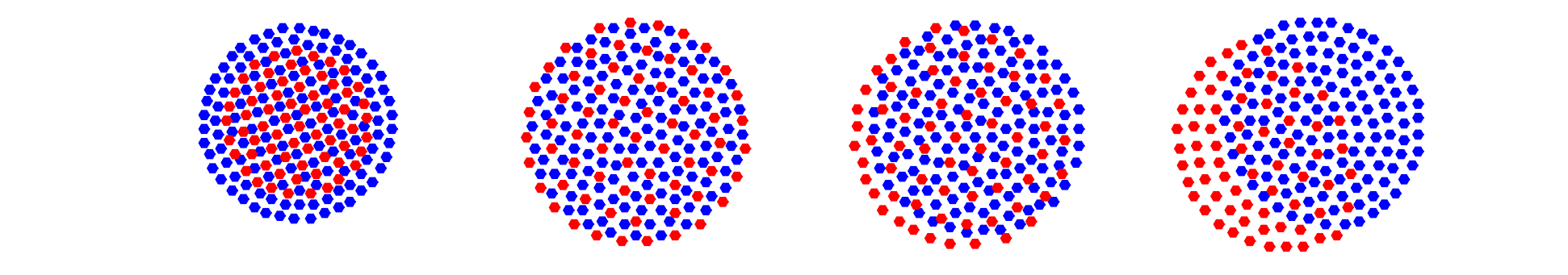}\\
& \hspace{0.15\textwidth} $t=0$ \hspace{0.1\textwidth} $t=50$ \hspace{0.1\textwidth} $t=220$ \hspace{0.1\textwidth} $t=500$
\end{tabular}
\caption{Steady state arising in parameter region $D_1$, starting from the overlap state (lighter species inside) that is stable in $D_6$.}%
\label{fig:unstab modes overlap}%
\end{figure}

In Figure \ref{fig:unstab modes overlap} we show a series of snapshots. The initial configuration is the overlap state that we find as the long-time steady state for $(A,B)=(0.5,1)$. Next we change the parameters to $(A,B)=(0.5,0.4)$, that is, we cross the boundary between $D_6$ and $D_1$. We observe that the radial symmetry is broken, and the system attains a state in which the supports of the two species partially overlap. It turns out that, at least asymptotically close to $(A,B)=(0,0)$, we can quantify this effect, and in particular we can find the distance between the species' centres of mass. This is done in Section \ref{sec:weak cross}.

%
%
%
%
%
\section{Weak cross-interactions}\label{sec:weak cross}
Consider the case in which the cross-interactions are much weaker than the self-interactions. We consider a small parameter $0<\eta\ll 1$ and substitute $K_c$ by $\eta\,K_c$ in \eqref{eqn:model}.
This system exhibits a `regular' timescale and a slow timescale (we also observe this in numerics; see Figure \ref{fig:timescales}). We will now examine this separation of timescales and its implications for the steady state.

Introduce a two-scale expansion in \eqref{eqn:model} given by the variables $t$ and $s:=\eta\,t$. Taking the transformation $\frac{\partial}{\partial t}\mapsto \frac{\partial}{\partial t} + \eta\,\frac{\partial}{\partial s}$ into account, the two-scale model equations are
\begin{subequations}
\label{eqn:model 2scale}
\begin{gather}
\frac{\partial \rho_{1}}{\partial t}+\eta\frac{\partial \rho_{1}}{\partial s}+\nabla\cdot(\rho_1 v_1)=0, \quad v_1=-\nabla K_s\ast\rho_1 - \eta\nabla K_c \ast \rho_2,\\
\frac{\partial \rho_{2}}{\partial t}+\eta\frac{\partial \rho_{2}}{\partial s}+\nabla\cdot(\rho_2 v_2)=0, \quad v_2=-\eta\nabla K_c\ast\rho_1 - \nabla K_s \ast \rho_2,
\end{gather}
\end{subequations}
and they can be separated into
\begin{align}
\Ord(\eta^0):&\qquad \frac{\partial \rho_{i}}{\partial t}+\nabla\cdot(\rho_i (-\nabla K_s\ast\rho_i))=0, &i=1,2,\label{eqn:2scale order 0}\\
\Ord(\eta^1):& \qquad \frac{\partial \rho_{i}}{\partial s}+\nabla\cdot(\rho_i (-\nabla K_c\ast\rho_j))=0, &i=1,2,\quad j\neq i.\label{eqn:2scale order 1}
\end{align}
Setting $\frac{\partial}{\partial t}\rho_i=0$ and $\frac{\partial}{\partial s}\rho_i=0$ in \eqref{eqn:2scale order 0}--\eqref{eqn:2scale order 1}, we obtain the following conditions for a steady state:
\begin{align}
-\nabla K_s\ast\bar\rho_i=0,& \quad \text{ on } \supp \bar\rho_i, \text{ for each }i=1,2,\label{eqn:steady single}\\
-\nabla  K_c\ast\bar\rho_j=0,& \quad \text{ on } \supp \bar\rho_i, \text{ for each }i=1,2,\quad j\neq i,\label{eqn:steady cross}
\end{align}
where $\bar\rho_1$ and $\bar\rho_2$ denote the steady state densities. By the first equation \eqref{eqn:steady single}, we know that each species independently attains the steady state of a single, isolated species corresponding to the potential $K_s$. In fact, the zeroth-order equation \eqref{eqn:2scale order 0} suggests that each species approaches this steady state at the `regular' timescale $t$. Note that these steady states are determined per species up to translation of the centre of mass.

The equilibrium distance between the centres of mass can be derived from \eqref{eqn:steady cross}. Integration of $-\nabla K_c\ast\bar\rho_2=0$ over $\supp \bar\rho_1$ yields
\begin{equation}\label{eqn:cond cross weak}
\int_{\supp \bar\rho_1}\int_{\supp \bar\rho_2} \nabla  K_c(x-y)  \bar\rho_2(y) \bar\rho_1(x)\,dy dx = 0.
\end{equation}
We note that due to the assumed antisymmetry of $\nabla K_c$, the same condition is obtained if we integrate $-\nabla K_c\ast\bar\rho_1=0$ over $\supp \bar\rho_2$. This condition \eqref{eqn:cond cross weak} holds for general cross-interaction kernel $ K_c$.

\subsection{Newtonian-quadratic interactions}
We now take the same interaction potentials $K_s$ and $K_c$ as in \eqref{eqn:intpot}. 
By taking $\eta$ small, we are zooming in at the origin in Figure \ref{fig:stability}. Analogous to \eqref{eqn:ABM}, the dimensionless numbers $A$ and $B$ are defined as the ratios of the cross- and self-interaction parameters. Here we take into account that the cross-interactions are pre-multiplied by $\eta$, and we have $A:=\eta\, a_c/a_s$ and $B:=\eta\, b_c/b_s$.

From \eqref{eqn:steady single} and \eqref{eqn:selfpot} it follows due to \cite{BertozziLaurentLeger, FeHuKo11} that the steady state densities are (to leading order) of the form:
\begin{equation}\label{eq:st st dens weak}
\bar\rho_1(x)= \dfrac{b_s\,M_1}{\pi\,a_s}\chi_{B(x_0,R)}(x), \qquad \bar\rho_2(x)= \dfrac{b_s\,M_2}{\pi\,a_s}\chi_{B(\bar x_0,R)}(x),
\end{equation}
for some $x_0, \bar x_0 \in \R^2$ and with $R^2:=a_s/b_s$. Here, $\chi$ is the characteristic function. Note that both supports have the same radius, even though the masses $M_1$ and $M_2$ are in general unequal. Without loss of generality, take $\bar x_0 = x_0 + (d,0)^T$ for some constant $d>0$. The condition \eqref{eqn:cond cross weak} can now be written as
\begin{equation}\label{eqn:cond cross weak NewtLin}
\int_{B(x_0,R)}\int_{B(\bar x_0,R)} \left[ a_c\, \dfrac{x-y}{|x-y|^2} - b_c\,(x-y) \right]\,dy dx = 0.
\end{equation}
and we now show that this can be reduced to a relation between $d$ and the model parameters. The attraction part is evaluated explicitly as
\begin{equation}\label{eqn:double int attr}
-b_c\,\int_{B(x_0,R)}\int_{B(\bar x_0,R)} (x-y) \,dydx= b_c\,\pi^2R^4(\bar x_0 - x_0)=\dfrac{a_s^2\,b_c\,\pi^2}{b_s^2}\left(\!\!\begin{array}{c} d\\0 \end{array}\!\! \right).
\end{equation}
For the repulsion part, we distinguish between the following cases:
\begin{description}
\item[Case 1: $d>2R$.] In this case $B(x_0,R)\cap B(\bar{x}_0,R)=\emptyset$, hence \eqref{eqn:intderiv} implies that for all $x\in B(x_0,R)$ 
\begin{equation}
\int_{B(\bar{x}_0,R)}\dfrac{x-y}{|x-y|^2}\,dy = \pi\,R^2\,\dfrac{x-\bar{x}_0}{|x-\bar{x}_0|^2}.
\end{equation}
Subsequently, \eqref{eqn:intderiv} yields that
\begin{equation}
\int_{B(x_0,R)}\pi\,R^2\,\dfrac{x-\bar{x}_0}{|x-\bar{x}_0|^2}\,dx = \pi^2\,R^4\,\dfrac{x_0-\bar{x}_0}{|x_0-\bar{x}_0|^2},\label{eqn:double int rep d>2R}
\end{equation}
because $\bar{x}_0\notin B(x_0,R)$. Noting that $x_0-\bar{x}_0=-d\,\mathbf{e}_1$, we conclude for $d>2R$ that \eqref{eqn:cond cross weak NewtLin} is equivalent to
\begin{equation}
-a_c\,\pi^2R^4\dfrac{1}{d}\,\mathbf{e}_1 + b_c\,\pi^2\,R^4\,d\,\mathbf{e}_1  =0,
\end{equation}
where we substituted \eqref{eqn:double int attr} and \eqref{eqn:double int rep d>2R}. We recall that $R^2=a_s/b_s$, $A:=\eta\, a_c/a_s$ and $B:=\eta\, b_c/b_s$. Hence,
\begin{equation}
\dfrac{d}{R} = \sqrt {\dfrac{A}{B}}. \label{eqn:rel d ab for d large}
\end{equation}
Therefore, \eqref{eqn:rel d ab for d large} implies that \emph{complete separation} of the two species (that is, $d>2R$), takes place for
\[A/B>4.\]

\item[Case 2: $R<d\leqs 2R$.] Define $B_1:= B(x_0,R)\cap B(\bar{x}_0,R)$ and $B_2:=B(x_0,R)\setminus B(\bar{x}_0,R)$. For the repulsion part of \eqref{eqn:cond cross weak NewtLin}, it holds --cf.~\eqref{eqn:intderiv}-- that
\begin{align}
\nonumber \int_{B(x_0,R)}\int_{B(\bar{x}_0,R)} \dfrac{x-y}{|x-y|^2}\,dy\,dx =&\, \int_{B_1} \pi\,(x-\bar{x}_0)\,dx + \int_{B_2} \pi\,R^2\,\dfrac{x-\bar{x}_0}{|x-\bar{x}_0|^2}\,dx\\
\nonumber =&\, \int_{B_1} \pi\,(x-\bar{x}_0)\,dx\\
&+ \int_{B(x_0,R)} \pi\,R^2\,\dfrac{x-\bar{x}_0}{|x-\bar{x}_0|^2}\,dx - \int_{B_1}\pi\,R^2\,\dfrac{x-\bar{x}_0}{|x-\bar{x}_0|^2}\,dx \label{eqn:decomp int rep R<d<2R}
\end{align}
The area of $B_1$ is $2R^2\arccos(d/(2R))-d/2\,\sqrt{4R^2-d^2}$ while, by construction, its centre of mass is $x_0+(d/2,0)^T$. Therefore
\begin{equation}
\int_{B_1} \pi\,(x-\bar{x}_0)\,dx = -\pi d \left[R^2\arccos\!\left(\dfrac{d}{2R}\right)-\dfrac d4\,\sqrt{4R^2-d^2}\right]\,\mathbf{e}_1.\label{eqn:int x-barx0 on B1}
\end{equation}
For $R<d\leqs 2R$, we have that $\bar{x}_0\notin B(x_0,R)$. Thus, it follows from \eqref{eqn:intderiv} that
\begin{equation}
\int_{B(x_0,R)} \pi\,R^2\,\dfrac{x-\bar{x}_0}{|x-\bar{x}_0|^2}\,dx = \pi^2\,R^4\,\dfrac{x_0-\bar{x}_0}{|x_0-\bar{x}_0|^2}.\label{eqn:int whole ball x0}
\end{equation}
It remains to evaluate the integral over $B_1$ in the last line of \eqref{eqn:decomp int rep R<d<2R}. For symmetry reasons, this integral is a vector in the direction of $x_0-\bar{x}_0$, that is, in the direction of $\mathbf{e}_1$. Consider therefore
\begin{align}
\nonumber \int_{B_1}\dfrac{x-\bar{x}_0}{|x-\bar{x}_0|^2}\,dx \cdot \mathbf{e}_1 =&\,  \int_{B_1}\nabla\ln\left(\dfrac{|x-\bar{x}_0|}{R}\right)\cdot \mathbf{e}_1\,dx \\
=&\, \int_{\partial B_1} \ln\left(\dfrac{|x-\bar{x}_0|}{R}\right)\,\,\mathbf{e}_1\cdot \hat{n}(x)\,dS(x),\label{eqn:int B1 rep is bdry int}
\end{align} 
where the last step follows from Gauss' theorem. The boundary of $B_1$ consists of two circular segments, that are subsets of $\partial B(x_0,R)$ and $\partial B(\bar{x}_0,R)$, respectively. Call these segments $\partial B_\alpha\subset\partial B(x_0,R)$ and $\partial B_\beta\subset\partial B(\bar{x}_0,R)$, such that $\partial B_\alpha \cup \partial B_\beta = \partial B_1$. Note that for $x\in\partial B_\beta$ it holds that $|x-\bar{x}_0|=R$, hence $\ln(|x-\bar{x}_0|/R)=0$, and therefore the contribution of the integration over $\partial B_\beta$ in \eqref{eqn:int B1 rep is bdry int} is zero. Thus
\begin{align}
\nonumber \int_{B_1}\dfrac{x-\bar{x}_0}{|x-\bar{x}_0|^2}\,dx\cdot \mathbf{e}_1 =&\, \int_{\partial B_\alpha} \ln\left(\dfrac{|x-\bar{x}_0|}{R}\right)\,\,\mathbf{e}_1\cdot \hat{n}(x)\,dS(x),
\end{align}
while for $x\in\partial B_\alpha$, we have $x=x_0+R(\cos\theta,\sin\theta)^T$, $\hat{n}(x)=(\cos\theta,\sin\theta)^T$ and $dS(x)=R\,d\theta$ with $-\gamma\leqs \theta\leqs \gamma$ and $\gamma:=\arccos(d/(2R))$. Consequently,
\begin{align}
\nonumber \int_{\partial B_\alpha} \ln\left(\dfrac{|x-\bar{x}_0|}{R}\right)\,\,\mathbf{e}_1\cdot \hat{n}(x)\,dS(x) =&\, R\int_{-\gamma}^{\gamma} \ln \bigg(\bigg|  \left(\!\begin{array}{c}\cos\theta\\\sin\theta\end{array}\!\right)+ \dfrac1R\underbrace{(x_0-\bar{x}_0)}_{=-d\,\mathbf{e}_1} \bigg|\bigg)\,\cos\theta\,d\theta\\
=&\, R\int_{-\gamma}^{\gamma} \ln \left(\sqrt{1+\frac{d^2}{R^2}-\frac{2d}{R}\,\cos\theta}\right)\,\cos\theta\,d\theta. \label{eqn:int partial B alpha}
\end{align}
A combination of \eqref{eqn:cond cross weak NewtLin}, and \eqref{eqn:double int attr}, \eqref{eqn:decomp int rep R<d<2R}, \eqref{eqn:int x-barx0 on B1}, \eqref{eqn:int whole ball x0} and \eqref{eqn:int partial B alpha} yields that $d/R$ is implicitly defined as a function of $A/B$ by
\begin{align}
\nonumber 0 =&\, -a_c\pi\, d \left[R^2\arccos\!\left(\dfrac{d}{2R}\right)-\dfrac d4\,\sqrt{4R^2-d^2}\right] - \dfrac{a_c\,\pi^2\,R^4}{d} \\
&\, -a_c\pi \,R^3\int_{-\gamma}^{\gamma} \ln\! \left(\sqrt{1+\frac{d^2}{R^2}-\frac{2d}{R}\,\cos\theta}\right)\,\cos\theta\,d\theta +b_c\pi^2\,R^4\,d, \label{eqn:impl eqn d for d intermediate} 
\end{align}
or
\begin{equation}
\dfrac{A}{B} = \dfrac{ \dfrac{d^2}{R^2}}{1+\dfrac{ d^2}{\pi R^2}\left[\gamma-\dfrac{d}{4R}\,\sqrt{4-\dfrac{d^2}{R^2}}\right] +  \dfrac{d}{2\pi\,R}\displaystyle\int_{-\gamma}^{\gamma} \ln\! \left(1+\frac{d^2}{R^2}-\frac{2d}{R}\cos\theta\right)\,\cos\theta\,d\theta}. \label{eqn:impl eqn d A/B for d intermediate} 
\end{equation}
Here we used that $R^2=a_s/b_s$ and we recall that $\gamma=\gamma(d/R)=\arccos(d/(2R))$.

\item[Case 3: $d\leqs R$.] In this case $\bar{x}_0\in B(x_0,R)$. The only difference with the case $R<d\leqs2R$ is therefore the evaluation of the integral over $B(x_0,R)$ in \eqref{eqn:decomp int rep R<d<2R}. Thus, we replace \eqref{eqn:int whole ball x0} by
\begin{equation}
\int_{B(x_0,R)} \pi\,R^2\,\dfrac{x-\bar{x}_0}{|x-\bar{x}_0|^2}\,dx = \pi^2\,R^2\,(x_0-\bar{x}_0)=-\pi^2\,d\,R^2\,\mathbf{e}_1.\label{eqn:int whole ball x0 d small}
\end{equation}
For $d\leqs R$, the analogue of \eqref{eqn:impl eqn d for d intermediate}--\eqref{eqn:impl eqn d A/B for d intermediate} is
\begin{align}
\nonumber 0 =&\, -a_c\pi\, d \left[R^2\arccos\!\left(\dfrac{d}{2R}\right)-\dfrac d4\,\sqrt{4R^2-d^2}\right] - a_c\,\pi^2\,d\,R^2 \\
&\, -a_c\pi \,R^3\int_{-\arccos\left(\frac{d}{2R}\right)}^{\arccos\left(\frac{d}{2R}\right)} \ln\! \left(\sqrt{1+\frac{d^2}{R^2}-\frac{2d}{R}\,\cos\theta}\right)\,\cos\theta\,d\theta +b_c\pi^2\,R^4\,d, \label{eqn:impl eqn d for d small} 
\end{align}
or
\begin{equation}
\dfrac{A}{B} = \dfrac{ \dfrac{\pi\,d}{R}}{\dfrac{\pi\,d}{R}+\dfrac{ d}{R}\left[\gamma-\dfrac{d}{4R}\,\sqrt{4-\dfrac{d^2}{R^2}}\right] + \dfrac12 \displaystyle\int_{-\gamma}^{\gamma} \ln\! \left(1+\frac{d^2}{R^2}-\frac{2d}{R}\,\cos\theta\right)\,\cos\theta\,d\theta}. \label{eqn:impl eqn d A/B for d small} 
\end{equation}
We recall that $A:=\eta \, a_c/a_s$ and $B:=\eta \, b_c/b_s$, while $\gamma=\arccos(d/(2R))$.\\
\\
Taking the limit $d\downarrow 0$ in the right-hand side of \eqref{eqn:impl eqn d A/B for d small}, we obtain $A/B=1$. This is the threshold value for full mixing.

\end{description}

Note that in \eqref{eq:st st dens weak} we concluded that (to leading order) the values of the steady state densities are the same as for the single species model. That is: $\bar\rho_i=b_s\,M_i/(\pi\,a_s)$. In \eqref{eqn:equilibria} we identified the steady state densities in a two-species model, and we can use these expressions to verify \eqref{eq:st st dens weak}. In the setting with weak cross-interaction parameters $\eta\,a_c$ and $\eta\,b_c$, we find analogously to \eqref{eqn:equilibria} that the equilibrium values for the density are
\begin{align}
\nonumber \bar\rho_1 =&\,  \frac{\bs M_1 + \eta\,\bc M_2}{\pi \as} \qquad \text{in regions where only one species exists,}\\
\nonumber \bar\rho_1 =&\, \frac{(\as \bs - \eta^2\,\ac \bc)M_1 + \eta\,(\as \bc - \ac \bs)M_2}{\pi(\as^2 - \eta^2\ac^2)}  \qquad \text{in overlap regions;}   
\end{align}
and similar expressions for species 2. \emph{To leading order}, these two expresions are however the same, and both are equal to the single species value $\bar\rho_1=b_s\,M_1/(\pi\,a_s)$. The deviations are only $\Ord(\eta)$.

\subsection{Numerical illustration}
\begin{figure}
\centering
\includegraphics[width=0.9\textwidth]{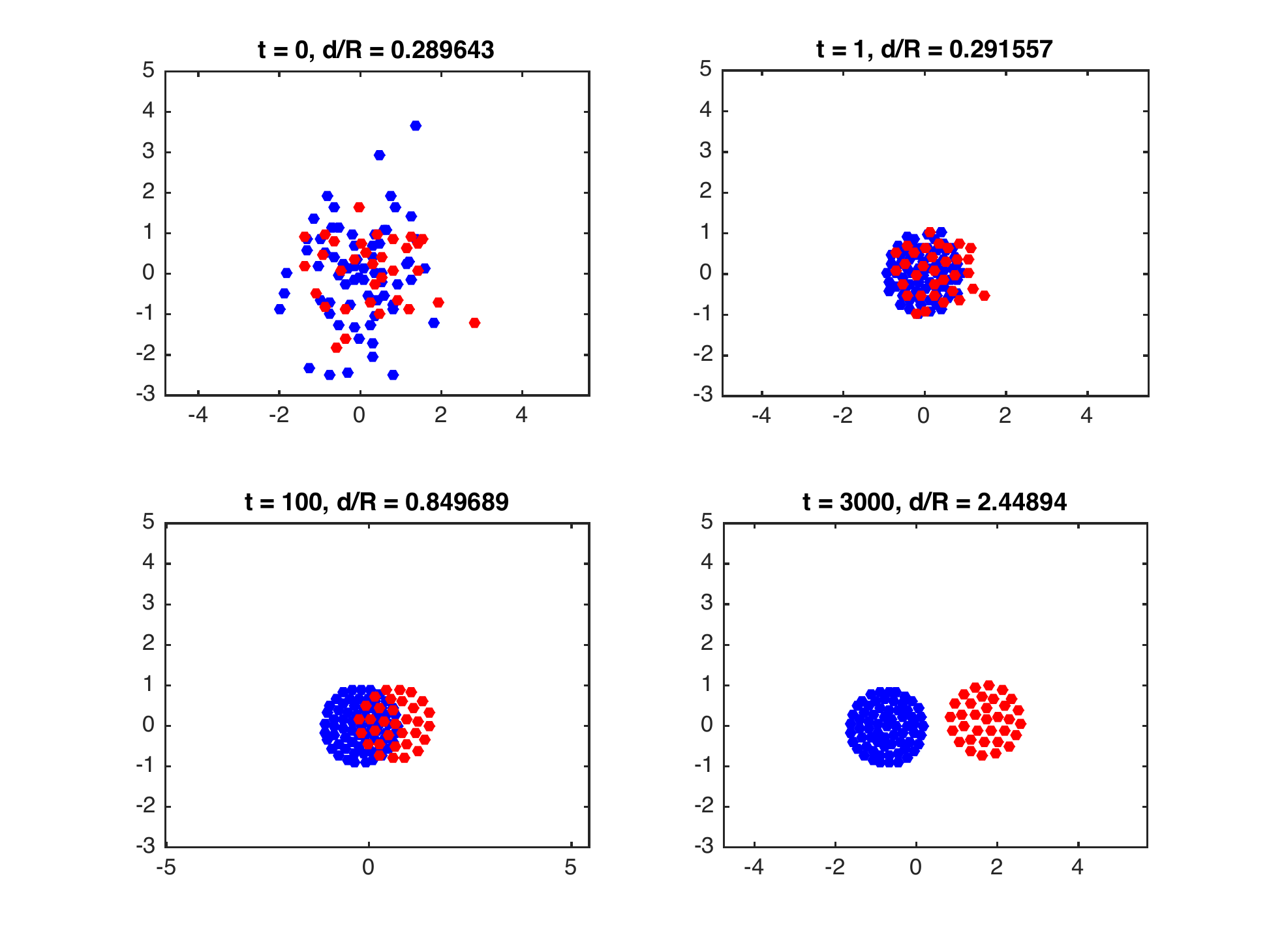}
\vspace{-1cm}\caption{Time evolution of the particle system. There are 100 particles, while $M=2$ (particles are distributed 67:33). We took $a_s=b_s=b_c=1$ and $a_c=6$. For each plot we indicate $d/R$, which is the distance between the two centres of mass divided by $R=\sqrt{a_s/b_s}$. This value approaches the theoretically derived prediction $\sqrt{A/B}=\sqrt{6}\approx 2.45$.}
\label{fig:timescales}%
\end{figure}

We start with an illustration of the two timescales present in the model with weak cross-interactions. We use a particle system of 100 particles with $M=2$. That is, we have 67 particles of species 1 and 33 particles of species 2. We take $\eta=0.001$, and furthermore $a_s=b_s=b_c=1$ and $a_c=6$. The particles are initialized randomly. Their time evolution is shown in Figure \ref{fig:timescales}. We clearly see that each species self-organizes fast into a circular shape. The distance between the two circles equilibrates at a much larger timescale. The distance between the centres of mass of the two species is computed and divided by $R=\sqrt{a_s/b_s}$ to get an estimate for $d/R$ (indicated above each plot in Figure \ref{fig:timescales}). The values obtained slowly approach the limit value $\sqrt{A/B}=\sqrt{6}\approx 2.45$, as predicted by our asymptotic analysis; see \eqref{eqn:rel d ab for d large}.

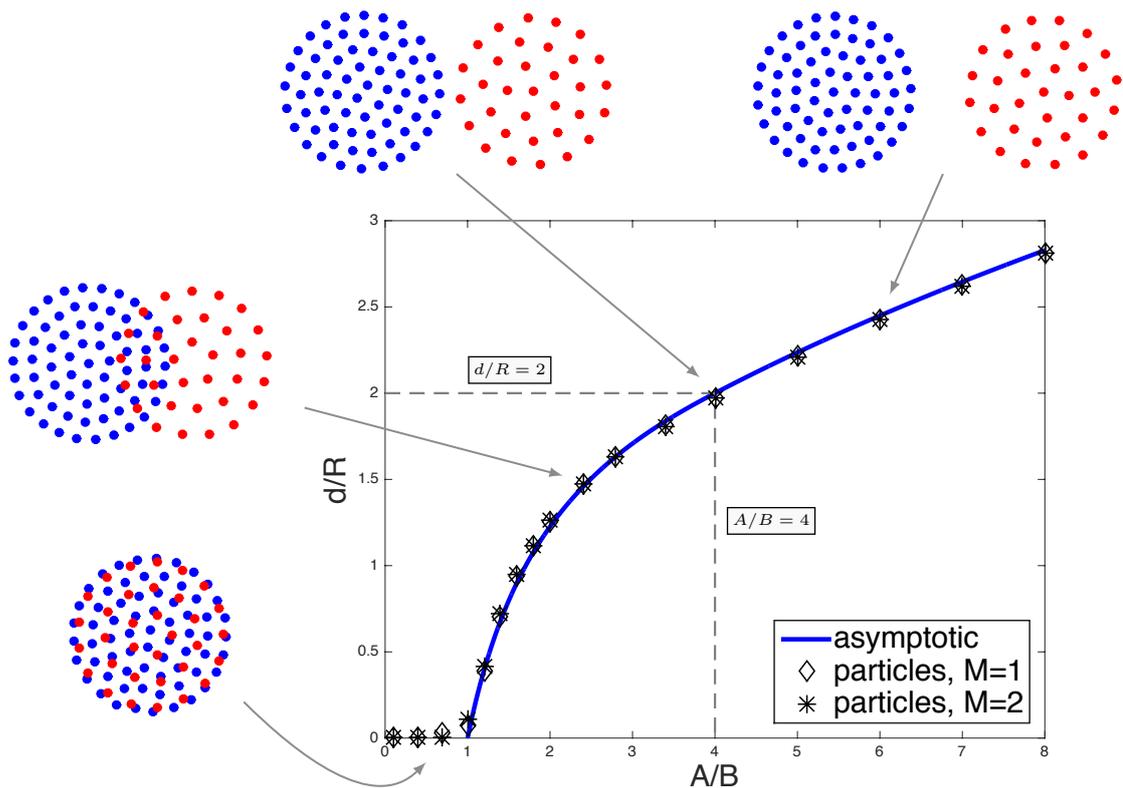
\begin{figure}[h]
\centering
\begin{tikzpicture}[>= latex]
\begin{scope}[xscale=2,yscale=2]
\draw (-1.6,2.6) node {\includegraphics[width=.4\textwidth]{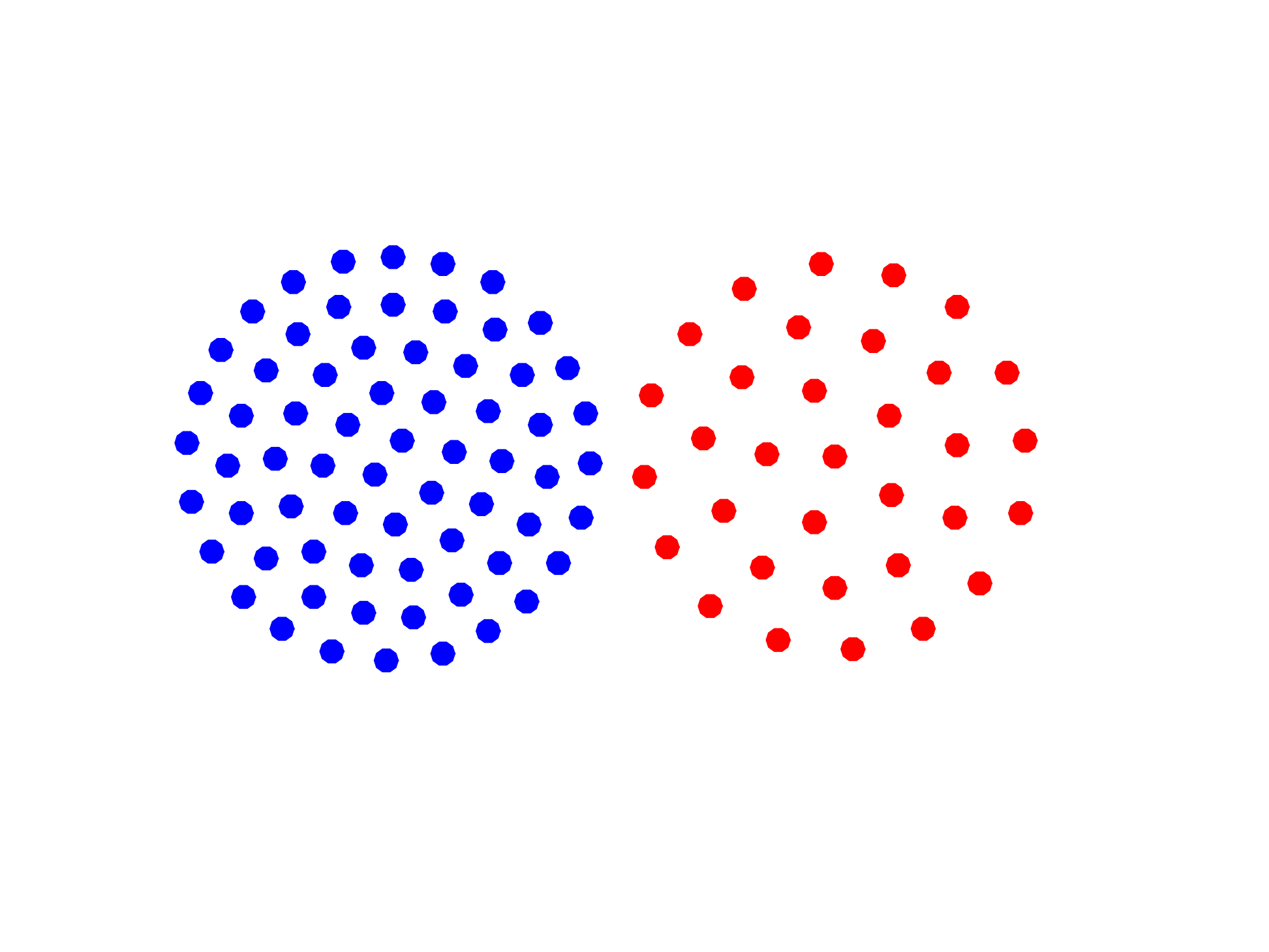}} ;
\draw (1.5,2.6) node {\includegraphics[width=.4\textwidth]{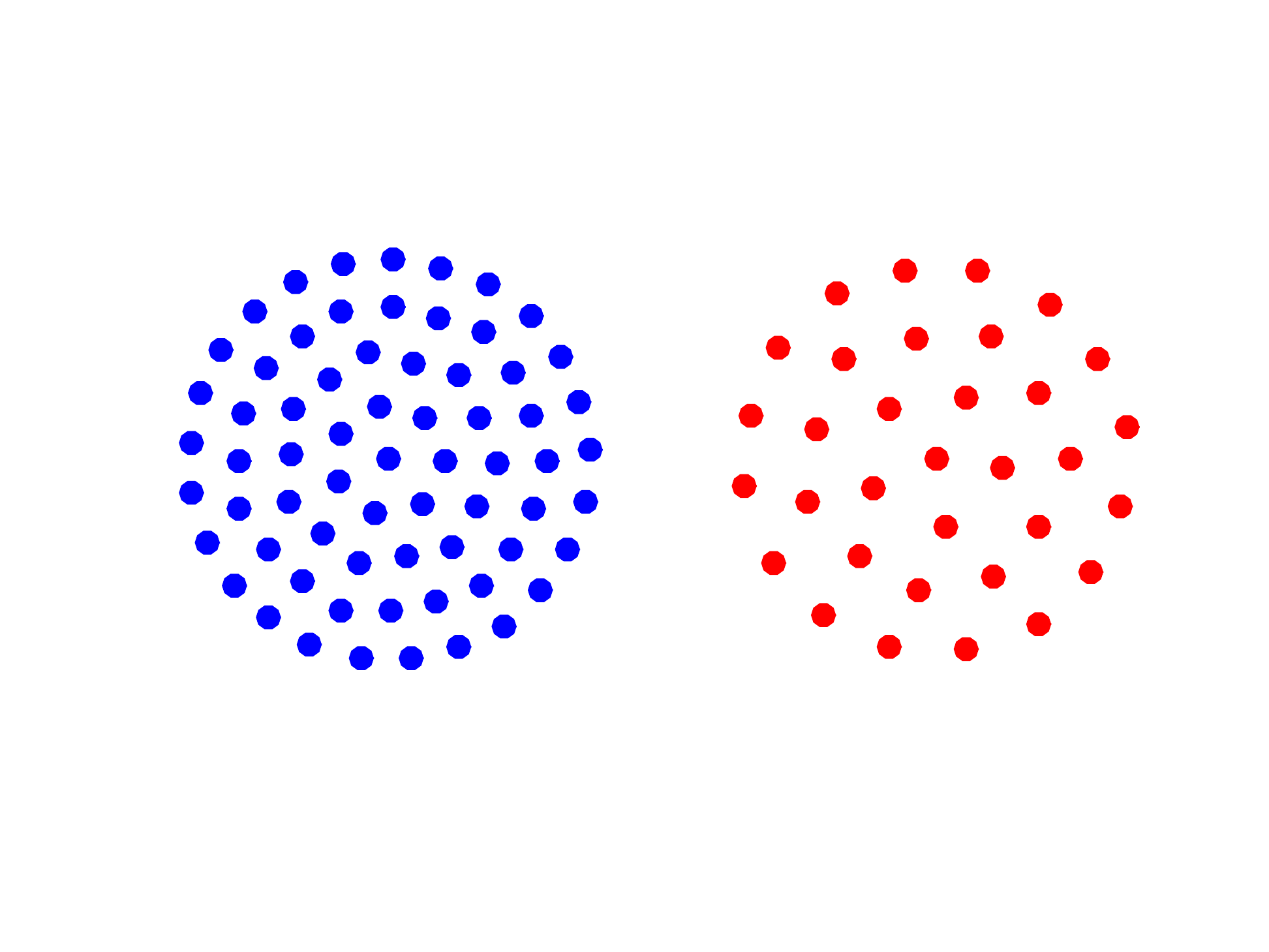}};
\draw  (-3.4,0.8) node {\includegraphics[width=.4\textwidth]{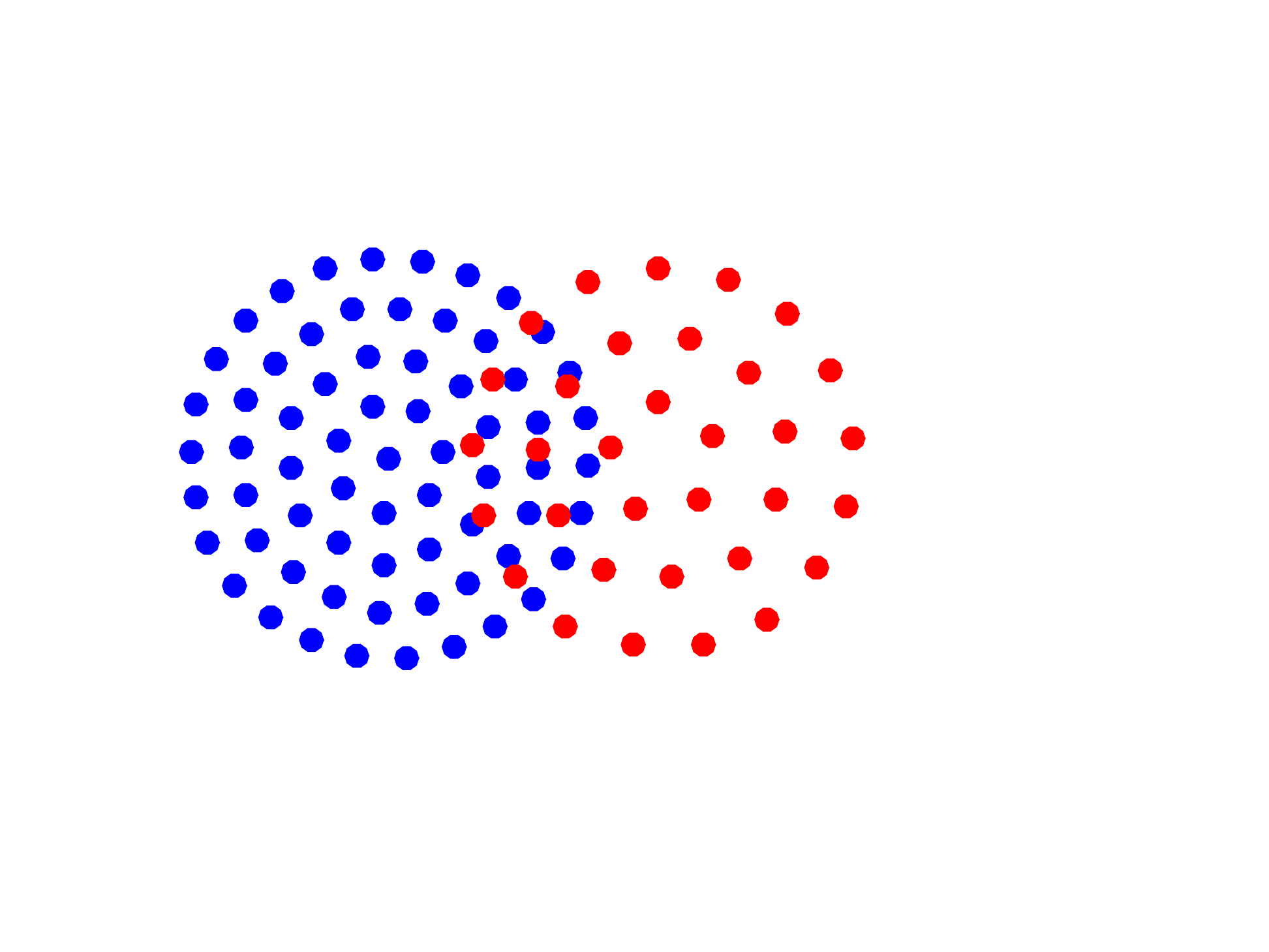}};
\draw (-3,-1) node {\includegraphics[width=.4\textwidth]{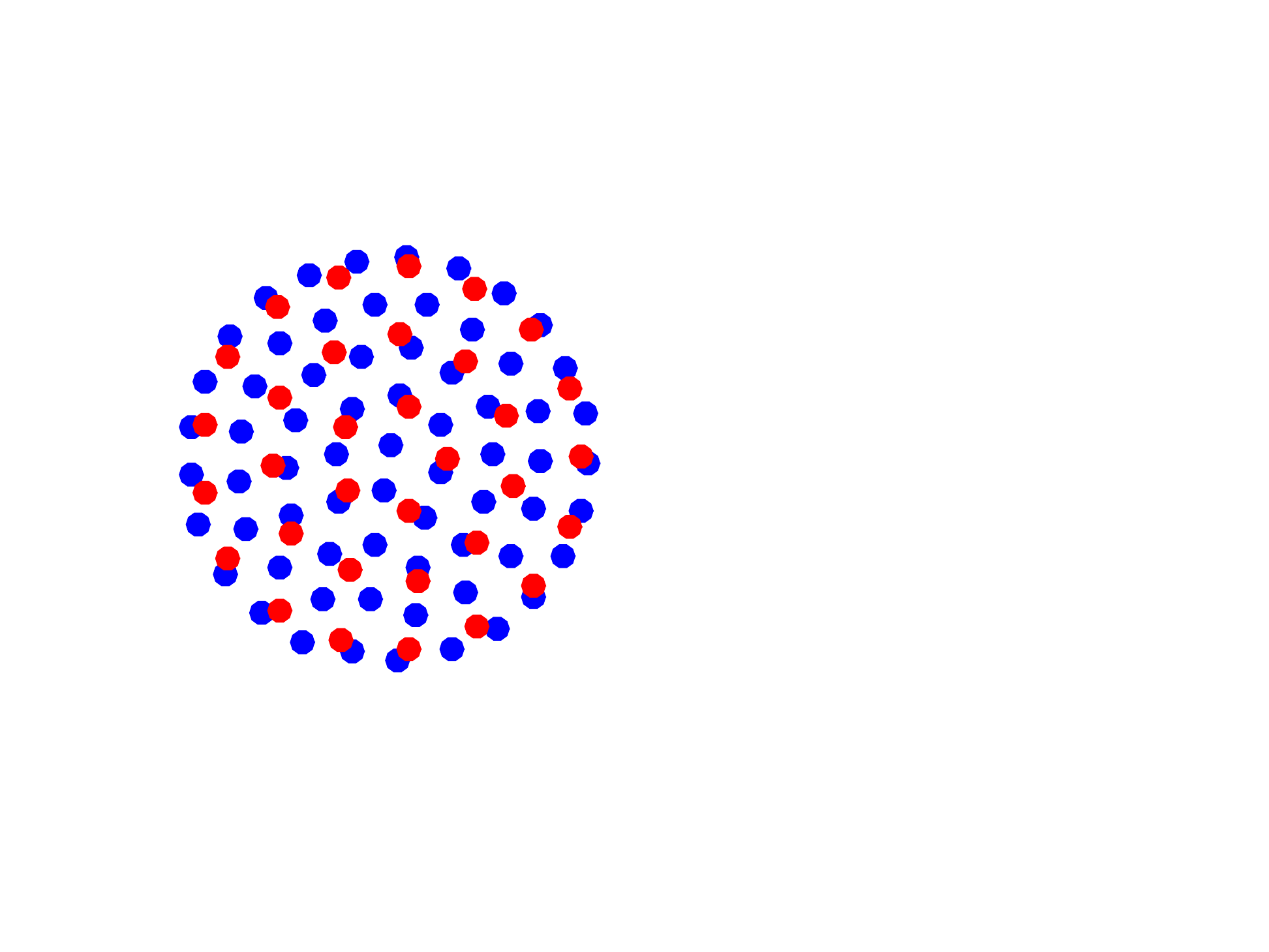}};
\draw (0,0) node {\includegraphics[width=0.7\textwidth]{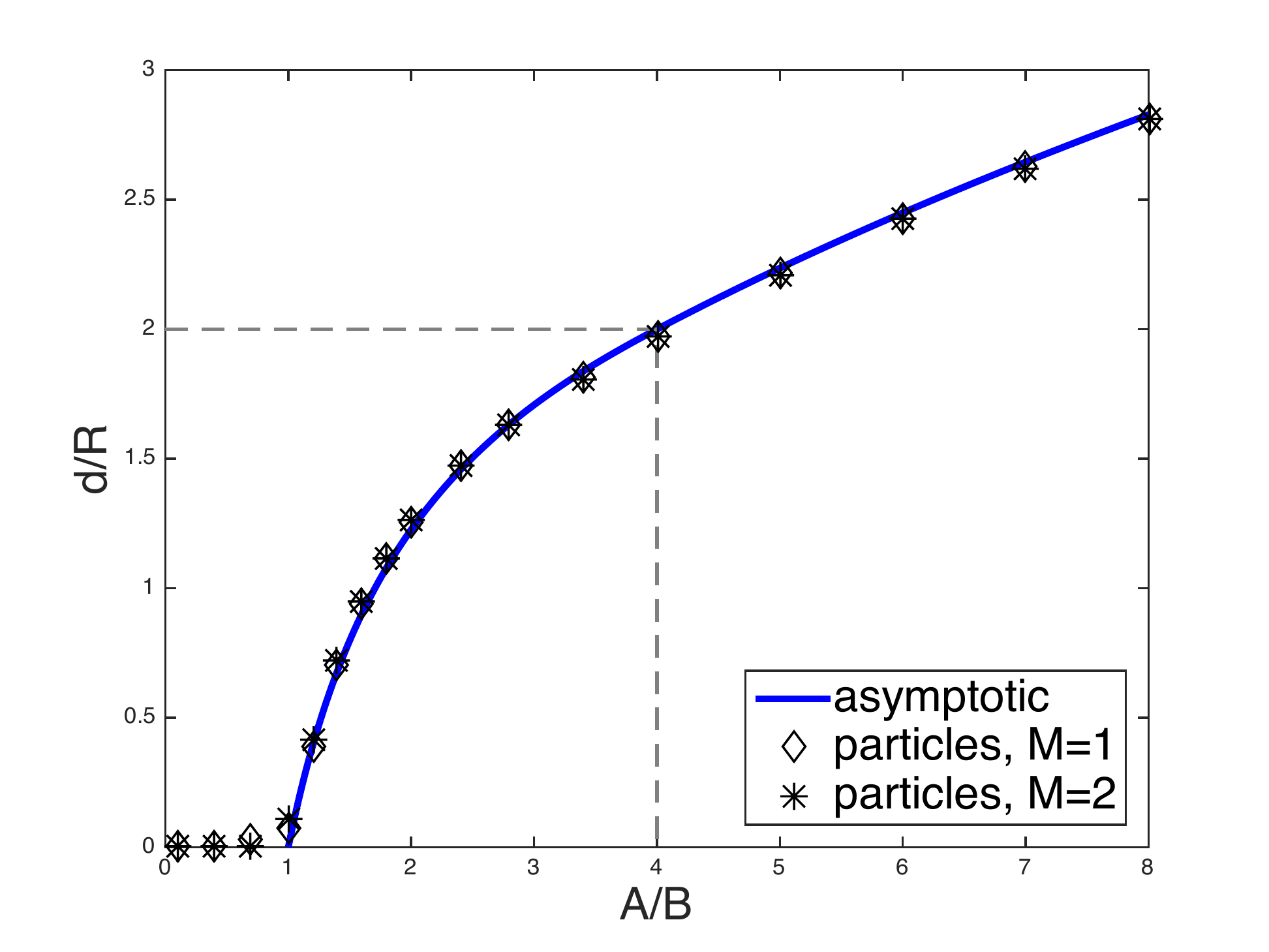}};

\draw [line width=0.75,->,color=gray!90] (1.6,2.1)  -- (1.25,1.3);
\draw [line width=0.75,->,color=gray!90] (-1.6,2.1) -- (0,0.75);
\draw [line width=0.75,->,color=gray!90] (-2.6,0.55) -- (-0.9,0.1);
\draw [line width=0.75,->,color=gray!90] (-3,-1.4)  to[out=-45,in=225] (-1.8,-1.8);

\begin{scope} 
\tikzstyle{every node} = [draw,rectangle, fill=gray!5, line width = 0.5,inner sep=2pt]
\draw (-1.25,0.8) node {\begin{tiny}$d/R=2$\end{tiny}};
\draw (0.47,-0.2) node {\begin{tiny}$A/B=4$\end{tiny}};
\end{scope}
\end{scope}
\end{tikzpicture}
\caption{Blue curve: $d/R$ as a function of $A/B$ provided by \eqref{eqn:rel d ab for d large}, and the implicit relations \eqref{eqn:impl eqn d A/B for d intermediate} and \eqref{eqn:impl eqn d A/B for d small}. Data points: estimate of $d/R$ based on particle simulations at time $t=3000$, for $\eta=0.05$ and several values of $A/B$. Black diamonds: $M=1$, hence 100 particles per species. Black stars: similar calculations for $M=2$, i.e.~133 particles of species 1, and 67 particles of species 2.
Typical particle configurations are given to illustrate full mixing, partial overlap, tangential disks (at the point where $A/B=4$ and $d/R=2$), and full separation.}
\label{fig:num}%
\end{figure}

We verify the relation between $A/B$ and $d/R$ provided by \eqref{eqn:rel d ab for d large}, \eqref{eqn:impl eqn d A/B for d intermediate} and \eqref{eqn:impl eqn d A/B for d small}. In Figure \ref{fig:num}, the blue curve is composed of three segments corresponding to the derived expressions for $d>2R$, $R<d\leqs2R$ and $d\leqs R$, respectively.
The black diamonds and stars are based on the evolution of a particle system of 200 particles. Initially they are distributed randomly. An estimate for $d$ is obtained from the long-time configuration (at $t=3000$), in which we compute the distance between the centres of mass of both species. We took $\eta=0.05$, $a_s=b_s=b_c=1$ and varied the value of $a_c$. Note that hence, $R=1$.

\begin{figure}[h]
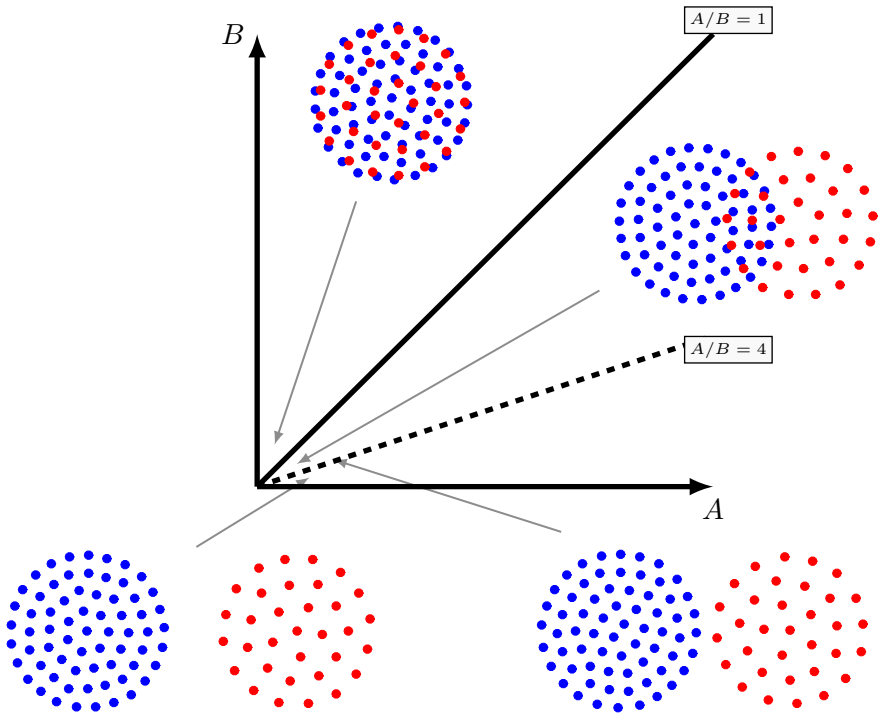

\centering
\begin{tikzpicture}[>= latex]
\begin{scope}[xscale=2,yscale=2]
\pgfmathsetmacro\M {2};
\draw (1.5,2.5) node {\includegraphics[width=.4\textwidth]{A05B1_t3000.pdf}};
\draw  (3.5,1.7) node {\includegraphics[width=.4\textwidth]{A2B1_t3000.pdf}};
\draw (3,-1) node {\includegraphics[width=.4\textwidth]{A4B1_t3000.pdf}} ;
\draw (-0.5,-1) node {\includegraphics[width=.4\textwidth]{A6B1_t3000.pdf}};
\draw [line width=0.75,->,color=gray!90] (30:2.6)  -- (30:0.3);
\draw [line width=0.75,->,color=gray!90] (71:2) -- (67.5:0.3);
\draw [line width=0.75,<-,color=gray!90] (0.5,0.18) -- (2,-0.3);
\draw [line width=0.75,->,color=gray!90] (-0.4,-0.4)  -- (10:0.35);

\draw[->,line width=2] (0,0) -- (0,3) node[left] {$B$};
\draw[->,line width=2] (0,0) -- (3,0) node[below] {$A$};
\draw[line width=2] (0,0) -- (3,3);
\draw[line width=2,dashed] (0,0) -- (3,1);

\begin{scope} 
\tikzstyle{every node} = [draw,rectangle, fill=gray!5, line width = 0.5,inner sep=2pt]
\draw (3.1,3) node[above] {\begin{tiny}$A/B=1$\end{tiny}};
\draw (3.1,1) node[below] {\begin{tiny}$A/B=4$\end{tiny}};
\end{scope}
\end{scope}
\end{tikzpicture}
\vspace{-1.5cm}
\caption{Steady states for weak cross-interactions depend on the value of $A/B$. The figure should be considered as the phase plane Figure \ref{fig:stability} zoomed in close to the origin. The lines $A/B=1$ and $A/B=4$ represent the theoretically derived asymptotic thresholds for full mixing and complete separation. We observe full mixing for $A/B<1$, partial mixing for $1<A/B<4$, two tangential circles for $A/B=4$, and separation for $A/B>4$.}
\label{fig:phase plane small cross}%
 \end{figure}

To show that the relations between $d/R$ and $A/B$ are independent of the mass ratio $M=M_1/M_2$, we perform the numerical calculations for $M=1$ (diamonds) and $M=2$ (stars). Both cases are nearly identical and coincide with the blue curve, hence confirm our prediction based on the asymptotic analysis. The blue curve also shows that $d=0$ is attained at $A/B=1$, the point at which a pitchfork bifurcation takes place. For $A/B<1$, the particle system calculations exhibit full mixing, i.e.~$d\approx 0$.
Figure \ref{fig:num} contains typical examples of particle configurations for full mixing, partial overlap, tangential disks, and full separation. We will discuss these regimes now once more, using the phase plane $A$ versus $B$. The weak cross-interactions regime $\eta\ll 1$ corresponds to an area infinitesimally close to the origin in the $(A,B)$ plane. In Figure \ref{fig:phase plane small cross} we `zoom in' near the origin and indicate which steady states we can expect to occur where. We emphasize that our considerations only hold asymptotically as $\eta\downarrow 0$.

For $A/B<1$, we concluded that total overlap is to be expected. This happens in the area above the line $B=A$ in Figure \ref{fig:phase plane small cross}. In the top configuration the two densities are supported on the same disk of radius $R$. Their densities may differ, depending on $M$. In the figure we have $M=2$ and the density of species 1 is therefore larger than the density of species 2.

For $A/B>1$ there is a bifurcation and steady states other than complete overlap come into existence. In Figure \ref{fig:num} we show the (scaled) distance between the two centres of mass being larger than zero in this case. In the phase plane (see Figure \ref{fig:phase plane small cross}, below the line $B=A$) we consequently see a non-radially symmetric state in which the two species are each supported on a disk, but the centres of the disk do not coincide. See the top right configuration. There is still a region of overlap, though, as long as $A/B<4$.

The threshold value $A/B=4$ denotes the transition from partial overlap to full separation. For $A/B$ we observe two tangential circular states (bottom right configuration in Figure \ref{fig:phase plane small cross}). For $A/B>4$ (that is: below the line $B=A/4$ in Figure \ref{fig:phase plane small cross}) the two disks are fully separated; see the bottom left configuration. The distance between the centres of mass is predicted by the relation \eqref{eqn:rel d ab for d large}, i.e.~$d/R=\sqrt{A/B}$.

In Figure \ref{fig:phase plane small cross}, the lines $A/B=1$ and $A/B=4$ are there only for their slopes to indicate the \emph{asymptotic} threshold values. If $(A,B)$ is taken $\Ord(1)$ away from the origin, it remains to be investigated if phase plane boundaries between geometrically different steady states can be found, and whether or not these are straight lines.

%
%
%
\section{Discussion}
In this paper, we produced a catalogue of steady states for model \eqref{eqn:model} with interaction kernels \eqref{eqn:intpot}, as presented in Figure \ref{fig:num st st}.

We argued (see Section \ref{sect:target}) by means of linear stability analysis that the target with the lighter species inside exists and is stable for parameters $(A,B)$ in regions $D_4$ and $D_5$. The target with the heavier species inside is (if it exists at all) not stable with respect to small perturbations of the boundaries. Both results agree with what we observe numerically. Our variational approach in both cases tends to predict stability regions that are larger than those for the linear stability analysis. We expect this to be a result of the fact that in some parameter regions, due to the difficulties explained in Section \ref{subsect:prelim-var}, we have not considered a certain type of perturbations, namely those supported in the support of the equilibrium (referred in the paper as class $\CalA$ perturbations). Nevertheless, in the region of the parameter space where $A<1$ and $B>1$ we were able to do a complete variational argument to show that the the target with the lighter species inside is a global minimizer.

The `overlap' state with the lighter species inside is numerically observed only in $D_6$. In Section \ref{sect:overlap} we verified this conjecture using our variational approach and found that it is a local minimizer of the energy with respect to class $\CalB$ perturbations. Moreover, in the subset of $D_6$ with $A<1$ and $B>1$ (which is all of $D_6$ except a bounded triangular region), a full variational analysis showed that the overlap equilibrium with the lighter species inside is a global minimizer of the energy. The overlap solution with the \emph{heavier} species inside is never observed numerically. By the variational approach however, we found that this overlap state is a local minimizer with respect to class $\CalB$ perturbations, if and only if parameters are taken from region $D_1$. Based on numerical results, it is our conjecture that in $D_1$ this state is \emph{not} a minimizer with respect to class $\CalA$ perturbations. 

The various non-radially symmetric steady states in Figure \ref{fig:num st st} were investigated numerically in more depth. This was shown in Figures \ref{fig:unstab modes targStab}, \ref{fig:unstab modes targUnstab} and \ref{fig:unstab modes overlap}. These states arose from initializing the system \eqref{eqn:part syst} in a radially symmetric configuration and taking parameters outside the corresponding stability region. We were able to observe the specific modes of instability.

Finally we shed light on the symmetry-breaking that happens when passing parameters from $D_6$ to $D_1$. In Section \ref{sec:weak cross} we examined the limit of weak cross interactions and obtained an (asymptotically valid) relation between the parameters $(A,B)$ and the distance between the centres of mass of the two species.

The choice of kernels \eqref{eqn:intpot} in this paper is quite specific. Our main motivation for taking these Newtonian-quadratic kernels is the fact that they leave room for obtaining results analytically. For instance, the fact that all steady states in this paper are ``piecewise" constant, is a direct result of the choice of kernels. Our method of performing linear stability analysis was inspired by this observation. However, from the point of view of the biological applications, one might prefer alternative kernels that induce interactions with a limited range. We acknowledge that kernels \eqref{eqn:intpot} even lead to \emph{increasing} attraction as the distance grows. One should be aware however that different kernels will lead inevitably to changes in the nature of the steady states.

A less radical way to increase realism is to remove the symmetry in the interactions between species 1 and 2. There is no direct biological reason why the two respective species among themselves would behave according to the same parameters. Neither is there a reason why species 1 would respond in the same way to species 2 as species 2 responds to species 1. Instead of having two repulsion parameters $\as$ and $\ac$, one could therefore introduce four parameters: $a_{11}$, $a_{22}$ for self-repulsion and $a_{12}$, $a_{21}$ for cross-repulsion. The analogue can be done for the attraction parameters $\bs$ and $\bc$. We note that such modifications may alter the dynamics of the system (e.g. introducing chasing dynamics). Moreover, the resulting system will in general no longer possess the gradient flow structure. Consequently, the variational approach of this paper may be not applicable anymore, and other methods for investigation of stability need to be designed.

Furthermore, making the step from a one-species model to a two-species model naturally opens the door to exploring multi-species models. Some preliminary numerical experiments for three species (not presented in the current paper) indicate that many interesting patterns may be expected.

\begin{appendix}
\section{Calculation of the basic integrals for perturbed boundaries}\label{app:pert integrals}
For the linear stability analysis, we evaluate integrals \eqref{eqn:pert int} on the perturbed boundaries; that is, for $x=p_j(\theta_0)$ for some index $j$. The latter integral in \eqref{eqn:pert int} becomes
\begin{align}
\int_{\Omega_\ell^\eps} (x-y)\,dy =&\; x\,\int_{\Omega_\ell^\eps}\,dy \; - \int_{\Omega_\ell^\eps}y\,dy
=\, \begin{cases}
    x\,\pi R_\ell^2 - \pi R_\ell^3\,\eps_{\ell,N}     & \quad \text{if } m=1,\\
    x\,\pi R_\ell^2  & \quad \text{if } m\geqs 2;\\
  \end{cases} 
  \label{eqn:eval attr int}
\end{align}
where second- and higher-order terms in $\underline{\eps}$ are omitted. Here we used the full expressions for the area and centre of mass of $\Omega_\ell^\eps$ as given in Section \ref{sect:prelim lin}. Note that by taking $x=p_j(\theta_0)$, we introduce extra $\Ord(\underline{\eps})$ terms in \eqref{eqn:eval attr int}.

To find the first-order approximation of \eqref{eqn:int ln Gauss}, we take $x=p_j(\theta_0)$, parameterize $\partial\Omega_\ell^\eps$ by $p_\ell(\theta)$ and compute that
\begin{equation}\label{eqn:n dS}
\hat{n}\,dS = R_\ell\, e^{i\theta}(1+(\eps_{\ell,N}+m\eps_{\ell,T})\cos(m\theta)+i\,(\eps_{\ell,T}+m\eps_{\ell,N})\sin(m\theta))\,d\theta.
\end{equation}
Subsequently, we obtain
\begin{multline}
-\int_{\partial\Omega_\ell^\eps} \ln|x-y|\,\hat{n}\,dS =\\
-\int_0^{2\pi} \ln(|p_j(\theta_0)-p_\ell(\theta)|)\, R_\ell\, e^{i\theta}(1+(\eps_{\ell,N}+m\eps_{\ell,T})\cos(m\theta)+i\,(\eps_{\ell,T}+m\eps_{\ell,N})\sin(m\theta))\,d\theta.\label{eqn:ln int}
\end{multline}
The integrand can be written as
\[\ln|p_j(\theta_0)-p_\ell(\theta)|=\frac12\ln|p_j(\theta_0)-p_\ell(\theta)|^2\]
and we expand it in terms of $\underline{\eps}$. First, we expand:
\begin{align}
\nonumber |p_j(\theta_0)-p_\ell(\theta)|^2 =&\, \overbrace{|R_j\,e^{i\theta_0}-R_\ell\,e^{i\theta}|^2}^{=R_j^2+R_\ell^2-2R_jR_\ell \cos(\theta-\theta_0)} \\
\nonumber &+ \eps_{j,N}2R_j\,\cos(m\theta_0) \cdot(R_j-R_\ell\cos(\theta-\theta_0)) \\
\nonumber &+ \eps_{\ell,N}2R_\ell\,\cos(m\theta)\cdot(R_\ell-R_j\cos(\theta-\theta_0)) \\
&- 2R_j R_\ell\sin(\theta-\theta_0)\cdot(\eps_{j,T}\sin(m\theta_0)-\eps_{\ell,T}\sin(m\theta))     + \Ord(|\underline{\eps}|^2).\label{eqn:expand arg ln}
\end{align}
Here, we use the generic notation $\Ord(|\underline{\eps}|^2)$ for higher-order terms in any $\eps_{j,N}$, $\eps_{j,T}$, $\eps_{\ell,N}$ or $\eps_{\ell,T}$. Introduce the notation $\alpha:=R_\ell/R_j$. We use \eqref{eqn:expand arg ln} and $\ln(X+\eps Y)\sim\ln(X)+\eps Y/X$ in \eqref{eqn:ln int}, we omit further $\Ord(|\underline{\eps}|^2)$ terms and write sines and cosines as complex exponentials. After expanding \eqref{eqn:ln int} in this way, there is a part containing a logarithm that consists of integrals of the form
\begin{equation}\label{eqn:log integral}
\int_{0}^{2\pi} \ln(1+\alpha^2-2\alpha \cos(\theta-\theta_0))\,e^{i\mu\theta}\,d\theta = \, -\dfrac{2\pi}{|\mu|}\, e^{i\mu\theta_0} \begin{cases}
    \alpha^{|\mu|}    & \quad \text{if } \alpha\leqs1,\\
    \alpha^{-|\mu|}  & \quad \text{if } \alpha>1;
  \end{cases} 
\end{equation}
for several values of $\mu\in\Z\setminus\{0\}$. The expression in \eqref{eqn:log integral} is derived in Appendix \ref{app:basic int} and is valid for any $\alpha>0$, including $\alpha=1$. Note that $\ln(R_j^2+R_\ell^2-2R_jR_\ell \cos(\theta-\theta_0))=\ln(R_j^2)+\ln(1+\alpha^2-2\alpha \cos(\theta-\theta_0))$, while one can show that the part of \eqref{eqn:ln int} containing $\ln(R_j^2)$ has zero contribution eventually. Consequently, the ``logarithmic part" of \eqref{eqn:ln int}, equals
\begin{align}
%
\begin{cases}
    \pi R_\ell e^{i\theta_0}\left(\beta +\dfrac12(\eps_{\ell,N}+\eps_{\ell,T})\beta^{2} e^{i\theta_0}\right)  & \quad \text{if } m=1,\\\\
       \pi R_\ell e^{i\theta_0}\left(\beta +\dfrac12(\eps_{\ell,N}+\eps_{\ell,T})\beta^{m+1} e^{im\theta_0}-\dfrac12(\eps_{\ell,N}-\eps_{\ell,T})\beta^{m-1}e^{-im\theta_0}\right) & \quad \text{if } m\neq1;
  \end{cases} \label{eqn:int log}
\end{align}
with $\beta:=\min\{R_j,R_\ell\}/\max\{R_j,R_\ell\}\leqs 1$.

There is a ``rational part" in \eqref{eqn:ln int} that, for $\alpha\neq1$, consists of contributions of the form
\begin{equation}
\int_{0}^{2\pi} \dfrac{e^{i\mu\theta}}{1+\alpha^2-2\alpha \cos(\theta-\theta_0)}\,d\theta 
=2\pi e^{i\mu\theta_0}\begin{cases}
       \dfrac{\alpha^{|\mu|}}{1-\alpha^2} & \quad \text{if } \alpha<1,\\
    \dfrac{\alpha^{-|\mu|}}{\alpha^2-1}  & \quad \text{if } \alpha>1.
  \end{cases} 
\label{eqn:rational integral alpha not 1}
\end{equation}
This expression is also derived in Appendix \ref{app:basic int} and it is valid for all $\mu\in\Z$.

The ``rational part" of \eqref{eqn:ln int} becomes (for $\alpha\neq1$):
\begin{equation}\footnotesize
\begin{cases}
       -\dfrac{R_\ell}{2}\pi\alpha e^{i\theta_0}\Big[\Big( \eps_{j,N}-\eps_{j,T}-\alpha(\eps_{\ell,N}-\eps_{\ell,T})\Big)e^{i\theta_0} +\left( \eps_{j,N}+\eps_{j,T}\right)e^{-i\theta_0}  \Big] & \quad \text{if } m=1 \text{ and } \alpha<1,\\\\
           -\dfrac{R_\ell}{2}\pi\alpha e^{i\theta_0}\Big[\Big( \eps_{j,N}-\eps_{j,T}-\alpha^m(\eps_{\ell,N}-\eps_{\ell,T})\Big)e^{im\theta_0}  &\\
           \hspace{0.25\linewidth}+\Big( \eps_{j,N}+\eps_{j,T}-\alpha^{m-2}(\eps_{\ell,N}-\eps_{\ell,T})\Big)e^{-im\theta_0}  \Big]  & \quad \text{if } m\neq1 \text{ and } \alpha<1,\\\\
           \dfrac{R_\ell}{2}\dfrac{\pi}{\alpha} e^{i\theta_0}\Big[\Big( \eps_{j,N}+\eps_{j,T}-\dfrac{1}{\alpha}(\eps_{\ell,N}+\eps_{\ell,T})\Big)e^{i\theta_0} +\Big( \eps_{j,N}-\eps_{j,T}-2\alpha\eps_{\ell,N}\Big)e^{-i\theta_0}  \Big]  & \quad \text{if } m=1 \text{ and } \alpha>1,\\\\
    \dfrac{R_\ell}{2}\dfrac{\pi}{\alpha} e^{i\theta_0}\Big[\Big( \eps_{j,N}+\eps_{j,T}-\alpha^{-m}(\eps_{\ell,N}+\eps_{\ell,T})\Big)e^{im\theta_0} &\\
    \hspace{0.25\linewidth}+\Big( \eps_{j,N}-\eps_{j,T}-\alpha^{2-m}(\eps_{\ell,N}+\eps_{\ell,T})\Big)e^{-im\theta_0}  \Big]   & \quad \text{if } m\neq1 \text{ and } \alpha>1.
  \end{cases} \label{eqn:int rat alpha neq 1}
\end{equation}

This approach is correct as long as $j\neq\ell$, since otherwise $R_j= R_\ell$, thus $\alpha=1$, and there is a singularity in the denominator. This case requires a partially different approach. If $j=\ell$, then the rational $\Ord(\underline{\eps})$ part becomes
\begin{multline}
-\dfrac{R_j}{2} \int_{0}^{2\pi} \left[\eps_{j,N}\cos(m\theta_0)+\eps_{j,N}\cos(m\theta)+ \eps_{j,T}\Big(\sin(m\theta)-\sin(m\theta_0)\Big)\dfrac{\sin(\theta-\theta_0)}{1-\cos(\theta-\theta_0)}\right] e^{i\theta}\,d\theta\\
= \begin{cases}
    -\dfrac{R_j}{2}\pi(\eps_{j,N} + \eps_{j,T})   & \quad \text{if } m=1,\\\\       
       -R_j\pi\,\eps_{j,T}\,e^{-i(m-1)\theta_0} & \quad \text{if } m\neq1.
  \end{cases} \label{eqn:int rat alpha 1}
\end{multline}
In the latter part of the integrand the singularity in the denominator is compensated by the terms in the numerator, which can be seen by expanding all sines and cosines in complex exponentials.
Note that if we take $j=\ell$, we see that the limits $\alpha\downarrow1$ and $\alpha\uparrow1$ in \eqref{eqn:int rat alpha neq 1} agree and are equal to the expressions in \eqref{eqn:int rat alpha 1} for both $m=1$ and $m\neq1$.

In conclusion, combining \eqref{eqn:int log}, \eqref{eqn:int rat alpha neq 1} and \eqref{eqn:int rat alpha 1}, we have the following:
\paragraph*{Case $R_\ell<R_j$, hence $\alpha<1$ and $\beta=\alpha=R_\ell/R_j$:} For all $m\geqs1$,
\begin{multline}
\int_{\Omega_\ell^\eps} \dfrac{x-y}{|x-y|^2}\,dy = \pi R_\ell e^{i\theta_0}\left[\beta -\left(\dfrac{1}{2}\beta(\eps_{j,N}-\eps_{j,T})-\beta^{m+1}\,\eps_{\ell,N}\right)e^{im\theta_0}-\dfrac{1}{2}\beta(\eps_{j,N}+\eps_{j,T})e^{-im\theta_0}\right]\\
= \pi R_\ell e^{i\theta_0}\left[ \beta + \left(-\beta\eps_{j,N}+\beta^{m+1}\eps_{\ell,N}\right)\,\cos(m\theta_0) +\left(\beta\eps_{j,T}+\beta^{m+1}\eps_{\ell,N}\right)\,i\sin(m\theta_0)\right].\label{eqn:pert int rep outside}
\end{multline}

\paragraph*{Case $\ell=j$, hence $R_\ell=R_j$ and $\beta=\alpha=1$:} For all $m\geqs1$,
\begin{align}
\nonumber \int_{\Omega_\ell^\eps} \dfrac{x-y}{|x-y|^2}\,dy =&\, \pi R_j e^{i\theta_0}\left[1 +\dfrac12(\eps_{j,N}+\eps_{j,T}) e^{im\theta_0}-\dfrac12(\eps_{j,N}+\eps_{j,T}) e^{-im\theta_0}\right]\\
=&\, \pi R_j e^{i\theta_0}\left[ 1 + \left(\eps_{j,N}+\eps_{j,T}\right)\,i\sin(m\theta_0)\right]. \label{eqn:pert int rep same}
\end{align}

\paragraph*{Case $R_\ell>R_j$, hence $\alpha>1$ and $\beta=1/\alpha=R_j/R_\ell$:} For all $m\geqs1$, 
\begin{multline}
\int_{\Omega_\ell^\eps} \dfrac{x-y}{|x-y|^2}\,dy = \pi R_\ell e^{i\theta_0}\left[\beta +\dfrac12 \beta(\eps_{j,N}+\eps_{j,T}) e^{im\theta_0} + \left( \dfrac12\beta(\eps_{j,N}-\eps_{j,T})-\beta^{m-1}\eps_{\ell,N} \right)  e^{-im\theta_0} \right]\\
= \pi R_\ell e^{i\theta_0}\left[ \beta + \left(\beta\eps_{j,N}-\beta^{m-1}\eps_{\ell,N}\right)\,\cos(m\theta_0) +\left(\beta\eps_{j,T}+\beta^{m-1}\eps_{\ell,N}\right)\,i\sin(m\theta_0)\right]. \label{eqn:pert int rep inside}
\end{multline}

Note that the above expressions are consistent in the sense that if we set $j=\ell$ (hence, $\beta=1$) in either \eqref{eqn:pert int rep outside} or \eqref{eqn:pert int rep inside} we obtain \eqref{eqn:pert int rep same}. Note also that if we set $\eps_{j,N}=\eps_{j,T}=\eps_{\ell,N}=\eps_{\ell,T}=0$ in \eqref{eqn:pert int rep outside}-\eqref{eqn:pert int rep same}-\eqref{eqn:pert int rep inside}, the expressions are consistent with \eqref{eqn:intderiv}.

\section{Evaluation of integrals \eqref{eqn:log integral} and  \eqref{eqn:rational integral alpha not 1}}\label{app:basic int}
The integral on the left-hand side of \eqref{eqn:rational integral alpha not 1} is treated as follows:
\begin{equation*}
\int_{0}^{2\pi} \dfrac{e^{i\mu\theta}}{1+\alpha^2-2\alpha \cos(\theta-\theta_0)}\,d\theta 
= e^{i\mu\theta_0}\int_{-\pi}^{\pi} \dfrac{e^{i\mu\phi}}{1+\alpha^2-2\alpha \cos(\phi)}\,d\phi,
\end{equation*}
where the substitution $\phi=\theta-\theta_0$ is used. The real part of the integrand on the right-hand side is an even function in $\phi$, while the imaginary part is odd. Hence,
\begin{equation*}
\int_{-\pi}^{\pi} \dfrac{e^{i\mu\phi}}{1+\alpha^2-2\alpha \cos(\phi)}\,d\phi = 2\int_{0}^{\pi} \dfrac{\cos(\mu\phi)}{1+\alpha^2-2\alpha \cos(\phi)}\,d\phi.
\end{equation*}
The latter integral is given in \cite[p.~253, No.~31]{Jeffrey}. Note that \cite{Jeffrey} only treats $\mu\geqs0$, but that the expression given therein is easily generalized to all $\mu\in\Z$. The above considerations combined yield \eqref{eqn:rational integral alpha not 1}.

To find \eqref{eqn:log integral}, we compute the derivative of the left-hand side with respect to $\theta_0$ in two different ways. On the one hand:
\begin{align}
\nn \dfrac{d}{d\theta_0}\int_{0}^{2\pi} \ln(1+\alpha^2-2\alpha \cos(\theta-\theta_0))\,e^{i\mu\theta}\,d\theta =&\, \int_{0}^{2\pi} \dfrac{-2\alpha\sin(\theta-\theta_0)}{1+\alpha^2-2\alpha \cos(\theta-\theta_0)}\,e^{i\mu\theta}\,d\theta\\
=&\, -2\alpha\, e^{i\mu\theta_0} \int_{0}^{2\pi} \dfrac{\sin(\phi)}{1+\alpha^2-2\alpha \cos(\phi)}\,e^{i\mu\phi}\,d\phi. \label{eqn:der log int 1}
\end{align}
Again, we used the substitution $\phi=\theta-\theta_0$. Writing $\sin\phi$ in terms of complex exponentials $\sin\phi=(\exp(i\phi)-\exp(-i\phi))/(2i)$, we can express the latter integral as the difference of two integrals of the form \eqref{eqn:rational integral alpha not 1}. By working out the result for all $\mu\in\Z\setminus\{0\}$, we obtain
\begin{equation}\label{eqn:der log int 1 eval}
-2\alpha\,e^{i\mu\theta_0} \int_{0}^{2\pi} \dfrac{\sin(\phi)}{1+\alpha^2-2\alpha \cos(\phi)}\,e^{i\mu\phi}\,d\phi = \,  \begin{cases}
    -2\pi i \, e^{i\mu\theta_0} \sgn(\mu)\, \alpha^{|\mu|}    & \quad \text{if } \alpha<1,\\
    -2\pi i \, e^{i\mu\theta_0} \sgn(\mu)\, \alpha^{-|\mu|}  & \quad \text{if } \alpha>1;
  \end{cases} 
\end{equation}
On the other hand:
\begin{align}
\nn \dfrac{d}{d\theta_0}\int_{0}^{2\pi} \ln(1+\alpha^2-2\alpha \cos(\theta-\theta_0))\,e^{i\mu\theta}\,d\theta          
\nn =&\, i\,\mu\, e^{i\mu\theta_0} \,\int_{0}^{2\pi} \ln(1+\alpha^2-2\alpha \cos(\phi))\,e^{i\mu\phi}\,d\phi\\
=&\, i\,\mu \int_{0}^{2\pi} \ln(1+\alpha^2-2\alpha \cos(\theta-\theta_0))\,e^{i\mu\theta}\,d\theta \label{eqn:der log int 2}
\end{align}
Together, \eqref{eqn:der log int 1}, \eqref{eqn:der log int 1 eval} and \eqref{eqn:der log int 2} yield the result of \eqref{eqn:log integral} in all cases except for $\alpha=1$.

For $\alpha=1$ (and $\mu\neq0$), note that the imaginary part of the integrand $\ln(2-2\cos(\phi))\,\exp(i\mu\phi)$ is an odd function, and therefore yields no contribution. Therefore, 
\begin{align}
\nn \int_{0}^{2\pi} \ln(2-2\cos(\theta-\theta_0))\,e^{i\mu\theta}\,d\theta =&\,  e^{i\mu\theta_0}\,\int_{0}^{2\pi} \ln(2-2 \cos(\phi))\,\cos(\mu\phi)\,d\phi
\end{align}
and integration by parts yields
\begin{align}
\int_{0}^{2\pi} \ln(2-2 \cos(\phi))\,\cos(\mu\phi)\,d\phi =&\, \left[ \dfrac1\mu \sin(\mu\phi) \ln\Big(2-2\cos(\phi)\Big) \right]_0^{2\pi} - \dfrac1\mu \int_0^{2\pi} \dfrac{\sin(\mu\phi)\,\sin(\phi)}{1-\cos(\phi)}d\phi.
\end{align}
The boundary terms vanish, and this can be shown by introducing series expansions around $0$ and $2\pi$ and applying l'H\^{o}pital's rule. Moreover, the integrand on the right-hand side has no singularities. To see this, expand all sines and cosines in complex exponentials (analogously to the arguments leading to \eqref{eqn:int rat alpha 1}). Some manipulations are required to observe that all apparent singularities are compensated by zeros in the numerator. Having done these manipulations, we can evaluate the integral exactly, multiplying by $\exp(i\mu\theta_0)$, we obtain the right-hand side of \eqref{eqn:log integral}, for $\alpha=1$. 

\end{appendix}

\bibliographystyle{plain}
\bibliography{lit}

\end{document}